\documentclass{amsart}
\usepackage{color}

\usepackage{amsmath} 
\usepackage{amssymb}
\usepackage{mathrsfs}
\usepackage[hidelinks]{hyperref}

\newtheorem{theorem}{Theorem}[section] 
\newtheorem{claim}[theorem]{Claim}

\newtheorem{conclusion}[theorem]{Conclusion}
\newtheorem{convention}[theorem]{Convention}
\newtheorem{observation}[theorem]{Observation}

\theoremstyle{definition}
\newtheorem{definition}[theorem]{Definition}
\newtheorem{example}[theorem]{Example}
\newtheorem{fact}[theorem]{Fact}
\newtheorem{thesis}[theorem]{Thesis}

\newtheorem{discussion}[theorem]{Discussion}

\newtheorem{hypothesis}[theorem]{Hypothesis}

\theoremstyle{remark}
\newtheorem{remark}[theorem]{Remark}
\newtheorem{question}[theorem]{Question}
\newtheorem{notation}[theorem]{Notation}
\newtheorem{context}[theorem]{Context}

\usepackage[shortlabels]{enumitem}
\setlist[enumerate,1]{label={(\Alph*)}}
\setlist[enumerate,2]{label={(\alph*)}}
\setlist[enumerate,3]{label={$\bullet_{\arabic*}$}}

\newenvironment{PROOF}[2][\proofname.]
   {\begin{proof}[#1]\renewcommand{\qedsymbol}{\textsquare$_{\rm #2}$}}
   {\end{proof}}

\newcommand{\cp}{\mathrm{c.p.}}
\newcommand{\fcp}{\mathrm{f.c.p.}}
\newcommand{\ncp}{\mathrm{n.c.p.}}

\newcommand{\lft}{\mathrm{l.f.t.}}
\newcommand{\lrp}{\mathrm{l.r.p.}}
\newcommand{\flp}{\mathrm{l.r.p.}}
\newcommand{\luft}{\mathrm{l.u.f.t.}}
\newcommand{\lup}{\mathrm{l.u.p.}}
\newcommand{\uflp}{\mathrm{u.f.l.p.}}

\newcommand{\lcr}{\mathrm{lcr}}

\newcommand{\Card}{\mathrm{Card}}
\newcommand{\chr}{\mathrm{chr}}
\newcommand{\Cr}{\mathrm{Cr}}
\newcommand{\coin}{\mathrm{co\text{-}in}}
\newcommand{\comp}{\mathrm{comp}}
\newcommand{\edge}{\mathrm{edge}}
\newcommand{\ful}{\mathrm{ful}}
\newcommand{\JEP}{\mathrm{JEP}}
\newcommand{\rfil}{\mathrm{rfil}}
\newcommand{\ruf}{\mathrm{ruf}}
\newcommand{\sep}{\mathrm{sep}}
\newcommand{\smooth}{\mathrm{smooth}}
\newcommand{\spec}{\mathrm{spec}}
\newcommand{\st}{\mathsf{st}}

\newcommand{\no}{\mathrm{no}}
\newcommand{\yes}{\mathrm{yes}}

\newcommand{\blacktrianglelefteq}{%
  \mathrel{%
    \vphantom{\blacktriangleleft}%
    \smash{\vcenter{\doSmash}}%
  }%
}
\newcommand{\doSmash}{%
  \hbox{\ooalign{%
    \noalign{\kern.5ex}
    $\blacktriangleleft$\cr
    \noalign{\kern1.1ex}
    \smash{$-$}\cr
}}}

\makeatletter
\providecommand*{\cupdot}{%
  \mathbin{%
    \mathpalette\@cupdot{}%
  }%
}
\newcommand*{\@cupdot}[2]{%
  \ooalign{%
    $\m@th#1\cup$\cr
    \hidewidth$\m@th#1\cdot$\hidewidth
  }%
}
\makeatother

\newcommand{\muult}{\mathrm{mult}}

\newcommand{\tthh}{\mathrm{th}}

\newcommand{\arity}{\mathrm{arity}}

\newcommand{\cf}{\mathrm{cf}}

\newcommand{\dom}{\mathrm{dom}}

\newcommand{\id}{\mathrm{id}}

\newcommand{\otp}{\mathrm{otp}}

\newcommand{\Rang}{\mathrm{Rang}}

\newcommand{\qf}{\mathrm{qf}}

\newcommand{\Dom}{\mathrm{Dom}}

\newcommand{\ec}{\mathrm{ec}}

\newcommand{\eq}{\mathrm{eq}}

\newcommand{\fil}{\mathrm{fil}}

\newcommand{\iif}{\mathrm{if}}
\newcommand{\lex}{\mathrm{lex}}

\newcommand{\LST}{\mathrm{LST}}

\newcommand{\Mod}{\mathrm{Mod}}

\newcommand{\Th}{\mathrm{Th}}
\newcommand{\tp}{\mathrm{tp}}

\newcommand{\uf}{\mathrm{uf}}

\newcommand{\PC}{\mathrm{PC}}

\newcommand{\Reg}{\mathrm{Reg}}

\newcommand{\bfa}{\mathbf{a}}

\newcommand{\bfb}{\mathbf{b}}

\newcommand{\bfh}{\mathbf{h}}
\newcommand{\bfI}{\mathbf{I}}
\newcommand{\bfi}{\mathbf{i}}

\newcommand{\bfj}{\mathbf{j}}

\newcommand{\bfL}{\mathbf{L}}
\newcommand{\bfM}{\mathbf{M}}

\newcommand{\bfS}{\mathbf{S}}
\newcommand{\bfs}{\mathbf{s}}

\newcommand{\bft}{\mathbf{t}}

\newcommand{\bfV}{\mathbf{V}}

\newcommand{\bfx}{\mathbf{x}}

\newcommand{\bbB}{\mathbb{B}}

\newcommand{\bbL}{\mathbb{L}}

\newcommand{\bbQ}{\mathbb{Q}}

\newcommand{\mn}{\medskip\noindent}
\newcommand{\sn}{\smallskip\noindent}
\newcommand{\bn}{\bigskip\noindent}

\newcommand{\cA}{\mathscr{A}}

\newcommand{\cD}{\mathscr{D}}
\newcommand{\cE}{\mathscr{E}}
\newcommand{\cF}{\mathscr{F}}

\newcommand{\cL}{\mathscr{L}}
\newcommand{\cM}{\mathscr{M}}

\newcommand{\cP}{\mathscr{P}}

\newcommand{\cU}{\mathscr{U}}

\newcommand{\clH}{\mathcal{H}}

\newcommand{\clT}{\mathcal{T}}

\newcommand{\gA}{\mathfrak{A}}
\newcommand{\gB}{\mathfrak{B}}

\newcommand{\gk}{\mathfrak{k}}

\newcommand{\eps}{\varepsilon}
\newcommand{\lh}{{\ell\kern-.25ex g}}
\newcommand{\rest}{\restriction}

\newcommand{\caret}{{\char 94}}
\newcommand{\LL}{\langle}
\newcommand{\RR}{\rangle}

\newcommand{\olsi}[1]{\,\overline{\!{#1}}} 

\newcommand*{\defeq}{\mathrel{\vcenter{\baselineskip0.5ex \lineskiplimit0pt\hbox{\scriptsize.}\hbox{\scriptsize.}}}=}

\newcount\skewfactor
\def\mathunderaccent#1#2 {\let\theaccent#1\skewfactor#2
\mathpalette\putaccentunder}
\def\putaccentunder#1#2{\oalign{$#1#2$\crcr\hidewidth
\vbox to.2ex{\hbox{$#1\skew\skewfactor\theaccent{}$}\vss}\hidewidth}}

\newbox\noforkbox \newdimen\forklinewidth
\setbox0\hbox{$\textstyle\bigcup$}
\setbox1\hbox to \wd0{\hfil\vrule width \forklinewidth depth \dp0
                        height \ht0 \hfil}
\setbox\noforkbox\hbox{\box1\box0\relax}
\def\unionstick{\mathop{\copy\noforkbox}\limits}
\def\nonfork#1#2_#3{#1\unionstick_{\textstyle #3}#2}
\def\nonforkin#1#2_#3^#4{#1\unionstick_{\textstyle #3}^{\textstyle
    #4}#2}
\setbox0\hbox{$\textstyle\bigcup$}
\setbox1\hbox to \wd0{\hfil{\sl /\/}\hfil}
\setbox2\hbox to \wd0{\hfil\vrule height \ht0 depth \dp0 width
                                \forklinewidth\hfil}
\newbox\doesforkbox
\setbox\doesforkbox\hbox{\box1\box0\relax}
\def\nunionstick{\mathop{\copy\doesforkbox}\limits}

\def\fork#1#2_#3{#1\nunionstick_{\textstyle #3}#2}
\def\forkin#1#2_#3^#4{#1\nunionstick_{\textstyle #3}^{\textstyle
    #4}#2}

\newcommand{\stickT}{%
\setbox255=\hbox{\raise1ex\hbox{$\hspace{0.2pt}\,\bullet\,$}}
\mathord{\rlap{\hbox to\wd255{\hss\hbox{$|$}\hss}}
\box255}
}
\newcommand{\stickS}{%
\setbox255=\hbox{\raise0.6ex\hbox{$\scriptstyle\bullet$}}
\mathord{\rlap{\hbox to\wd255{\hss\hbox{$\scriptstyle|$}\hss}}
\box255}
}

\begin{document}
\makeatletter\def\shfiuwefootnote{\gdef\@thefnmark{}\@footnotetext}\makeatother\shfiuwefootnote{Version 2023-06-06. See \url{https://shelah.logic.at/papers/1019/} for possible updates.}

\title[Model Theory]{Model theory for a compact cardinal \\
1019}
\author {Saharon Shelah}
\address{Einstein Institute of Mathematics\\
Edmond J. Safra Campus, Givat Ram\\
The Hebrew University of Jerusalem\\
Jerusalem, 9190401, Israel\\
 and \\
 Department of Mathematics\\
 Hill Center - Busch Campus \\ 
 Rutgers, The State University of New Jersey \\
 110 Frelinghuysen Road \\
 Piscataway, NJ 08854-8019 USA}
\email{shelah@math.huji.ac.il}
\urladdr{http://shelah.logic.at}
\thanks{For earlier versions up to 2019, the author thanks Alice Leonhardt for the beautiful typing. The author would like to thank the Israel Science Foundation
for partial support of this research (Grant No. 1053/11).  
First typed May 10, 2012.\\
For the latest version, the author would like to thank the NSF and BSF for partially supporting this research --- NSF-BSF 2021: grant with M. Mailliaris, NSF 2051825, BSF 3013005232 (2021/10--2026/09). The author is also grateful to an individual who wishes to remain anonymous for generously funding typing services, and thanks Matt Grimes for the careful and beautiful typing. References like \cite[Th0.2=Ly5]{Sh:950} mean that the internal label of Theorem 0.2 in \cite{Sh:950} is \textsf{y5}. The reader should note that the version on my website is usually more up to date than the one on arXiv.}

\makeatletter
\@namedef{subjclassname@2020}{\textup{2020} Mathematics Subject Classification}
\makeatother


\date {June 5, 2023}

\begin{abstract}
We would like to develop classification theory for $T$, a complete theory in
$\bbL_{\theta,\theta}(\tau)$ when $\theta$ is a compact cardinal.  
We already have bare bones stability theory and it seemed we can go
no further.  Dealing with ultrapowers (and ultraproducts) naturally we
restrict ourselves to ``$D$ a $\theta$-complete ultrafilter 
on $I$, probably $(I,\theta)$-regular".  The basic theorems of model
theory work and
can be generalized (like \L os' theorem), but
can we generalize deeper parts of model theory?

The first section is trying to sort out what occurs to the notion of
``stable $T$" for complete $\bbL_{\theta,\theta}$-theories $T$.  We
generalize several properties of complete first order $T$, equivalent
to being stable (see \cite{Sh:c}) and find out which implications hold
and which fail.

In particular, can we generalize stability enough to 
generalize \cite[Ch.VI]{Sh:c}?  Let us concentrate on saturation in the
local sense (types consisting of instances of one formula).
We prove that at least we can
characterize the $T$-s (of cardinality $\le \theta$ for simplicity) 
which are minimal for appropriate cardinal $\lambda \ge 2^\kappa +
|T|$ in each of the following two senses.  
One is generalizing Keisler order $\triangleleft$ which
measures how saturated are ultrapowers.  Another generalizes the
results on $\triangleleft^*$. That is, we ask: ``Is there an
$\bbL_{\theta,\theta}$-theory $T_1 \supseteq T$ of cardinality $|T| +
2^\theta$ such that for every model $M_1$ of $T_1$ of cardinality 
$> \lambda$, the $\tau(T)$-reduct $M$ of $M_1$ is 
$\lambda^+$-saturated?"  Moreover, the two versions of
stable used in the characterization are different.
\end{abstract}

\maketitle
\numberwithin{equation}{section}
\setcounter{section}{-1}
\newpage

\centerline {Annotated Content}

\noindent
\S0 \quad Introduction \hfill pg. \pageref{0}
\begin{enumerate}
    \item[\S(0A)] Background and results \hfill pg.\pageref{0A}

    \item[\S(0B)] Preliminaries \hfill pg. \pageref{0B}
\end{enumerate}
\medskip

\noindent
\S1 \quad Basic Stability \hfill pg. \pageref{1}
\mn
\begin{enumerate}
\item[${{}}$]  [We try to sort out several natural generalizations of
  ``$T$ is stable" and give examples to show they are different.]
\end{enumerate}
\medskip

\noindent
\S2 \quad Saturation of ultrapowers \hfill pg. \pageref{2}
\mn
\begin{enumerate}
    \item[${{}}$]  [We define versions of saturation and give examples 
    to illustrate the difference with the first order case.  
    Then we define the generalization of $\triangleleft$, Keisler's order to $\bbL_{\theta,\theta}$; we also generalize $\triangleleft^*$.]
\end{enumerate}
\medskip

\noindent
\S3 \quad The $\ncp$ and Local Minimality \hfill pg.\pageref{3}
\mn
\begin{enumerate}
    \item[${{}}$]  [We characterize the $T$-s which are minimal in several senses, where $T$ is a complete $\bbL_{\theta,\theta}$-theory with no model of cardinality $< \theta$.  
    The first version is: there is $T_1 \supseteq T$ of cardinality 
    $\le |T| + \theta$ such that for every $M_1 \models T_1$, 
    $M_1 \rest \tau(T)$ is locally 
    $(\|M\|,\theta,\bbL_{\theta,\theta})$-saturated.  
    The second version is: the ultra-power $M^I/D$ is locally $(\lambda^+,\theta,\bbL_{\theta,\theta})$-saturated for 
    every model $M$ of $T$ and $\theta$-complete 
    $(\lambda,\theta)$-regular ultrafilter $D$ on $\lambda$.  
    We also give an example to show that those two properties
    are not equivalent.  Above, ``locally" means types involving 
    instances $\varphi(\bar x,\bar a)$ of just one formula 
    $\varphi(\bar x,\bar y) \in \bbL_{\theta,\theta}$.]
\end{enumerate}
\medskip

\noindent
\S4 \quad Global $\cp$ and Full Minimality \hfill  pg.\pageref{4}
\mn
\begin{enumerate}
\item[${{}}$]  [We replace ``local" in \S3 by ``full types" and restrict
  ourselves to the case $|T| \le \theta$, parallel to demanding ``$T$
  countable" in the first order case.]
\end{enumerate}
\newpage

\section {Introduction} \label{0}

\subsection {Background and results} \label{0A}
\bigskip

In a model theory class in Winter 2012, I tried to explain a
position I have held for a long time: model theory can deal 
extensively with $\bbL_{\lambda^+,\aleph_0}$-classes and a.e.c. 
\underline{However}, while we can generalize basic model theory to 
$\bbL_{\lambda,\kappa}$-classes with $\lambda \ge \kappa > \aleph_0$ 
(see \cite{Dic85}), we cannot do considerably more.  
The latter logics are known to have downward LST theorems 
and various connections to large cardinals and
 consistency results, and only rudimentary stability theory (see
\cite{Sh:300a}).  Note that, e.g. there is
$\psi \in \bbL_{\aleph_1,\aleph_1}$ such that $M \models \psi$ iff $M$
is isomorphic to $(\bfL_\alpha,\in)$ for some ordinal $\alpha$
such that $\beta < \alpha \Rightarrow [\bfL_\beta]^{\le \aleph_0}
\cap \bfL \subseteq \bfL_\alpha$.  Hence, assuming $\bfV = \bfL$, 
if $\mu > \cf(\mu) = \aleph_0$ then
$\psi$ has a model of cardinality $\mu$ and every
model $M$ of $\psi$ of cardinality $\mu$ is isomorphic to 
$(\bfL_\mu,\in)$.  It follows that, e.g. for every second order sentence
$\varphi$ there is $\psi \in \bbL_{\aleph_1,\aleph_1}$ which is
categorical in the cardinal $\lambda$ iff 
$(\exists \mu)[\bfL_\mu \models \varphi$ and $\lambda = \mu^{+ \omega}]$: so the categoricity spectrum
is not so nice.  Similar results hold if, e.g. $0^\#$ does not exist,
noting that if $\ge \theta^\#$ and $\mu > \cf(\mu) = \aleph_0$ then
for some real $r$, $\mu$ has cofinality $\aleph_0$ in $\bfL[r]$.
  Such views have been quite general --- see V\"a\"an\"anen's 
book \cite{Van11}.

This work is dedicated to starting to try to disprove this for the logic
$\bbL_{\theta,\theta}$ for $\theta > \aleph_0$ a compact cardinal.
Still, \L os' theorem on ultra-products was known to generalize, so let
us review the background in this direction.

In the sixties, ultra-products were very central in model theory, see
e.g. the books \cite{BeSl69} and \cite{CK73}.

Concerning isomorphisms of ultrapowers see Keisler \cite{Ke61} and then Shelah
\cite{Sh:13}; later for infinitary logics see
Hodges-Shelah \cite{Sh:109}.

In \cite{Sh:797}, the logic $\bbL^1_\theta$ is introduced.  By
\cite{Sh:1101}, elementary equivalence for $\bbL^1_\theta$ is
characterized by isomorphic ultra-limits; this was originally part of
the present paper (it was called \S3).

Here we deal with the logic $\bbL_{\theta,\theta}$ itself.
We are mainly interested in generalizations of \cite[Ch.VI]{Sh:c}, on
Keisler order $\triangleleft$ and saturation of ultra-powers and the
order $\triangleleft^*$ from \cite{Sh:500}. See history there, in 
\cite{Sh:c} and recent works with Malliaris (\cite{Sh:996}, \cite{Sh:997},
\cite{Sh:998}) dealing with unstable $T$-s and lately
\cite{Sh:1050}, \cite{Sh:1051}, \cite{Sh:1069}, \cite{Sh:1070}. 

\bn
In particular, after \cite[Ch.VI]{Sh:c} the picture was:

\begin{theorem}\label{v1}
Assume $T$ is a complete countable first order theory.

\sn
1) The following conditions are equivalent, for any 
$\lambda \ge 2^{\aleph_0}$:
\mn
\begin{enumerate}
    \item[$(a)''$]  If $D$ is a regular ultrafilter on $\lambda$ and 
    $M$ is a model of $T$ \underline{then} $M^\lambda/D$ is $\lambda^+$-saturated.\footnote{We can use ``$2^\lambda$-saturated".}
\sn
    \item[$(b)''$]  There is a first order theory $T_1 \supseteq T$ 
    such that $M_1 \models T_1 \Rightarrow M_1 \rest \tau(T)$ is 
    locally saturated (i.e. for types 
    $\subseteq \{\varphi(\bar x,\bar a) : a \in {}^{\lh(\bar y)}(M_1)\}$ for some $\varphi = \varphi(\bar x,\bar y)$.)
\sn
    \item[$(c)''$]   $T$ is stable\footnote{For first order $T$, stability follows from ``without the f.c.p."} 
    without the $\fcp$
\sn
    \item[$(d)''$]   Like $(b)''$, but $|T_1| = \aleph_0$.
\end{enumerate}

\noindent
2) The following conditions are equivalent:
\mn
\begin{enumerate}
    \item[$(a)$]  If $\bfx = \LL D_\alpha : \alpha < \delta\RR$, where
    $\delta$ is a limit ordinal and $D_\alpha$ is a regular ultrafilter
    on a cardinal $\lambda_\alpha$ for each $\alpha < \delta$,
    \underline{then} for any (equivalently, `some') model $M$ of $T$,
    $M_\delta$ is $\sup\{2^{\lambda_\alpha} : \alpha < \delta\}$-saturated, where $M_\delta$ is ultra-limit of $M$ by $\bfx$ 
    (i.e. $M_\alpha$ is $\prec$-increasing continuous for 
    $\alpha \le \delta$, $M_0 = M$, and
    $M_{\alpha +1} = M^{\lambda_\alpha}_\alpha/D_\alpha$).
\sn
    \item[$(b)$]  There is a first order theory $T_1 \supseteq T$ 
    such that $M_1 \models T_1 \Rightarrow M_1 \rest \tau(T)$ is saturated.
\sn
    \item[$(c)$]  $T$ is superstable without the $\fcp$
\sn
    \item[$(d)$]  Like $(b)$, but $|T_1| = 2^{\aleph_0}$.
\end{enumerate}
\mn
3) The following conditions are equivalent:
\mn
\begin{enumerate}
    \item[$(b)'$]  Like $(b)$, but $|T_1| = \aleph_0$.
\sn
    \item[$(c)'$]  $T$ is $\aleph_0$-stable without the $\fcp$
\end{enumerate}
\mn
See more in \cite{Sh:570} and \cite{Sh:500}.
\end{theorem}

\bn
Our main topic is generalizing results like \ref{v1}, replacing first order
logic with $\bbL_{\theta,\theta}$, so ``countable" is replaced by ``of
cardinality $\le\theta$".  More specifically, one aim is to
characterize the complete $\bbL_{\theta,\theta}$-theories $T$ such that
for some $\bbL_{\theta,\theta}$-theory $T_1$ extending $T$, for every
model $M_1$ of $T_1$, the $\tau(T)$-reduct of the model $M_1$ is
(locally) saturated. Such $T$ will be called (locally) minimal.  The
main conclusions are \ref{a33}, \ref{a33g}, \ref{a54}.

Note that $(a)'' \Leftrightarrow (c)''$ of Theorem \ref{v1}(1)
characterizes when $T$ is $\triangleleft_\lambda$-minimal and even
$\triangleleft$-minimal (but not $\triangleleft_\lambda$-minimal in
 the case $\aleph_0 < \lambda < 2^{\aleph_0}$, on it see \cite{Sh:c}).
There is much more to be said on this order.

Analogously, $(b) \Leftrightarrow (c)$ of Theorem \ref{v1} is related to
the partial orders $\triangleleft^*,\triangleleft^*_\lambda$ implicitly
investigated in \cite[Ch.VI]{Sh:c} but introduced in \cite{Sh:500}; see
more on them in Dzamonja-Shelah \cite{Sh:692}, Shelah-Usvyatsov
\cite{Sh:844} and lately Malliaris-Shelah (\cite{Sh:1051}); related
is Baldwin-Grossberg-Shelah \cite{Sh:570}.

But in our context trying to generalize Theorem \ref{v1} (i.e. the
minimal case) was hard enough.  In fact, there is a problem already 
in generalizing the notion of being stable.  In \S1 we suggest some 
reasonable definitions and try to map their relations.  Note that 
those generalizations are really very different in the present 
context (though equivalent for the first order case).  
For some versions, some ``unstable" $T$-s are categorical in 
all relevant $\lambda$-s; while other ``unstable" versions 
imply the maximal number of models up to isomorphism in relevant
cardinalities, and some ``stable $T$-s" have an intermediate 
behaviour (i.e. $\dot I(\lambda,T) = \lambda^+$).

To get sufficient conditions on $T$ for having many models, 
we may consider the tree ${}^{\theta \ge}\!\lambda$ and try to
combine it with the identities for $(\aleph_1,\aleph_0)$ (see
\cite{Sh:74}) which is a kind of the relevant indiscernibility; 
we hope to deal with this in \cite{Sh:F1396}.

Originally we were interested in generalizing the characterization of
the minimal theories in Keisler order 
$(\triangleleft,\triangleleft_\lambda)$, where $T$ is 
$\triangleleft_\lambda$-bigger if, for fewer
regular ultrafilters $D$ on the cardinal $\lambda$, $M^\lambda/D$
is $\lambda^+$-saturated for some (equivalently, `any') model of $T$.  

The earlier version was flawed but we succeed in characterizing the
$\blacktriangleleft^*_{\lambda,\theta}$-minimal ones (see \S3).  
Later we also get the characterization of the
$\blacktriangleleft_{\lambda,\theta}$-minimal ones where
$\blacktriangleleft_{\lambda,\theta}$ is defined below,
\underline{but} we use a different version of stable.

Of course, before all this we have to define saturation and local
saturation.  This is straightforward, but `unfortunately' two wonderful
properties true in the first order case are missing: existence and uniqueness.

The main achievements are in \S3,\S4: first (in \ref{a33}), a 
characterization of the (locally) minimal theories as stable with 
$\theta$-n.c.p. under reasonable definitions (see Definition \ref{a9}).
But unlike the first order case, some stable theories (even just
theories of one equivalence relation) are maximal.  In fact we get two
characterizations: one for the local version (dealing with types
containing formulas
$\varphi(\bar x_{[\eps]},\bar a)$ only for one $\varphi$ and various
$\bar a$-s) and another for the global one (naturally for theories 
$T$ with $|T| = \theta$).  Second (in \ref{a33g}), we characterize the
$\blacktriangleleft_{\lambda,\theta}$-minimal $T$ as definably 
stable with the $\theta$-n.c.p.

We may hope this will help us to resolve the categoricity spectrum.
It is natural to try to first prove that having long linear orders 
implies many models.  But this is not so --- see \ref{y45} --- so the 
situation has a marked difference from the first order case. 
We hope to continue this in \cite{Sh:F1396} and see the related
\cite{Sh:1064}; note that criteria for ``there is no universal model of
$T$ in $\lambda$" help to prove the non-categoricity of 
$T \subseteq \bbL_{\theta,\theta}$ in $\lambda$.  
See survey \cite{Dj05} and the recent \cite{Sh:F1808}.

This work was presented in a lecture in MAMLS meeting, Fall 2012 and
in courses in The Hebrew University, Spring 2012 and 2013.

We thank Doron Shafrir for (in late 2013) proof-reading, pointing out
several problematic claims (subsequently some were withdrawn, some
changed, some given a full proof) and rewriting the proof of
\ref{a23}(3).

We thank the referee for many helpful remarks.
\bigskip

\centerline {$* \qquad * \qquad *$}

\medskip

\begin{discussion}\label{v2}
1) We may wonder, for $\theta > \aleph_0$ a compact cardinal what about
$\bbL_{\theta,\aleph_0}$-theories?

\sn
2) Recall the logic from \cite[\S2]{Sh:271}; that is, given two compact
cardinals $\kappa > \theta > \aleph_0$, a logic
$\bbL_{\kappa/\theta,\kappa/\theta}$ is defined and proved to be
``nice", e.g. it is $\lambda$-compact for $\lambda < \theta$, has
interpolation, has downward $\LST$ property down to $\kappa$ and the
upward $\LST$ property for models of cardinality $\ge \lambda$ but is
not $\theta^+$-compact.

\sn
3) On the classical results on $\bbL_{\lambda,\kappa}$ see
e.g. \cite{Dic85}; on ``when for given $M_1,M_2$ there are $I$ and $D \in
\uf_\theta(I)$ such that $M^I_1/D \cong M^I_2/D$", see Hodges-Shelah
\cite{Sh:109}.

\sn
4) Recent close works are Malliaris-Shelah \cite{Sh:999} which deals with
$\kappa$-complete ultrafilters (on sets and relevant Boolean
algebras) on the way to understanding the amount of saturation of
ultra-powers by regular ultrafilters.  On reduced power, see
\cite{Sh:1064}.

\sn
5) Concerning dependent (non-elementary) classes, see also
Kaplan-Lavi-Shelah \cite{Sh:1055}.

\sn
6) Is the lack of uniqueness of saturation a sign this is a bad choice?
It does not seem so to me.  

\sn
7) If we insist on ``union of $\prec_{\cL}$-increasing 
countable chain" is an $\prec_{\cL}$-extension, we can
restrict ourselves to $\bbL^1_\theta$, but what about unions of length
$\kappa \in \Reg \cap (\aleph_0,\theta)$?  If we restrict our logic
as in $\bbL^1_\theta$ for all those $\kappa < \theta$ maybe we get close to
a.e.c., or get an interesting new logic with EM models (as indicated in
\cite{Sh:797}, \cite{Sh:893}).

\sn
8) Presently, our intention here is to show 
$\bbL_{\theta,\theta}$ has a model theory; in
particular, classification theory.  At this point having found
significant dissimilarities to the first order case on the one hand,
and solving the parallel of serious theorems on the other hand, there
is no reason to abandon this direction.
\end{discussion}

\noindent
We may wonder

\begin{question}\label{v4}
Characterize the (first order complete) $T$ such that 
$M^\lambda/D$ is not $\lambda^+$-saturated
whenever $M$ is a $\lambda$-saturated model of $T$, 
$\lambda \ge \theta > \aleph_0$, $D$ a $(\lambda,\theta)$-regular
$\theta$-complete ultrafilter on $\lambda$.
\end{question}

\begin{question}\label{v7}
Can we prove nice things on the following logics?

\mn
\begin{enumerate}
    \item[(A)]  let $\bbL^*_\kappa$ be 
    \begin{align*}
        \big\{\psi : &\ \text{for every $\mu < \kappa$ large, enough we have $\psi \in \bbL_{\mu^+,\mu^+}$,}\\ 
        &\text{ and if } \LL M_s : s \in I\RR 
        \text{ is $\prec_{\bbL_{\mu^+,\mu^+}}$-increasing,}\\ 
        &\text{ $I$ a directed partial order then }
        \textstyle\bigcup\limits_{s} M_s \models \psi \Leftrightarrow \bigwedge\limits_{s} M_s \models \psi \big\}.
    \end{align*}
    How close is
    $\bbL^*_\kappa$ to a.e.c. when $\kappa$ is a compact cardinal?
\sn
    \item[(B)]   As above, but $I$ is linearly ordered.
\end{enumerate}
\end{question}

\bigskip

\subsection {Preliminaries} \label{0B} \

\bigskip

\begin{hypothesis}\label{w0}
$\theta$ is a compact uncountable cardinal (of course, we use only
restricted versions of this).
\end{hypothesis}

\begin{notation}  \label{w2}
1) Let $\varphi(\bar x)$ mean: $\varphi$ is a formula of 
$\bbL_{\theta,\theta}$, $\bar x$ is a sequence of
variables with no repetitions including the variables occurring
freely in $\varphi$, also $\lh(\bar x) < \theta$ if not said
 otherwise.  We use $\varphi,\psi,\vartheta$ to denote formulas and
 for a logical statement $\{\st\}$ let
 $\varphi^{\st}$ or $\varphi^{[\st]}$ or $\varphi^{\iif(\st)}$ 
be $\varphi$ if $\st$ is true or 1 and be $\neg\varphi$ if $\st$ is false or 0.

\sn
2) For a set $u$, usually of ordinals, let $\bar x_{[u]}
 = \LL x_\eps:\eps \in u \RR$, now $u$ may be an
 ordinal but, e.g. if $u = [\alpha,\beta)$ we may write $\bar
 x_{[\alpha,\beta)}$; similarly for $\bar y_{[u]},\bar z_{[u]}$; let
 $\lh(\bar x_{[u]}) = u$.

\sn
3) $\tau$ denotes a vocabulary, i.e. a set of predicates and function symbols
   each with $< \theta$ places.

\sn
4) $T$ denotes a theory in $\bbL_{\theta,\theta}$; usually complete 
in the vocabulary $\tau_T$ and with a model of cardinality 
$\ge \theta$ if not said otherwise.  

\sn
5) Let $\Mod_T$ be the class of models of $T$.

\sn
6) For a model $M$ let its vocabulary be $\tau_M$.
\end{notation}

\begin{notation}  
\label{w4}
1) $\eps,\zeta,\xi$ are ordinals $< \theta$.

\sn
2) For a linear order $I$ let $\comp(I)$ be its completion.
\end{notation}

\begin{definition}  
\label{w8}
1) Let $\uf_\theta(I)$ be the set of $\theta$-complete ultrafilters on
   $I$, non-principal if not said otherwise.
Let $\fil_\theta(I)$ be the set
of $\theta$-complete filters on $I$; mainly we use
$(\theta,\theta)$-regular ones (see below).

\sn
2) The filter $D \in \fil_\theta(I)$ is called 
$(\lambda,\theta)$-regular \underline{when}  there is a
   witness $\olsi w = \LL w_t:t \in I\RR$ which means: $w_t \in
[\lambda]^{< \theta}$ for $t \in I$ and $\alpha < \lambda \Rightarrow
   \{t:\alpha \in w_t\} \in D$.

\sn
3) Let $\ruf_{\lambda,\theta}(I)$ be the set of $(\lambda,\theta)$-regular 
$D \in \uf_\theta(I)$; let $\rfil_{\lambda,\theta}(I)$ be the 
set of $(\lambda,\theta)$-regular $D \in \fil_\theta(I)$;
when $\lambda= |I|$ we may omit $\lambda$; so necessarily $\lambda \le \theta$.

\sn
4) For $S \subseteq \Card \cap \theta$ with $\sup(S) = \theta$ 
and $D \in \uf_\theta(I)$ which is not $\theta^+$-complete let
$\lcr(S,D) = \min\{\mu:\mu \ge \theta$ and 
for some $f \in {}^I S$ we have $\mu =
|\prod\limits_{s \in I} f(s)/D|\}$ and let $\Cr(S,D) =
\{\mu$: for some $f \in {}^I S$ the cardinality of $\prod\limits_{s \in I}
f(s)/D$ is $\mu\}$. 
\end{definition}

\noindent
Note that
\begin{observation}
\label{w10}
If $S = \Card \cap \theta$ and $D \in \uf_\theta(I)$ and $\mu$ is the
cardinal $\theta^I/D$ \underline{then}  $\lcr(S,D)$ is $\theta$ and $\Cr(S,D)$
 is $\Card \cap \mu^+$ or $\Card \cap \mu$.  Moreover, if $D$ is
 $(\lambda,\theta)$-regular then $\Cr(S,D) \nsubseteq 2^\lambda$ hence
 $|I| = \lambda \Rightarrow 2^\lambda \in \Cr(S,D)$; and so
$|I| = \lambda \Rightarrow 2^\lambda = \max(\Cr(S,D))$.
\end{observation}

\begin{PROOF}{\ref{w10}}
E.g., concerning the second sentence assume that $D$ is
$(\lambda,\theta)$-regular and choose $\olsi w = \LL w_s:s \in
I\RR$ witnessing it, i.e. $w_s \in [\lambda]^{< \theta}$ and
$\alpha < \lambda \Rightarrow A_\alpha \defeq \{s \in I:\alpha \in w_s\}$
belongs to $D$.  We define $f \in I_S$ by $f(m) = \min(S \setminus
2^{|w_s|})$, hence $f(s) \in S$ and let $\LL u_{s,i}:i <
2^{|w_s|}\RR$ list $\cP(w_s)$.

Now for every $u \subseteq \lambda$ let $f_u \in \prod\limits_{s \in
  I} f(s)$ be defined by: $f_u(s)$ is the $i < 2^{|w_s|} < f(s)$ such
that $u \cap w_t = u_{s,i}$.

So
\mn
\begin{enumerate}
\item[(a)]  $\{f_u/D:u \subseteq \lambda\}$ is a subset of
  $\prod\limits_{s} f(s)/D$ and
\sn
\item[(b)]  if $u_1 \ne u_2 \subseteq \lambda$ then $f_{u_1}/D \ne
  f_{u_)2}/D$.
\end{enumerate}
\mn
[Why?  Choose $\alpha \in u_1 \triangle u_2$, hence $\{s \in
I:f_{u_1}(s) \ne f_{u_2}(s)\} \supseteq \{s:\alpha \in w_s\} \in D$.]

Together we are done proving $\Cr(S,D) \nsubseteq 2^\lambda$.  Lastly,
if $I= \lambda$ then $g \in {}^I S \Rightarrow |\prod\limits_{s \in I}
g(s)/D| \subseteq |\prod\limits_{s} g(s)| \le \theta^{|I|} =
\theta^\lambda = 2^\lambda$ well assuming $0 \notin S$ for transparency.
\end{PROOF}

\begin{notation}  
\label{x0}
1) A vocabulary $\tau$ means with $\arity(\tau) \le \theta$ if not said
otherwise, where $\arity(\tau) = \aleph_0 + \sup\{|\arity(P)|^+:P$ 
is a predicate (or function symbol) from $\tau\}$, of course, where
$\arity(P)$ is the number of places of $P$.

\sn
2) If $A \subseteq N,\bar a \in {}^\eps\! N$
   and $\Delta \subseteq \bbL_{\theta,\theta}(\tau_M)$ \underline{then} 
$\tp_\Delta(\bar a,A,N) = \{\varphi(\bar x_{[\eps]},\bar b):
\varphi(\bar x_{[\eps]},\bar y) \in \Delta,N \models
``\varphi[\bar a,\bar b]"$ and $\bar b \in {}^{\lh(\bar y)}\!M\}$.

\sn
3) $\bfS^\eps_\Delta(A,M) = \{\tp_\Delta(\bar a,A,N)$: for
   some $N,M \prec_{\bbL_{\theta,\theta}} N$ and $\bar a \in
   {}^\eps\! N\}$.

\sn
4) If $\Delta = \bbL_{\theta,\theta}$ then we may omit $\Delta$.

\sn
4A) If $\Delta$ is the set of quantifier free formulas from
$\bbL(\tau_N)$, we may write $\tp_{\qf}$ instead of $\tp_\Delta$.
\end{notation}

\begin{definition}  
\label{x1}
1) $\bbL_{\theta,\theta}(\tau)$ is the set of formulas of
   $\bbL_{\theta,\theta}$ in the vocabulary $\tau$.

\sn
2) For $\tau$-models $M,N$ let $M \prec_{\bbL_{\theta,\theta}} N$
mean: if $\varphi(\bar x) \in
\bbL_{\theta,\theta}(\tau_M)$ and $\bar a \in {}^{\lh(\bar x)}\!M$
 then $M \models \varphi[\bar a] \Leftrightarrow N \models
   \varphi[\bar a]$. 
\end{definition}

\begin{definition}  \label{x2}
For a set $v$ of ordinals, a sequence $\bar u = \LL u_\alpha:\alpha
\in v\RR$ and models $M_1,M_2$ of the same vocabulary $\tau$ and
$\Delta \subseteq \bbL_{\theta,\theta}(\tau)$ a set of formulas we define a
game $\Game = \Game_{\Delta,\bar u}(M_1,M_2)$ but when $(\forall
\alpha \in v)(u_\alpha = u)$ we may write
$\Game_{\Delta,u,v}(M_1,M_2)$:
\begin{enumerate}
    \item  A play lasts some finite number of moves not known in advance.
\sn
    \item  In the $n^\tthh$ move the antagonist chooses
    \begin{enumerate}
        \item $\alpha_n \in v$ such that 
        $m < n \Rightarrow \alpha_n < \alpha_m$.

        \item A sequence $\LL a_{n,i,\ell}:(n,i,\ell) \in I\RR$, 
        where
        \begin{itemize}
            \item  $I = \{(n,i,\ell_{n,i}) : i \in u_{\alpha_n}\}$

            \item $\ell_{n,i} = \ell(n,i) \in \{1,2\}$

            \item $ a_{n,i,\ell(n,i)} \in M_{\ell_{n,i}}$
        \end{itemize}
    \end{enumerate}
\sn
    \item  In the $n^\tthh$ move (after the antagonist's move) 
    the protagonist chooses 
    $a_{n,i,3 - \ell(n,i)} \in M_{3-\ell(n,i)}$ for 
    $i \in u_{\alpha_n}$.
\sn
    \item  The play ends when the antagonist cannot choose $\alpha_n$.
\sn
    \item  The protagonist wins a play \underline{when}:
    \begin{itemize}
        \item The set $\{(a_{m,i,1},a_{m,i,2}):i \in u_{\alpha_m}$ and the $m^\tthh$ move was done$\}$ is a function,

        \item it is a partial one-to-one function from $M_1$ into $M_2$, and

        \item it preserves satisfaction of $\Delta$-formulas and their negations.
    \end{itemize}
\end{enumerate}
\end{definition}

\noindent
We know (see, e.g. \cite{Dic85})

\begin{fact}\label{x4}
The $\tau$-models $M_1,M_2$ are $\bbL_{\theta,\theta}$-equivalent 
\underline{iff} for every set $\Delta$ of $< \theta$ atomic formulas 
and $\alpha,\beta < \theta$ the protagonist wins in 
the game $\Game_{\Delta,\alpha,\beta}(M_1,M_2)$.
\end{fact}

\noindent
And, of course

\begin{fact}\label{x6}
For a complete $T \subseteq \bbL_{\theta,\theta}(\tau)$,

\sn
1) $(\Mod_T,\prec_{\bbL_{\theta,\theta}})$ has amalgamation 
and the joint embedding property (JEP), that is:
\mn
\begin{enumerate}
    \item  \textbf{Amalgamation}: if $M_0
  \prec_{\bbL_{\theta,\theta}} M_\ell$
  for $\ell=1,2$ then there are $M_3,f_1,f_2,M'_1,M'_2$ such that
\sn
    \begin{enumerate}
        \item[$\bullet$]  $M_0 \prec_{\bbL_{\theta,\theta}} M_3$
\sn
        \item[$\bullet$]  For $\ell=1,2$, $f_\ell$ is a $\prec_{\bbL_{\theta,\theta}} $-embedding of $M_\ell$ 
        into $M_3$ over $M_0$. That is, for some $\tau_T$-models 
        $M'_\ell$ for $\ell=1,2$ we have 
        $M'_\ell \prec_{\bbL_{\theta,\theta}} M_3$ and $f_\ell$ 
        is an isomorphism from $M_\ell$ onto $M'_\ell$ over $M_0$.
    \end{enumerate}
\sn
    \item  \textbf{JEP}: if $M_1,M_2$ are 
    $\bbL_{\theta,\theta}$-equivalent $\tau$-models \underline{then}
    there is a $\tau$-model $M_3$ and 
    $\prec_{\bbL_{\theta,\theta}}$-embedding 
    $f_\ell$ of $M_\ell$ into $M_3$ for $\ell=1,2$.
\end{enumerate}
\mn
2) Types are well defined (see \cite{Sh:300b}); i.e. the orbital type
\textbf{tp} and the types as a set of formula
$\tp_{\bbL_{\theta,\theta}}$ are essentially equivalent. That is:
\mn
\begin{enumerate}
    \item[$(*)$]  If $M_0 \prec_{\bbL_{\theta,\theta}} M_\ell$, 
    $\zeta < \theta$, $\bar a_\ell \in {}^\zeta|M_\ell|$ for $\ell=1,2$ and
so $\tau = \tau(M_\ell)$ for $\ell=0,1,2$ \underline{then}  the 
following conditions are equivalent:
\sn
    \begin{enumerate}
        \item  The set of formulas (= type) 
        $\tp_{\bbL_{\theta,\theta}}(\bar a_1,M_0,M_1)$ is equal to
        $\tp_{\bbL_{\theta,\theta}}(\bar a_2,M_0,M_2)$ 
        (see \ref{x0}(2)). That is, if $\xi < \theta$, 
        $\bar b \in {}^\xi(M_0)$, and 
        $\varphi(\bar x_{[\zeta]},\bar y_{[\xi]}) \in \bbL_{\theta,\theta}(\tau)$ then $M_1 \models \varphi[\bar a_1,\bar b] \Leftrightarrow M_2 \models \varphi[\bar a_2,\bar b]$.
\sn
        \item  (orbital types) There are $M_3,f_1,f_2$ as in 
        \ref{x6}(1)(a) such that $f_1(\bar a_1) = f_2(\bar a_2)$.
    \end{enumerate}
\end{enumerate}
\end{fact}

\noindent
The well known generalization of \L os' theorem (see e.g. \cite{J}
or \cite{Sh:109}) is:

\begin{theorem}\label{x7}
If $\varphi(\bar x_{[\zeta]}) \in \bbL_{\theta,\theta}(\tau)$, 
$D \in \uf_\theta(I)$, $M_s$ is a $\tau$-model for $s \in I$, and
$f_\eps \in \prod\limits_{s \in I} M_s$ for $\eps < \zeta$ 
\underline{then}  
$M \models \varphi[\ldots,f_\eps/D,\ldots]_{\eps < \zeta}$ 
\underline{iff} the set 
$$\{s \in I : M_s \models \varphi[\ldots,f_\eps(s),\ldots]_{\eps < \zeta}\}$$ belongs to $D$.
\end{theorem}

\noindent
Recall

\begin{fact}\label{x8}
Assume $D \in \uf_\theta(I)$ is not $\theta^+$-complete and $\gB =
(\clH(\chi),\in,\theta)^I/D$.

\sn
1) If $\cf(\chi) \ge \theta$ and $a_\alpha \in \gB$ for 
$\alpha < \theta$ \underline{then} there is $\bar b \in \gB$ such
that $\gB \models ``\bar b$ is a sequence of length $< \theta$ 
with the $\alpha^\tthh$ element being $a_\alpha$"
for\footnote{We are identifying elements of $\clH(\chi)$ with ones of
    $\gB$ naturally, see \ref{x18}(2).  Alternatively, expand 
    $\gA = (\clH(\chi),\in,\theta)$ by having 
    $c^{\gA^+}_\alpha = \alpha$, so $c_\alpha \in \tau(\gA^+)$ 
    is an individual constant for $\alpha < \lambda$, so 
    $\gB^+ = (\gA^+)^I/D$ is an expansion of $\gB$ and
    $\gB^+ \models ``a_\alpha$ is the $c_\alpha$-th element 
    of the sequence $b$".}
every $\alpha < \theta$.

\sn
2) If $\cf(\chi) > \lambda$ and $D$ is $(\lambda,\theta)$-regular 
and $a_\alpha \in \gB$ for $\alpha < \lambda$ \underline{then} 
there is $w \in \gB$ such that 
$\alpha < \lambda \Rightarrow \gB \models ``|w| < \theta$
and $a_\alpha \in w"$, (in fact, also the inverse holds).

\sn
3) For some function $h$ from $I$ onto $\theta$, 
$D/h = \{u \subseteq \theta : h^{-1}(u) \in D\}$ 
is a normal ultrafilter on $\theta$.
\end{fact}

\begin{PROOF}{\ref{x8}}  
1) Let $a_\alpha = f_\alpha/D$ where $f_\alpha \in {}^I (\clH(\chi))$.
Let $F:I \rightarrow \theta$ be such that $\alpha < \theta
\Rightarrow \{s:\alpha \le F(s)\} \in D$, such a function $F$ 
exists by the assumption on $D$.  We define $g:I \rightarrow \clH(\chi)$ by
$$g(s) = \LL f_\alpha(s):\alpha < F(s)\RR.$$
Now $g/D$ is as required: check.

\sn
2) Similarly using $\olsi w = \LL w_s:s \in I\RR$ from \ref{w8}, so
$$g(s) = \{f_\alpha(s):\alpha \in w_s\}.$$

\noindent
3) See, e.g. \cite{Je03}.
\end{PROOF}

\noindent
Recall (see history \cite[\S1]{Sh:950}) in the literature usually 
we say ``strongly convergent" instead of ``convergent" to distinguish
from other versions; but here this is not needed.

\begin{definition}\label{x11}
Assume $\Delta \subseteq \bbL_{\theta,\theta}(\tau_M)$, $I$ is a
linear order, $\bar{\bfa} = \LL \bar a_t : t \in I\RR$ and
$t \in I \Rightarrow \bar a_t \in {}^u\! M$, and 
$$\bar\theta = \big\LL \theta_\varphi = 
\theta_{\varphi(\bar x_{[u]},\bar y)} : \varphi  =
\varphi(\bar x_{[u]},\bar y) \in \Delta \big\RR$$ 
where $\theta_\varphi$ is a cardinal $\le \theta$; if 
$\bigwedge\limits_{\varphi \in \Delta} \theta_\varphi = \sigma$ 
we may write $\sigma$; if $\sigma = \theta$ we may omit it. 

\sn
1) We say $\bar{\bfa}$ is a $(\Delta,\bar\theta)$-convergent sequence 
in $M$ \underline{when} for every 
$\varphi(\bar x_{[u]},\bar y) \in \Delta$ and 
$\bar b \in {}^{\lh(\bar y )}\!M$ there is $J \subseteq \comp(I)$ of cardinality $< \sigma$ \underline{or} 
$< \theta_{\varphi(\bar x_{[u]},\bar y)} < \theta$ respectively, 
such that:
\mn
\begin{itemize}
    \item  If $s,t \in I$ and 
    $\tp_{\qf}(s,J,\comp(I)) = \tp_{\qf}(t,J,\comp(I))$ then 
    $M \models ``\varphi[\bar a_s,\bar b] \equiv \varphi[\bar a_t,\bar b]"$.
\end{itemize}
\mn
1A) We say $\bar{\bfa}$ is a middle $(\Delta,\sigma)$-convergent sequence
\underline{when}  $\bar{\bfa}$ is $(\Delta,\bar\theta)$-convergent
for some $\bar\theta = \LL \theta_\varphi:\varphi \in
\Delta\RR$ satisfying $\varphi \in \Delta \Rightarrow
\theta_\varphi < \sigma$.  If $\sigma = \theta$ then we may omit it.

\sn
2) We say ``strictly $(\Delta,\bar\theta)$-convergent" \underline{when} 
we demand ``$J \subseteq I$;" similarly in the other variant.
\end{definition}

\begin{definition}\label{x13}
For a linear order $I$:

\sn
1) $I^*$ is its inverse, $\cf(I)$ is the cofinality of $I$ (so $0,1$
or a regular cardinal) and $\coin(I)$ is the co-initiality of
$I$ (that is, the cofinality of its inverse).

\sn
2) A cut is a pair $(C_1,C_2)$ such that $C_1$ is an initial segment
   of $I$ and $C_2 = I \setminus C_1$.

\sn
3) The cofinality $(\kappa_1,\kappa_2)$ of the cut $(C_1,C_2)$ is 
the pair $(\kappa_1,\kappa_2)$ of regular cardinals (or 0 or 1) 
such that $\kappa_1 = \cf(I \rest C_1)$, $\kappa_2 = \coin(I \rest C_2)$.

\sn
4) We say $(C_1,C_2)$ is a pre-cut of $I$ [of cofinality
$(\kappa_1,\kappa_2)$] \underline{when}  $C_1,C_2 \subseteq I$ and
$$\Big(\big\{s \in I:(\exists t \in C_1)[s \le_I t]\big\},
\big\{s \in I:(\exists t \in C_2)[t \le_I s]\big\}\Big)$$ 
is a cut of $I$ [of cofinality $(\kappa_1,\kappa_2)$].
\end{definition}

\begin{definition}
\label{x17}
0) We say $X$ respects $E$ \underline{when}  
$E$ is an equivalence relation on some set $I \supseteq X$ and 
$s\ E\ t \Rightarrow (s \in X \leftrightarrow t \in X)$.

\sn
1) We say $\bfx = (I,D,\cE)$ is a $(\kappa,\sigma)$-$\luft$ 
(limit-ultrafilter-iteration triple) \underline{when}:
\begin{enumerate}[(a)]
    \item  $D$ is a filter on the set $I$.
\sn
    \item  $\cE$ is a family of equivalence relations on $I$.
\sn
    \item  $(\cE,\supseteq)$ is $\sigma$-directed; i.e. if 
    $\alpha(*) < \sigma$ and $E_i \in \cE$ for $i < \alpha(*)$ 
    \underline{then} there is $E \in \cE$ refining $E_i$ for 
    every $i < \alpha(*)$.
\sn
    \item  If $E \in \cE$ then $D/E$ is a $\kappa$-complete 
    ultrafilter on $I/E$, where \\
    $D/E \defeq \{X/E:X \in D$ and $X$ respects $E\}$.
\end{enumerate}
\mn
1A) We say $\bfx$ is a $(\kappa,\sigma)$-$\lft$ \underline{when}  
above we weaken clause (d) to:
\mn
\begin{enumerate}
    \item[(d)$'$]  If $E \in \cE$ then $D/E$ is a $\kappa$-complete filter.
\end{enumerate}
\mn
2) Omitting ``$(\kappa,\sigma)$" means $(\theta,\aleph_0)$, recalling
$\theta$ is our fixed compact cardinal.

\sn
3) Let $(I_1,D_1,\cE_1) \le^1_h (I_2,D_2,\cE_2)$ mean that:
\mn
\begin{enumerate}[(a)]
    \item  $h$ is a function from $I_2$ onto $I_1$
\sn
    \item  If $E \in \cE_1$ \underline{then} $h^{-1} \circ E \in \cE_2$, where
    $$h^{-1} \circ E = \big\{(s,t) \in I_2 \times I_2 : h(s)\ E\ h(t) \big\}.$$

    \item  If $E_1 \in \cE_1$ and $E_2 = h^{-1} \circ E_1$ \underline{then} $D_1/E_1 = h(D_2/E_2)$.
\end{enumerate}
\end{definition}

\begin{remark}\label{x17f}
Note that in \ref{x17}(3), if $h = \id_{I_2}$ then $I_1 = I_2$. 
\end{remark}

\begin{definition}\label{x16}
Assume $\bfx = (I,D,\cE)$ is a $(\kappa,\sigma)$-l.u.f.t.

\sn
1) For a function $f$ let $\eq(f) = \{(s_1,s_2) : f(s_1) = f(s_2)\}$.
If $\bar f = \LL f_i:i < i_*\RR$ and $i < i_* \Rightarrow
\dom(f_i) = I$ then $\eq(\bar f) = \bigcap\{\eq(f_i) : i < i_*\}$.

\sn
2) For a set $U$ let $U^I|\cE = \{f \in {}^I U:\eq(f)$ 
is refined by some $E \in \cE\}$.

\sn
3) For a model $M$ let $\flp_{\bfx}(M) = M^I_D|\cE =$ 
$$(M^I/D) \rest \{f/D : f \in {}^I\! M \text{ and $\eq(f)$ is refined by some } E \in \cE\}.$$
Pedantically (as $\arity(\tau_M)$ may be $> \aleph_0$), 
$M^I_D|\cE = \bigcup\{M^I_D \rest E : E \in \cE\}$. 
($\lrp$ stands for \emph{limit reduced power}.)

\sn
4) If $\bfx$ is $\luft$ we may in (3) write $\uflp_{\bfx}(M)$.
\end{definition}

\sn
We now give the generalization of Keisler \cite{Ke63}; Hodges-Shelah
\cite[Lemma 1,pg.80]{Sh:109} is the case $\kappa = \partial$.

\begin{theorem}\label{x18}
1) If $(I,D,\cE)$ is $(\kappa,\partial)$-$\luft$, 
$\varphi = \varphi(\bar x_{[\zeta]}) \in \bbL_{\kappa,\partial}(\tau)$ 
(so $\zeta < \partial$), and $f_\eps \in M^I|\cE$ for $\eps < \zeta$
\underline{then}  $M^I_D|\cE \models \varphi[\ldots,f_\eps/D,\ldots]$ 
\underline{iff} $$\{s \in I : M \models \varphi[\ldots,f_\eps(s),\ldots]_{\eps < \zeta}\} \in D.$$

\sn
2) Moreover, $M \prec_{\bbL_{\kappa,\partial}} M^I_D/\cE$ (pedantically,
$\bfj = \bfj_{M,\bfx}$ is a $\prec_{\bbL_{\kappa,\partial}}$-elementary 
embedding of $M$ into $M^I_D/\cE$, where $\bfj(a) = \LL a:s \in I\RR/D$).

\sn
3) We define $(\prod\limits_{s \in I} M_s)^I_D|\cE$ similarly
\underline{when} the equivalence relation 
$$\{(s,t) \in I \times I : M_s = M_t\}$$ 
is refined by some $E \in \cE$.
\end{theorem}

\begin{convention}\label{x21}
1) Abusing notation in $\prod\limits_{s \in I} M_s/D$, we allow $f/D$ 
for $f \in \prod\limits_{s \in S} M_s$ when $S \in D$.

\sn
2) For $\bar c \in {}^\gamma \big(\prod\limits_{s \in I} M_s/D \big)$, 
we can choose $\LL \bar c_s : s \in I\RR$ such that 
$\bar c_s \in {}^\gamma(M_s)$ and 
$\bar c =\LL \bar c_s : s \in I \RR/D$. (This means if 
$i < \lh(\bar c)$ \underline{then} $c_{s,i} \in M_s$ and 
$c_i = \LL c_{s,i} : s \in I\RR/D$.)
\end{convention}

\begin{remark}\label{x24}
1) Why the ``pedantically" in \ref{x16}(3)?  Otherwise 
if $\bfx$ is a $(\theta,\sigma)$-$\luft$, $(\cE_{\bfx},
\supseteq)$ is not $\kappa^+$-directed, and $\kappa < \arity(\tau)$
\underline{then} defining $\uflp_{\bfx}(M)$, we have freedom: if 
$R \in \tau$ and $\arity_\tau(R) \ge \kappa$; i.e. on 
$$R^N \rest \{\bar a \in {}^{\arity(P)}\!N : \text{ no } E \in \cE
\text{ refines } \eq(\bar a)\}$$ we have no restrictions.

\sn
2) So, e.g. for categoricity we better restrict ourselves to vocabularies $\tau$ such that $\arity(\tau) = \aleph_0$.
\end{remark}

\begin{definition}\label{x28}
We say that $M$ is a $\theta$-complete model \underline{when} for
every $\eps < \theta$, $R_* \subseteq {}^\eps\! M$, and 
$F_*:{}^\eps\! M \to M$ there are $R,F \in \tau_M$ such that 
$R^M = R_* \wedge F^M = F_*$. 
\end{definition}

\begin{observation}\label{x31}
1) If $M$ is a $\tau$-model of cardinality $\lambda$ \underline{then}
there is a $\theta$-complete expansion $M^+$ of $M$ so 
$\tau(M^+) \supseteq \tau(M)$ and $\tau(M^+)$ has cardinality 
$|\tau_M| + 2^{(\|M\|^{< \theta})}$.

\sn
2) For models $M \prec_{\bbL_{\theta,\theta}} N$ and $M^+$ as 
above the following conditions are equivalent:

\begin{enumerate}
    \item[$(a)$]  $N = \uflp_{\bfx}(M)$ up to isomorphisms over 
    $M$ identifying $a \in M$ with $\bfj_{\bfx}(a) \in N$, for 
    some $(\theta,\theta)$-$\luft$ $\bfx$.
\sn
    \item[$(b)$]  There is $N^+$ such that 
    $M^+ \prec_{\bbL_{\theta,\theta}} N^+$ and $N^+ \rest \tau_M$ 
    is isomorphic to $N$ over $M$.
\end{enumerate}
\mn
3) For a model $M$, if $(P^M,<^M)$ is a $\theta$-directed 
partial order and $\chi = \cf(\chi) \ge \theta$ and $\lambda = 
\lambda^{\|M\|} + \chi$ \underline{then}  for some 
$(\theta,\theta)$-$\luft$ $\bfx$, the model $N \defeq \lup_{\bfx}(M)$ 
satisfies $(P^N,<^N)$ has a cofinal increasing sequence of
length $\chi$ and $|P^N| = \lambda$.
\end{observation}

\begin{PROOF}{\ref{x31}}
Easy, e.g. 

\sn
3) Let $M^+$ be as in part (1).  
Note that $M^+$ has Skolem functions for formulas $\varphi(\bar x,\bar
y) \in \bbL_{\theta,\theta}(\tau_{M^+})$ and let 
$T' \defeq Th_{\bbL_{\theta,\theta}}(M^+) \cup
\{P(\sigma(x_{\eps_0},\ldots,x_{\eps_i},\ldots)_{i<i(*)})
\rightarrow \sigma(x_{\eps_0},\dotsc,x_{\eps_i},
\ldots)_{i < i(*)} < x_\eps: \sigma$ is a 
$\tau(M^+)$-term so $i(*) < \theta$
and $i < i(*) \Rightarrow 
\eps_i < \eps < \lambda \cdot \chi\}$.  

Clearly
\mn
\begin{enumerate}
\item[$(*)$]   $T'$ is $(< \theta)$-satisfiable in $M^+$.
\end{enumerate}
\mn
[Why?  Because if $T'' \subseteq T'$ has cardinality $< \theta$ then
the set $u = \{\eps < \lambda \cdot \chi:x_\eps$ appears in 
$T''\}$ has cardinality $< \theta$ and let $i(*) = \otp(u)$; now 
for each $\eps \in u$ the set $\Gamma_\eps = T' \cap
\{P(\sigma(x_{\eps_0},\ldots)) \rightarrow
\sigma(x_{\eps_0},\dotsc,x_{\eps_i},\ldots)_{i<i(*)} <
x_\eps:i(*) < \theta$ and $\eps_i < \eps$ for
$i<i(*)\}$ has cardinality $< \theta$.  Now we choose $c_\eps
\in M$ by induction on $\eps \in u$ such that the
assignment $x_\zeta \mapsto c_\zeta$ for $\zeta \in \eps \cap
u$ in $M^+$ satisfies $\Gamma_\eps$, possible because
$|\Gamma_\eps| < \theta$ and $(P^M,<^M)$ is $\theta$-directed.
So the model $M^+$ with the assignment $x_\eps \mapsto c_\eps$
for $\eps \in u$ is a model of $T''$, so $T'$ is ($<
\theta)$-satisfiable indeed.]

Recalling that $|M| = \{c^{M^+}:c \in \tau(M^+)$ an individual
constant$\}$, $T'$ is realized in some
$\prec_{\bbL_{\theta,\theta}}$-elementary extension $N^+$ of $M^+$ by
the assignment $x_\eps \mapsto a_\eps$ (for $\eps 
< \lambda \cdot \chi$).  Without loss of generality, $N^+$ is the Skolem hull of
$\{a_\eps:\eps < \lambda \cdot \chi\}$, so $N \defeq N^+
\rest \tau(M)$ is as required.  Now $\bfx$ 
as required exists by part (2). 
\end{PROOF}

\begin{observation}
\label{x32}
1) If $\bfx$ is a non-trivial $(\theta,\theta)$-$\luft$ and $\chi =
\cf(\uflp(\theta <))$ \underline{then}  $\chi =\chi^{< \theta}$.

\sn
2) Also $\mu = \mu^{< \theta}$ when $\mu$ is the cardinality of
   $\lup(\theta,<)$. 
\end{observation}

\begin{PROOF}{\ref{x32}}
1) By the choice of $\bfx$ clearly $\chi = \cf(\chi) \ge \theta$.  
As $\chi$ is regular $\ge \theta$ by a theorem of Solovay \cite{So74} we have
$\chi^{< \theta} = \chi$.

\sn
2) See the statement and the proof of \ref{a27c}.
\end{PROOF}
\newpage

\section {Basic stability} \label{1}

For a complete first order $T$, being stable has many equivalent
definitions; see \cite{Sh:c}.  We define the parallel properties for a
complete $\bbL_{\theta,\theta}$-theory and try to sort out the implications.

A difference with the first order case which may be confusing is that
the existence of long orders is not so strong and does not imply other
versions of unstability, see in particular \ref{y45}.

In Definition \ref{y2} below, defining the notions ``$\iota$-unstable"
generally demand more when $\iota$ increase; it seems reasonable that
we shall order the parts of \ref{y2} in increasing order by $\iota$,
but we deviate putting ``4-unstable" just after ``1-unstable" as it is
more easy to define than 2/3-unstable.

\begin{definition}  
\label{y2}
Let $T \subseteq \bbL_{\theta,\theta}$, not necessarily complete;
below ``$T$ is $\iota$-stable" is the negation of ``$T$ is
$\iota$-unstable"; below if $\Delta = \bbL_{\theta,\theta}(\tau_T)$
then we may omit $\Delta$ except in parts (4),(4A).

\sn
1) $T$ is 1-unstable \underline{iff}  for some $\eps,\zeta < \theta$
   and formula $\varphi(\bar x_{[\eps]},\bar y_{[\zeta]}) \in
   \bbL_{\theta,\theta}$ there is a model $M$ of $T$ and $\bar
   a_\alpha \in {}^\eps\! M,\bar b_\alpha \in {}^\zeta\! M$ for
   $\alpha < \theta$ such that $M \models 
\varphi[\bar a_\alpha,\bar b_\beta]^{\iif(\alpha < \beta)}$ 
for $\alpha,\beta < \theta$.

\sn
2) We say $T$ is 4-unstable \underline{when}  there are $\varphi(\bar x,\bar y) \in
\bbL_{\theta,\theta}$ and a model $M$ of $T$ and $\bar b_\eta \in {}^{\lh(\bar y)}\!M$ for $\eta \in {}^{\theta >}2$ such that $p_\eta(\bar x) 
= \{\varphi(\bar x,\bar b_{\eta \rest \alpha})^{\iif(\eta(\alpha))}:
\alpha < \theta\}$ is a type in $M$ for every $\eta \in {}^\theta 2$,
i.e. every subset of cardinality $< \theta$ is realized.

\sn
3) For a class $\bfI$ of linear orders we say $T$ is 
$\bfI$-unstable \underline{when}  for some 
$\varphi(\bar x,\bar y) \in \bbL_{\theta,\theta}$ for every $I \in
\bfI$ there are $M$ and $\LL (\bar a_s,\bar b_s):s \in
   I\RR$ is as in part (1).  If $\bfI = \{I\}$ we may write
   $I$-unstable.  We say $T$ is $(\Delta,\bfI)$-unstable when above
   $\varphi(\bar x,\bar y) \in \Delta$.

\sn
4) We say $T$ is strongly $(\Delta,\bfI)$-unstable
\underline{when}\footnote{The difference
   between \ref{y2}(3) and \ref{y2}(4) is the ``convergent".  In part
 (5) for the applications we have in mind it 
is enough to restrict ourselves to the case $\bfI_2 =
\{\sum\limits_{i<(*)} \delta_i$: \underline{where}  
$\delta_i \in\{\theta,\theta^+\},i(*)$ an ordinal$\}$.} \, 
for some $\varphi(\bar x,\bar y) \in \Delta$ satisfying 
$\lh(\bar x) = \lh(\bar y)$ 
for every linear order $I \in \bfI$ there are $M \models T$ and
sequence $\LL \bar a_s \caret \bar b_s:s \in I\RR$ in $M$
such that:
\mn
\begin{enumerate}
\item[(a)]  $M \models \varphi[\bar a_s,\bar b_t]^{\iif(s<t)}$ for $s,t
 \in I$, 
\sn
\item[(b)] $\LL \bar a_s:s \in I\RR$ is strictly $\varphi(\bar
  x_{[\eps]},\bar y_{[\zeta]})$-convergent where $\lh(\bar a_s) =
  \eps$
\sn
\item[(c)] $\LL \bar b_s : s \in I\RR$ is strictly 
$\psi(\bar x_{[\zeta]},\bar y_{[\eps]})$-convergent where 
$\lh(\bar b_s) = \zeta$ and $\psi(\bar x_{[\zeta]},\bar
y_{[\eps]}) = \varphi(\bar y_{[\eps]},\bar x_{[\zeta]})$ also called
$\varphi^+(\bar x_{[\zeta]},\bar y_{[\eps]})$
\end{enumerate}
\mn
recalling Definition \ref{x11}(1),(2).  Let the default value
 of $\Delta$ be $\{\varphi(\bar x_{[\eps]},\bar y_{[\zeta]}),
\psi(\bar x_{[\zeta]},\bar y_{[\eps]})\}$.

\sn
4A) We say $T$ is middle $\Delta$-unstable \underline{when}  in part (4) we
replace ``strictly $\Delta$-convergent" by ``strictly middle 
$\Delta$-convergent", see Definition \ref{x11}(1),(2).  The default
value of $\Delta$ is as in part (4).

\sn
5) We say $T$ is 3-unstable \underline{when}  it is strongly $\bfI_2$-unstable
where $\bfI_2 = \{\sum\limits_{i<i(*)} I_i:i(*)$ an ordinal and for
each $i,I_i$ is anti-isomorphic to some ordinal $\delta_i,
\cf(\delta_i) \ge \theta\}$.

\sn
6) We say $T$ is 2-unstable \underline{iff}  it is $\bfI_2$-unstable.

\sn
7) We say $T$ is 5-unstable if it is $({}^{\theta >}2,<_{\lex})$-unstable.
\end{definition}

\begin{remark}  
\label{y2f}
We shall clarify all implications between 
``$\iota$-unstable" and definably stable which is defined below; this
is summed up in \ref{y62}.
\end{remark}

\begin{definition}  
\label{y3}
Let $T$ be as in \ref{y2}.

\sn
1) $T$ is definably stable (definably unstable is the negation) \underline{when} :
if $\varphi = \varphi(\bar x_{[\eps]},\bar y_{[\zeta]}) \in 
\bbL_{\theta,\theta}$ \underline{then}  there is $\psi(\bar y_{[\zeta]},
\bar z_{[\xi]}) \in \bbL_{\theta,\theta}$ such that:
\mn
\begin{enumerate}
\item[$(*)$]  if $M \prec_{\bbL_{\theta,\theta}} N$ are models of $T$
and $\bar a \in {}^\eps\! N$ \underline{then}  there is 
$\bar c \in {}^\xi M$ satisfying:
$\psi(\bar y_{[\zeta]},\bar c)$ defines $\tp_\varphi(\bar a,M,N)$, that is:
\sn
\begin{enumerate}
\item[$\bullet$]   if $\bar b \in {}^\zeta\! M$ \underline{then}  $N \models
  \varphi[\bar a,\bar b]$ iff $M \models \psi[\bar b,\bar c]$.
\end{enumerate} 
\end{enumerate}
\mn
2) We say $\varphi(\bar x,\bar y) \in 
\bbL_{\theta,\theta}(\tau_T)$ is 1-stable (for $T$) \underline{when}  \ref{y2}(1)
fails for $\varphi$ (and $T$).  Similarly for the other versions.  
We say $\varphi(\bar x,\bar y)$ is symmetrically 1-stable (for $T$) 
\underline{when}  it is 1-stable and also $\varphi^\perp(\bar y,\bar x)$ is
1-stable where $\varphi^\perp(\bar y,\bar x) = \varphi(\bar x,\bar y)$
is called the dual of $\varphi(\bar x,\bar y)$.

\sn
3) We say $T$ is $(\lambda,\Delta)$-stable \underline{when}  $\Delta \subseteq
\bbL_{\theta,\theta}(\tau_T)$ and for every model $M$ of $T$ and $A
\subseteq M$ of cardinality $\le \lambda$ and $\varphi = \varphi(\bar
x_{[\eps]},\bar y_{[\zeta]}) \in \Delta$ the set $\bfS^\eps_\varphi(A,M)$ has cardinality $\le \lambda$ \underline{where} 
$\bfS^\eps_\Delta(A,M) = \{\tp_\Delta(\bar a,A,N):N,
\bar a$ satisfy $M \prec_{\bbL_{\theta,\theta}} N,
\bar a \in {}^\eps\! N\}$.

\sn
4) We say $T$ is $\Delta$-stable \underline{when}  $T$ is
$(\lambda,\Delta)$-stable for every $\lambda = \lambda^{< \theta} +
\lambda^{|T|}$.

\sn
4A) In part 3) and 4) omitting $\Delta$ means $\Delta =
\bbL_{\theta,\theta}(\tau_T)$. 
\end{definition}

\begin{claim}  
\label{y5}
Let $T \subseteq \bbL_{\theta,\theta}$ (not necessarily complete),
$\tau = \tau(T)$ and let $\partial = (\theta + |T|)^{< \theta}$.  

\sn
1) We have $(a) \Rightarrow (b) \Rightarrow (c) \Rightarrow (x) 
\Rightarrow (f) \Rightarrow (g) \Rightarrow (h)
\Rightarrow (i) \Leftrightarrow (j)$ for $x = d,e$ where:
\mn
\begin{enumerate}
\item[$(a)$]  $T$ is 5-unstable, see \ref{y2}(7)
\sn
\item[$(b)$]   $T$ is 4-unstable, see \ref{y2}(2)
\sn
\item[$(c)$]  for some $\eps < \theta$ for every $\lambda \ge
\theta$ there are 
$A \subseteq M \models T,|A| = \lambda$ such that $\bfS^\eps(A,M) =
\{\tp_{\bbL_{\theta,\theta}}(\bar a,A,N):
M \prec_{\bbL_{\theta,\theta}} N,\bar a \in 
{}^\eps\! N\}$ has cardinality $> \lambda$
\sn
\item[$(d)$]   for some $\eps < \theta$, for every $\lambda =
  \lambda^\partial$ for\footnote{What if we ask for a fixed $\varphi$,
    not depending on $\lambda$?  This makes $(c) \Rightarrow (d)$ problematic.}
 some $\varphi = \varphi(\bar x_{[\eps]},\bar y_{[\zeta]}) 
\in \bbL_{\theta,\theta}$
 there are $A \subseteq M \models T,|A| = \lambda$ 
such that $\bfS^\eps_\varphi(A,M)$ has
cardinality $> \lambda$
\sn
\item[$(e)$]  like (c) but for some $\lambda = \lambda^\partial$
\sn
\item[$(f)$]  like (d) but for some $\lambda = \lambda^\partial$
\sn
\item[$(g)$]  $T$ is definably unstable
\sn
\item[$(h)$]  there are $\eps < \theta,M \models T,\varphi =
\varphi(\bar x_{[\eps]},\bar y_{[\zeta]}) \in
\bbL_{\theta,\theta}(\tau_T)$ and $\LL (\bar b_{\alpha,0},\bar
b_{\alpha,1},\bar c_\alpha):\alpha < \theta\RR$ such that:
\sn
\begin{enumerate}
\item[$\bullet$]  $\bar b_{\alpha,0},\bar b_{\alpha,1} \in {}^\zeta\! M$
  and $\bar c_\alpha \in {}^\eps\! M$
\sn
\item[$\bullet$]  $\tp(\bar b_{\alpha,0},\cup\{\bar b_{\beta,0},
\bar b_{\beta,1},\bar c_\beta:\beta < \alpha\},M) = 
\tp(\bar b_{\alpha,1},\cup\{\bar b_{\beta,0},\bar b_{\beta,1},\bar
  c_\beta:\beta < \alpha\},M)$
\sn
\item[$\bullet$]  $\{\varphi(\bar x_\eps,\bar
b_{\beta,1}),\neg \varphi(\bar x_\eps,\bar b_{\beta,0}):
\beta < \alpha\}$ is realized by $\bar c_\alpha$ in $M$
\end{enumerate}
\sn
\item[$(i)$]  $T$ is 2-unstable, see \ref{y2}(6)
\sn
\item[$(j)$]  $T$ is 1-unstable, see \ref{y2}(1).
\end{enumerate}
\mn
2) $T$ is 3-unstable $\Rightarrow T$ is definably unstable.

\sn
3) $T$ is 1-unstable \underline{iff}  $T$ is $\{(\lambda,<)\}$-unstable for every
(equivalently some) $\lambda \ge \theta$.

\sn
4) $T$ is 5-unstable \underline{iff}  $T$ is $\{I\}$-unstable for every linearly
ordered $I$.

\sn
5) $T$ is 2-unstable \underline{iff}  for every $\eps,\zeta < \theta$ it
is $\eps \times \zeta^*$-unstable.  

\sn
6) In Definition \ref{y2}(1), we can use $\bar a_\alpha = \bar
b_\alpha$ so $\eps = \zeta$.
\end{claim}

\begin{PROOF}{\ref{y5}}  
1) \underline{$(a) \Rightarrow (b)$}

Obvious; by clause (a) there is $\varphi = \varphi(\bar x,\bar y) \in
\bbL_{\theta,\theta}(\tau_T)$ which witnesses $T$ is $({}^{\theta
  >}2,<_{\lex})$-unstable, so there is a model $M$ of $T$ and $\bar
a_\eta \in {}^{\lh(\eta)}\!M$ for $\eta \in {}^{\theta >}2$ such that
$M \models ``\varphi[\bar a_\eta,\bar a_\nu]"$ \underline{iff}  $(\eta <_{\lex} \nu)$
for every $\eta,\nu \in {}^{\theta >}2$.  
Let $\bar y = \bar y_{[\zeta]},\bar y' = \bar y_{[\zeta + \zeta]}$ and 
let $\varphi' = \varphi'(\bar x,\bar y')$ be 
$(\varphi(\bar x,\bar y' \rest [0,\zeta)) \equiv \varphi(\bar x,\bar
  y' \rest [\zeta,\zeta + \zeta))$, easily $\varphi'$ witnesses $T$ is
 4-unstable as witnessed by $\LL \bar b_\eta:\eta \in {}^{\theta
   >}2\RR$ where $\bar b_\eta = \bar a_{\eta \caret \LL 0 
\RR} \caret \bar a_{\eta \caret \LL 1 \RR}$.
\medskip

\noindent
\underline{$(b) \Rightarrow (c)$}

Let $\varphi(\bar x_{[\eps]},\bar y_{[\zeta]})$ be as in
\ref{y2}(2).  Note that
\mn
\begin{enumerate}
\item[$(*)$]  in Definition \ref{y2}(2), without loss of generality
there are $\bar c_\eta \in {}^\eps\! M$ for $\eta \in
{}^\theta 2$ realizing $p_\eta(\bar x_{[\eps]})$.
\end{enumerate}
\mn
[Why?  There is a $\theta$-complete uniform ultrafilter $D$ on
  $\theta$ hence in $M^\theta/D$ there are such $\bar c_\eta$-s.]  

So by compactness for $\bbL_{\theta,\theta}$, for every $\lambda$
there are $M_\lambda \models T$ and $\bar a^\lambda_\nu \in
{}^\zeta(M_\lambda)$ for $\nu \in {}^{\lambda >} 2$ and $\bar
c^\lambda_\eta \in {}^\eps(M_\lambda)$ for $\eta \in {}^\lambda
2$ such that $M_\lambda \models \varphi[\bar c^\lambda_\eta,\bar
  a^\lambda_\nu]^{\iif(\eta(\lh(\nu)))}$ when $\nu \triangleleft
\eta \in {}^\lambda 2$.

For any cardinal $\lambda$ let $\mu = \min\{\mu:2^\mu > \lambda\}$
hence $\mu \le \lambda \wedge (\forall \partial < \mu)(2^\partial \le
\lambda)$ and so $2^{< \mu} \le \lambda$ hence $\mu \le \lambda$, 
let $A = \bigcup\{\bar a^\mu_\nu:\nu \in {}^{\mu >}2\}$, so 
$A \subseteq M_\mu$ has cardinality $\le 2^{<\mu} + \theta 
\le \lambda$ and $\bfS^\eps(A,M_\mu)$ has cardinality
$\ge |\{\tp(\bar c^\mu_\eta,A,M_\mu):\eta \in {}^\mu 2\}| 
\ge 2^\mu > \lambda$. 
\medskip

\noindent
\underline{$(c) \Rightarrow (d)$}

It suffices to prove $\neg(d) \Rightarrow \neg(c)$.  So assume that
$\neg(d)$ holds and note that clearly the
set $\{\varphi(\bar x_{[\eps]},\bar y_{[\zeta]}) \in
\bbL(\tau_T):\eps,\zeta < \theta\}$ has cardinality $\partial$
recalling $\partial = (|T|^{< \theta} + \theta)$.

Hence if $A \subseteq M \models T$ and $|A| \le \lambda$ then
\mn
\begin{itemize}
\item  $|\bfS^\eps(A,M)| \le
\Pi\{|\bfS^\eps_\varphi(A,M)|:\varphi = \varphi(\bar
x_{[\eps]},\bar y_{[\zeta]}) \in
\bbL_{\theta,\theta}(\tau_T)\} \le (\sup\{|\bfS^\eps_\varphi(A)|:
\varphi=\varphi(\bar x_{[\eps]},\bar y_{[\zeta]})
\in \bbL_{\theta,\theta}(\tau_T)\})^\partial \le \lambda^\partial =
\lambda$.
\end{itemize}
\mn
[Why?  First inequality by the definitions of $\bfS^\eps(-),\bfS^\eps_\varphi(-)$, second inequality because the number of relevant
$\varphi$-s is $\le \partial$, third inequality by the present
assumption $\neg(d)$; the last inequality by the meaning of $\neg(d)$.
But the deduced inequality means $\neg(c)$.]
\medskip

\noindent
\underline{$(c) \Rightarrow (e)$}

Easy as there are $\lambda = \lambda^\partial$.
\medskip

\noindent
\underline{$(d) \Rightarrow (f)$}

As there are cardinals $\lambda$ such that $\lambda = \lambda^\partial$.
\medskip

\noindent
\underline{$(e) \Rightarrow (f)$}

As in $(c) \Rightarrow (d)$.
\medskip

\noindent
\underline{$(f) \Rightarrow (g)$}

Clearly $\neg(g) \Rightarrow \neg(f)$ holds by counting.
\medskip

\noindent
\underline{$(g) \Rightarrow (h)$}

So by compactness for $\bbL_{\theta,\theta}$ 
for some $\eps < \theta$ and $M \models T$ and $p \in \bfS^\eps(M)$ and $\varphi = \varphi(\bar x_{[\eps]},\bar
y_{[\zeta]})$ there are no $\psi(\bar y_{[\zeta]},\bar z_{[\xi]})$ and
$\bar c \in {}^\xi M$ as in Definition \ref{y3}.  Again by compactness for
$\bbL_{\theta,\theta}$ without loss of generality $|\tau_T| < \theta$.

Recalling Definition \ref{y3}(2), for each $\kappa < \theta$ we try by induction on $\alpha < \kappa$ to
choose $\bar b^\kappa_{\alpha,0},\bar b^\kappa_{\alpha,1},\bar
c^\kappa_\alpha$ such that (recalling \ref{y2}(3)):
\mn
\begin{itemize}
    \item  $\bar b^\kappa_{\alpha,0},\bar b^\kappa_{\alpha,1} \in {}^\zeta\! M$ realize the same $\big\{\varphi^\perp (\bar x_{[\zeta]},\bar y_{[\eps]})\big\}$-type over 
    $$A^\kappa_\alpha \defeq \bigcup\{\bar b^\kappa_{\beta,0}, \bar b^\kappa_{\beta,1},\bar c^\kappa_\alpha : \beta < \alpha\}.$$

    \item   $\varphi(\bar x_{[\eps]},\bar
  b^\kappa_{\alpha,1}),\neg \varphi(\bar x_{[\eps]},\bar
b^\kappa_{\alpha,0}) \in p$
\sn
    \item  $\bar c^\kappa_\alpha$ realizes
$\{\varphi(\bar x_{[\eps]},\bar b^\kappa_{\beta,1}),
\neg \varphi(\bar x_{[\eps]},
\bar b^\kappa_{\beta,0}):\beta \le \alpha\}$.
\end{itemize}
\mn
\underline{Case 1}:  For every $\kappa$ we succeed to carry the induction.

Let $\bar c^\kappa \in {}^\eps\! M$ realize $\{\varphi(\bar
x_{[\eps]},\bar b^\kappa_{\alpha,1}) \wedge \neg \varphi(\bar
x_{[\eps]},\bar b^\kappa_{\alpha,0}):\alpha < \kappa\}$.  By
compactness for $\bbL_{\theta,\theta}$ we can get clause (h).
\medskip

\noindent
\underline{Case 2}:  For some $\kappa$ and $\alpha < \kappa$, we
cannot choose $\bar b^\kappa_{\alpha,0},b^\kappa_{\alpha,1}$ (but have
chosen $\LL \bar b^\kappa_{\beta,\ell}:\beta < \alpha,\ell <
2\RR$).

We can find $\psi$ contradicting our choice of $M,\varphi,p$.
\medskip

\noindent
\underline{$(h) \Rightarrow (j)$}

Let $\varphi(\bar x_{[\eps]},\bar y_{[\zeta]}),M,\bar b_{\alpha,0},
\bar b_{\alpha,1},\bar c_\alpha(\alpha < \theta)$ be as in clause (h)
and let $\varphi'$ be as in the proof of $(a) \Rightarrow (b)$.  Now
$\varphi',\LL(\bar c_\alpha,\bar b_{\alpha,0} \caret \bar
b_{\alpha,1}):\alpha < \theta\RR$ are as required in clause (j)
 because for $\alpha,\beta < \theta$ we have 
$M \models ``\varphi[\bar c_\alpha, 
\bar b_{\beta,0}] \equiv \varphi[\bar c_\alpha,\bar b_{\beta,1}]"$ iff
  $\beta > \alpha$.
\medskip

\noindent
\underline{$(j) \Rightarrow (i)$}

Let $I = \theta \times \theta^*$, i.e. $\{(\alpha,\beta):\alpha,\beta
< \theta\}$ ordered by $(\alpha_1,\beta_1) < (\alpha_2,\beta_2)$ iff
$\alpha_1 < \alpha_2$ or $\alpha_1 = \alpha_2 \wedge \beta_1 > \beta_2$.  

Let $\varphi(\bar x_{[\eps]},\bar y_{[\zeta]})$ witness $T$ is
1-unstable and $M,\LL (\bar a_\alpha,\bar b_\alpha):\alpha <
\theta\RR$ exemplify this.  Let $\bar x' = \bar x_{[\eps +
 \eps]},\bar y' = \bar y_{[\zeta + \zeta + \eps]}$ 
and for $\alpha,\beta < \theta$ let $\bar a'_{(\alpha,\beta)} = 
\bar a_\alpha \caret \bar a_\beta,\bar b'_{(\alpha,\beta)} = \bar
b_{\alpha} \caret \bar b_{\beta +1} \caret \bar a_\alpha$
and let $\varphi'(\bar x',\bar y')$ say $\varphi(\bar x' \rest
\eps,\bar y' \rest \zeta)$ or $(\bar x' 
\rest \eps = \bar y' \rest [\zeta + \zeta,\zeta + \zeta +
  \eps) \wedge \neg \varphi(\bar x' \rest 
[\eps,\eps + \eps),\bar y'[\zeta,\zeta + \zeta))$. 

Now $\varphi',M,\LL (\bar a'_\alpha,\bar b'_\alpha):\alpha <
 \theta\RR$ are as required in Definition \ref{y2}(3) by part (5)
 proved below. 
\medskip

\noindent
\underline{$(i) \Rightarrow (j)$}

Trivially.

\noindent
2) Note that ``3-unstable $\Rightarrow$ definably unstable" holds by 
recalling the Definitions \ref{x11}(1), \ref{y2}(5), \ref{y3}(1).

\noindent
3) Easy, too.

\noindent
4) First, the implication $\Rightarrow$ holds by ``$\theta$ is
compact" because every linear order $I$ is embeddable into $({}^\alpha
2,<_{\lex})$ for some ordinal $\alpha$.  Second, the 
implication $\Leftarrow$ is trivial.

\noindent
5) First, the implication $\Rightarrow$ holds as $\theta$ is a compact
cardinal.  Second, the implication $\Leftarrow$ is trivial.

\noindent
6) Easy, too, using enough dummy variables; i.e. let $\bar a' = \bar
a_\alpha \caret \bar b_\alpha$ and $\varphi''(\bar x_{[\eps +
  \zeta]},\bar y_{[\eps + \zeta]}) \defeq \varphi(\bar x_{[\eps +
  \zeta]} \rest [0,\eps),\bar y_{[\eps +  \zeta]} \rest 
[\eps,\eps + \zeta))$.
\end{PROOF}

\begin{conclusion}  
\label{y10}
1) Assume $T \subseteq \bbL_{\theta,\aleph_0}$ is (complete and)
$(\varphi(\bar x_{[\eps]},\bar y_{[\zeta]},\{\theta\})$-unstable for
some $\eps < \theta$ and $\varphi(\bar x_{[\eps]},\bar y_{[\eps]})
\in \bbL_{\theta,\aleph_0}$.

For every $\lambda \ge \theta^+ + |T|$, there are
$M_\alpha \in \Mod_T$ for $\alpha < 2^\lambda$ which are pairwise
non-isomorphic each of cardinality $\lambda$.

\noindent
2) If $T \subseteq \bbL_{\theta,\theta}$ is strongly 3-unstable and
$\lambda = \lambda^{< \theta} \ge \theta^+ + |T|$, \underline{then}  the
conclusion of part (1) holds.
\end{conclusion}

\begin{PROOF}{\ref{y10}}  
Follows by \cite[\S3]{Sh:E59} (which improve \cite[Ch.III]{Sh:300}) but we
explain the background.  By \cite[Ch.VIII]{Sh:c}, if $T \subseteq T_1$
are complete first order and $\lambda \ge |T_1| + \aleph_1$ and $T$
unstable \underline{then}  there are models $M_\alpha$ of $T$ of cardinality
$\lambda$ for $\alpha < 2^\lambda$, pairwise non-isomorphic each from
$\PC(T_1,T)$, i.e. each $M_\alpha$ can be expanded to a model of
$T_1$.  This was done mainly using E.M. models, i.e. for some $T_2
\supseteq T_1$ of cardinality $\le \lambda$ with Skolem functions each
$M_\alpha$ can be expanded to a model $N_\alpha$ of $T_2$ which is
generated by $\{\bar a^\alpha_t:t \in I_\alpha\},I_\alpha$ a linear
order $\bar{\bfa}_\alpha = \LL \bar a^\alpha_t:t \in
I_\alpha\RR$ is an indiscernible sequence in $N_\alpha$ and for
some $\varphi(\bar x,\bar y) \in \bbL(\tau_T),N_\alpha \models
\varphi[\bar a^\alpha_s,\bar a^\alpha_s]$ iff $s <_{I_\alpha} t$.

Now \cite[\S3]{Sh:E59} improve it by just requiring $\LL
(N_\alpha,\bar{\bfa}_\alpha):\alpha < 2^\lambda\RR$ and $\eps
< \theta,\varphi = \varphi(\bar x_{[\eps]},\bar y_{[\eps]}) \in
\bbL_{\theta,\aleph_0}(\tau_T)$ to have some of the properties of such
E.M. models (called there ``being $\kappa$-skeleton like").

This means here just (where $\lambda$ is regular for transparency):
\mn
\begin{enumerate}
\item[$(*)$]  $(A) \Rightarrow (B)$  where:
\sn
\begin{enumerate}
\item[(A)]
\begin{enumerate}
\item[(a)]  $\bar{\bfa}_\alpha = \LL \bar a_{\alpha,s}:s \in
  I_\alpha\RR,\bar{\bfb}_\alpha = \LL \bar b_{\alpha,s}:s \in
  I_\alpha\RR,\bar b_{\alpha,s} = \bar a_{\alpha,s} \in
  {}^\eps(M_\alpha),\zeta = \eps,M_\alpha,\varphi(\bar
  x_{[\eps]},\bar y_{[\eps]}),\varphi(\bar x_{[\eps]},\bar
  y_{[\eps]})$ are as in Definition \ref{y2}(4)
\sn
\item[(b)]  $I_\alpha = \sum\limits_{i < \lambda}
  I_{\alpha,i},S_\alpha \subseteq \lambda,I_{\alpha,\eps}$ is
  isomorphic to $(\theta,>)$ if $\eps \in S_\alpha$ and to
  $(\theta^+,>)$ if $\eps \in \lambda \setminus S_\alpha$
\end{enumerate}
\sn
\item[(B)]  $\{M_\alpha/\cong:\alpha < 2^\lambda\}$ has cardinality
  $2^\lambda$.
\end{enumerate}
\end{enumerate}
\mn
Why there are models as in (A)?  For part (2) by Definition
\ref{y2}(5) and see \ref{y2}(4).  For part (1) by the definition on
E.M. nodes.  Note that in \cite[\S3]{Sh:E59} we first deal with the
case $\eps$ is finite, but we are assuming $\lambda = \lambda^{<
  \theta}$ hence allowing $\eps \in [\omega,\theta)$ cause no
problem, see \cite[Th.3.28,pg.48,L3c.16]{Sh:E59}.
\end{PROOF}

\begin{question}  
\label{y12}
1) Can we add in \ref{y10} 
``pairwise not $\bbL_{\infty,\theta^+}$-equivalent"?

\noindent
2) Does the logic $\cL$ have interpolation when 
$\bbL_{\theta,\aleph_0} \subseteq \cL \subseteq \bbL_{\theta,\theta}$
and $\cL$ is defined by: $\psi \in \cL(\tau)$ \underline{iff}  $\psi
   \in \bbL_{\theta,\theta}(\tau)$ and for $\bft \in \{\yes,\no\}$
 the class of models of $\psi^{\bft}$ is closed under 
$M^I_D|\cE$ when $(I,D,\cE)$ is $(\theta,\aleph_0)$-complete, see
 Definition \ref{x16}.
\end{question}

\noindent
Now recall stability implies the existence of convergence
sub-sequences, specifically:
\begin{claim}  
\label{y15}
Assume $|T| \ge \lambda = \cf(\lambda)$ and $\mu < \lambda \Rightarrow
(\mu^{|T|})^+ + \theta < \lambda,|T|^{< \theta} < \partial =
\cf(\partial) < \lambda$.  If $T$ is 1-stable,
$\eps < \theta,M$ is a model of $T$ and $\bar a_\alpha \in
{}^\eps\! M$ for $\alpha < \lambda$ \underline{then}  for some stationary
$S \subseteq S^\lambda_\partial$ the sequence $\LL \bar
a_\alpha:\alpha \in S\RR$ is $(< \omega)$-indiscernible and strongly
$\bbL_{\theta,\theta}$-convergent in $M$, see Definition \ref{x11}(1).
\end{claim}

\begin{PROOF}{\ref{y15}}
By \cite{Sh:300a} but we explain the background.  First, we may find a
$\prec_{\bbL_{\theta,\theta}}$ increasing sequence $\LL
M_\alpha:\alpha \le \lambda\RR$ such that $M_\alpha
\prec_{\bbL_{\theta,\theta}} M,\|M_\alpha\| \le |\alpha|^{< \theta} + \theta + |T|^{< \theta}$ and
$\bar a_\alpha \subseteq M_{\alpha +1}$.

Second, for each $\alpha \in S_0 \defeq \{\delta < \lambda:\cf(\delta)
= \partial\}$ we can find $B_\alpha \subseteq M_\alpha$ of cardinality
$\le |T|^{< \theta}$ such that $\tp_{\bbL_{\theta,\theta}}(\bar
a_\alpha,M_\alpha,M)$ is definable over $B_\alpha$.

Third, by Fodour lemma there is a stationary $S_1 \subseteq S$
such that $\LL B_\delta:\delta \in S_1\RR$ is constantly $B$,
and even the definition scheme is the same.  We then prove $\LL
\bar a_\alpha:\alpha \in S_1\RR$ is $n$-indiscernible by induction
on $n$ (as there).

Lastly, for proving convergence, we fix $\bar b \in {}^{\theta >}\!M$
and use ``$\tp_{\bbL_{\theta,\theta}}(\bar b,M_\lambda,M)$ is definable.
\end{PROOF}

\noindent
The experience with first order classes says categoricity even for
PC-classes (see below) implies stability (also 
$\triangleleft_{\lambda,\theta}$-minimality) \underline{however} this
is not so here 
(where on $\blacktriangleleft_{\lambda,\theta}$, see Definition \ref{a9})
hence we now consider some examples (see also \ref{a22}).  In the rest
of this section we prove this and give other examples.
\begin{claim}
\label{y41}
$T$ being 1-unstable does not imply $T$ being definably unstable, and
does not imply satisfying \ref{y5}(h).
\end{claim}

\begin{PROOF}{\ref{y41}}
Let $M = (\theta,<)$ and $T = \Th_{\bbL_{\theta,\theta}}(M)$; 
clearly $T$ is 1-unstable and is definably stable.  
As for \ref{y5}(h), toward contradiction assume
$N \models T$ and $\varphi = \varphi(\bar x_{[\eps]},\bar
y_{[\zeta]}),\LL (\bar a_\alpha,\bar b_\alpha,\bar
c_\alpha):\alpha < \theta\RR$ are as in clause (h) of \ref{y5}.
As $\theta$ is a compact cardinal without loss of generality $\LL \bar a_\alpha
\caret \bar b_\alpha \caret \bar c_\alpha:\alpha < \theta\RR$
is an indiscernible sequence in $M$, i.e. $n$-indiscernible for
every $n$.  Now check.
\end{PROOF}

\begin{thesis}  
\label{y43}
A big difference with the first order, that is the 
$\theta = \aleph_0$ case, is:
\mn
\begin{enumerate}
\item[$(a)$]  long linear orders does not contradict categoricity, in
  particular see \ref{y40} below 
\sn
\item[$(b)$]  consider interpreting for 
$\partial < \theta$, a group isomorphic to
 the Abelian group $(\{\eta \in {}^A 2:
(\exists^{< \partial}a \in A)(\eta(a) = 1)\},\triangle)$ where $\triangle$
is the symmetric difference; it appears ``for
 free" (formally\footnote{Why?  E.g. for a model $M$ let
\mn
\begin{enumerate}
\item[$\bullet$]  the set of elements in $\varphi(M)$ where $\varphi =
  \varphi(\bar x_{[\omega]})$ says: $\bigwedge\limits_{n \ne m} (x_{2n} \ne
  x_{2n+1} \wedge x_{2m} \ne x_{2m+1} \rightarrow x_{2n} \ne x_{2m})$,
  let $\Rang^*(\bar x_{2n}) = \{x_{2n}:x_{2n} = x_{2n+1}\}$
\sn
\item[$\bullet$]  the congruence $\varphi_{\eq}(\bar x_{[\omega]},\bar
  y_{[\omega]})$ says $\Rang^*(\bar x_{[\omega]}) = \Rang(\bar
  y_{[\omega]})$
\sn
\item[$\bullet$]  $\varphi_{\muult}(\bar x_{[\omega]},\bar
  y_{[\omega]},\bar z_{[\omega]}) = \Rang^*(\bar x_{[\omega]}) \triangle
  \Rang^*(\bar y_{[\omega]}) = \Rang^*(\bar z_{[\omega]})$.
\end{enumerate}
\mn
For clause (c) of \ref{y43} - more cumbersome.}, if we allow equality
 for the group being just a congruence relation)
\sn
\item[$(c)$]  similarly for the group generated by $\{x_a:a \in A\}$ freely.
\end{enumerate}
\end{thesis}

\begin{example}  
\label{y40}
\noindent
1) There are $T$ and $T_1$ such that:
\mn
\begin{enumerate}
\item[$(a)$]   $T \subseteq \bbL_{\theta,\theta}(\{<\})$ is complete
\sn
\item[$(b)$]  $T_1 \subseteq \bbL_{\theta,\theta}(\tau_1)$ is complete,
$\tau_1$ finite and $<$ belongs to $\tau_1$
\sn
\item[$(c)$]   $T_1 \supseteq T$
\sn
\item[$(d)$]  models of $T$ are dense linear orders
\sn
\item[$(e)$]  $\PC(T,T_1)$ is categorical in every $\lambda \ge
  \theta$, recalling
\sn
\item[${{}}$]  $\bullet \quad \PC(T,T_1) = \{M_1 \rest \tau_T:M_1 \in
  \Mod_{T_1}\}$ 
\sn
\item[$(f)$]  $T$ is 1-unstable
\sn
\item[$(g)$]  $T$ is definably stable.
\end{enumerate}
\mn
2) Moreover $T = \Th_{\bbL_{\theta,\theta}}(N)$ where:
\mn
\begin{enumerate}
\item[(a)]
\begin{enumerate}
\item[($\alpha$)]  $N$ is a dense linear order
\sn
\item[($\beta$)]  $N$ is of cardinality $\theta$
\end{enumerate}
\sn
\item[(b)]
\begin{enumerate}
\item[($\alpha$)]  $N$ is the union of $\aleph_0$ well ordered sets
\sn
\item[($\beta$)]  $N$ has cofinality $\aleph_0$, also its inverse
\end{enumerate}
\sn
\item[(c)]
\begin{enumerate}
\item[($\alpha$)]    if $\sigma$ is regular uncountable, any increasing
  sequence of length $\sigma$ has no lub
\sn
\item[($\beta$)]   if $s \in N$ then $N_{<t} = N \rest \{s:s <_I t\}$
  has cofinality $\aleph_0$ and $N_{>t} = N \rest \{s:t <_I s\}$ has
co-initiality $\aleph_0$
\end{enumerate}
\sn
\item[(d)]  any two intervals of $N$ are isomorphic (note: $T$ cannot
  say this but $T_1$ can).
\end{enumerate}
\mn
3) Moreover $T_1$ extends $T$ and just says in addition only that
every two intervals of $N$ are isomorphic.
\end{example}

\begin{remark}
\label{y40g}
1) See \cite[\S2]{Sh:E62} as explained below.

\noindent
2) Hausdorff has introduced and investigated the class of scattered
linear orders.  Galvin and Laver, see \cite{Lv71}
investigate the class $\cM$ of 
linear orders which are a countable union of scattered
linear orders.  They were interested in linear orders up to
embeddability inside the class $\cM = \bigcup\{\cM_{\lambda,\mu_1,\mu_2}:
\mu_1,\mu_2$ are regular uncountable such that
$\lambda^+ = \mu_1 + \mu_2\}$ where $\cM_{\lambda,\mu_1,\mu_2}$ is the class
of linear orders from $\cM$ of cardinality $\lambda$ with 
no increasing sequences of
length $\mu_1$ and no decreasing sequences of length $\mu_2$.
Galvin defined $\cM_{\lambda,\mu_1,\mu_2}$ and prove existence of a
universal member.

Laver, solving a long standing conjecture of Fra\"isse, and using the
theory of better quasi orders of Nash Williams \underline{prove} the
following.  The class $\cM$ is well quasi ordered and even better 
quasi order under embeddability; this answers affirmatively 
Fra\"isse's conjecture which says that $\cM_{\aleph_1,\aleph_1} =$ the
class of countable linear orders, is well ordered.  
So categoricity (\ref{y40}(1)(e)) and
 clause (c) of \ref{y40}(2) were irrelevant there, the latter is crucial here 
for categoricity.  
In \cite[pp.308,309]{Sh:220}, this is continued being interested in
uniqueness.  We do more in \cite[\S2]{Sh:E62}.

\noindent
3) As requested we explain that in \cite[\S2]{Sh:E62}, we investigate
classes of $I^+$ of the form: a linear order, $I$ expanded by unary
relations $P^{I^+}_s(s \in S)$ such that $\LL P^{I^+}_s:s \in
S\RR$ is a partition of $I$ and if, e.g. $\LL t_i:i <
\kappa\RR$ is increasing with lub $t_\kappa,\kappa = \cf(\kappa) >
\aleph_0$ and $t_\kappa \in P^{I^+}_s$ then we know for a club
of $\delta < \kappa$, what is the co-initiality of $\{s \in I:(\forall
i < \delta)(t_i <_I s)\}$ and more.  It is proved there that 
under such restritions we get uniqueness for those expanded linear orders.
\end{remark}

\begin{PROOF}{\ref{y40}}  
We know (see \cite[\S2]{Sh:E62} and \ref{y40g} above)
\mn
\begin{enumerate}
\item[$(*)_1$]  there is a linear order $N$ satisfying Clauses
(a)-(d) of part (2)
\sn
\item[$(*)_2$]
\begin{enumerate}
\item[(a)]   choose $N_*$ as in $(*)_1$
\sn
\item[(b)]  let $T = \Th_{\bbL_{\theta,\theta}}(N_*)$
\sn
\item[(c)]   let $T_1$ be $T \cup \{\psi\}$, where $\psi$
  says that: if $x_1 < y_1,x_2 < y_2$ then 
$z \mapsto F(z,x_1,y_1,x_2,y_2)$ is an isomorphism from the interval
$(x_1,y_1)$ onto the interval $(x_2,y_2)$ for the linear order
\sn
\item[(d)]  note that the theory $T_1$ is consistent as we can expand 
$N_*$ to a model of $T_1$
\end{enumerate}
\sn
\item[$(*)_3$]
\begin{enumerate}
\item[(a)]   if $N$ is a linear order failing sub-clause ($\alpha$)
  of (b) of \ref{y40}(2) \underline{then}  there is $N_1 \subseteq N$ 
of cardinality $< \theta$ failing it, hence $N$ is not a model of $T$
\sn
\item[(b)]  similarly for $(b)(\beta),(c)(\beta)$ and even
  $(c)(\alpha)$ for $\sigma < \theta$.
\end{enumerate}
\end{enumerate}
\mn
[Why?  By $\theta$ being a compact cardinal.]

So easily
\mn
\begin{enumerate}
\item[$(*)_4$]
\begin{enumerate}
\item[(a)]  if $M$ is a model of $T$ then $M$
 satisfies Clauses $(a)(\alpha),(b),(c)$ of \ref{y40}(2)
\sn
\item[(b)]   if $M \in \PC(T,T_1)$, i.e. $M = M_1 \rest
 \{<\}$ where $M_1 \models T_1$ \underline{then}  $M$ satisfies 
Clauses $(a)(\alpha),(b),(c),(d)$ of \ref{y40}(2).
\end{enumerate}
\end{enumerate}
\mn
[Why?  Mainly by $(*)_3$, e.g. why $M$ satisfies clause $(c)(\alpha)$
of \ref{y40}(2)? let $\bar a =
  \LL a_\alpha:\alpha < \partial \RR$ be increasing, $\partial$
regular uncountable and we shall prove it has no lub.  If $\partial <
\theta$ this is said in $T$.  If $\partial \ge \theta$ or just $\partial
\ge \aleph_1$, then $\bar a$ is bounded (see
\ref{y40}(2)(b)($\beta$)) so there is a decreasing $b = \LL
  b_\beta:\beta < \kappa\RR$ such that $(\bar a,\bar b)$ is a
  pre-cut of $M$, see \ref{x13}(4) and $\kappa$ is 1 or a regular
  cardinal.  Now by \ref{y40}(2)$(b)(\alpha)$ necessarily $\kappa 
= \aleph_0$ or $\kappa=1$; but by $M \models T$ recalling
\ref{y40}$(c)(\beta),\kappa = 1$ is impossible.]

Also
\mn
\begin{enumerate}
\item[$(*)_5$]  $\PC(T,T_1)$ is categorical in every
  $\lambda \ge \theta$.
\end{enumerate}
\mn
[Why?  By \cite[\S2]{Sh:E62} and see \ref{y40g}(3).]

So $T$ satisfies all the clauses of \ref{y40}(1), e.g. we shall prove
that $T$ is definably stable; toward this assume
\mn
\begin{enumerate}
\item[$(*)_{6.1}$]   $M \prec_{\bbL_{\theta,\theta}} N$ are models of
  $T$ and we should prove that for $\bar a \in {}^{\theta >}
M,\tp_{\bbL_{\theta,\theta}}(\bar a,M,N)$ is definable (in $M$).
\end{enumerate}
\mn
Toward this for $a \in N \setminus M$ clearly $M_{>a} \defeq
\{b \in M:a <^N b\}$ has co-initiality 1 or $\aleph_0$ so let
$\bar b_{a,1}$ list a countable subset of $M_{>a}$ unbounded from
below in $M_{> a}$.

Let $M_{<a} = \{b \in M:b <^N a\}$ and let $\bar b_{a,2}$ be a sequence of
elements of $M_{<a}$ of length $< \theta$ which is unbounded in $N_{<
  a} \cap M$ if possible, empty otherwise.  Letting $\bar b = \bar b_{a,1}
\caret \bar b_{a,2}$ clearly it is a sequence of elements of $M$ of
length $< \theta$ (but actually $\bar b_2$ is not necessary).

So clearly it suffices to prove:
\mn
\begin{enumerate}
\item[$(*)_{6.2}$]  if $\bar a \in {}^{\theta >}\!N$ and $\bar b \in
{}^{\theta >}\!M$ includes $\bar b_{a_\eps}$ (or just 
$\bar b_{a_\eps,1}$) for every $\eps < \lh(\bar a)$ then
  $\tp_{\bbL_{\theta,\theta}}(\bar a,M,N)$ is definable over $\bar b$.
\end{enumerate}
\mn
For this it suffices to prove:
\mn
\begin{enumerate}
\item[$(*)_{6.3}$]   Assume $\partial \le \theta$ is regular and
e.g. inaccessible, $\eps < \partial$ and $\bar a_1,\bar a_2 \in
{}^\eps\! N$.  The following are equivalent: 
\sn
\begin{enumerate}
\item[(a)]  $\tp_{\bbL_{\theta,\tau}}(\bar a_1,M,N) =
  \tp_{\bbL_{\theta,\theta}}(\bar a_2,M,N)$
\sn
\item[(b)]
\begin{enumerate}
\item[($\alpha$)]   if $\xi,\zeta < \eps$
  then $a_{1,\xi} <_M a_{1,\zeta} \Leftrightarrow a_{2,\xi} <
  a_{2,\zeta}$ (in $M$)
\sn
\item[($\beta$)]  if $u \subseteq \eps$ \underline{then} 
the cofinalities of $\bigcap\limits_{\zeta \in u} 
M_{< a_{1,\zeta}},\bigcap\limits_{\zeta \in u} M_{< a_{2,\zeta}}$ are
 equal or are both $\ge \partial$
\sn
\item[($\gamma$)]   if $u \subseteq \eps$ \underline{then}  the
co-initialities of $\bigcap\limits_{\zeta \in u} M_{> a_{1,\zeta}},
\bigcap\limits_{\zeta \in u} M_{> a_{2,\zeta}}$ are equal or 
are both $\ge \partial$.
\end{enumerate}
\end{enumerate}
\end{enumerate}
\mn
This is easy to check.
\end{PROOF}

\begin{example}  
\label{y45}
0) $\Th_{\bbL_{\theta,\theta}}(\theta,<)$ is 1-unstable, definably
   stable.

\noindent
1) Let $T_2 = \Th(N),N$ is the linear order $\theta \times 
(\theta +1)^*$ ordered lexicographically expanded by $P^N = \theta
\times \{\theta +1\}$.

Then:
\mn
\begin{enumerate}
\item[$(a)$]  $T_2$ is 2-unstable as exemplified by a formula $\varphi
  = \varphi(x,y)$ but $T_2$ is 3-stable and stable
as well as 4-stable and 5-stable
\sn
\item[$(b)$]  $M$ is a model of $T_2$ when $M$ is $\sum\limits_{i <
  \delta} M_i,\delta$ an ordinal of cofinality $\ge \theta$ and each
  $M_i$ is isomorphic to $\delta_i + 1,\delta_i$ an ordinal of cofinality
  $\ge \theta$.
\end{enumerate}
\mn
2) Let $T_3 = \Th_{\bbL_{\theta,\theta}}(N),N$ is the linear order
$\theta \times \theta^*$.  

Then
\mn
\begin{enumerate}
\item[$(a)$]  $T_3$ is 3-unstable but stable hence 4-stable and 5-stable
\sn
\item[$(b)$]  like \ref{y45}(1)(b) but $M_i \cong \delta_i$.
\end{enumerate}
\mn
3) Let $T_4 = \Th_{\bbL_{\theta,\theta}}({}^{\theta
  >}2,\triangleleft)$
\mn
\begin{enumerate}
\item[$(a)$]  $T_4$ is 4-unstable but 5-stable and 3-stable
\sn
\item[$(b)$]  $M$ is a model of $T$ iff it is isomorphic to
  $(\clT,\triangleleft)$ where for some ordinal $\alpha$ of cofinality
  $\ge \theta,\clT$ is a subset of ${}^{\alpha >} 2$, closed 
under initial segments, $\eta \in \clT \Rightarrow \eta \caret 
\LL 0 \RR \in \clT \wedge \eta
  \caret \LL 1 \RR \in \clT$ and $\clT$ is closed under
  increasing unions of length $< \theta$.
\end{enumerate}
\mn
4) Let $T_5$ be the $\bbL_{\theta,\theta}$-theory of any dense linear
order which is $\theta$-saturated in the first order sense (so with
neither first nor last element), can use also
$\Th_{\bbL_{\theta,\theta}}({}^{\theta >}2,<_{\lex})$
\mn
\begin{enumerate}
\item[$(a)$]  $T_5$ is $\iota$-unstable, for $\iota = 1,\dotsc,5$.
\end{enumerate}
\mn
5) Let $T_6 = \Th_{\bbL_{\theta,\theta}}(M)$ where $M = ({}^{\theta
  >}2,\triangleleft,P^M),P^M = \{\eta \caret \LL 1 \RR:\eta
\in {}^{\theta >}2\}$ so $\tau_M = \{<,P\}$ so $<,P$ are two-place,
one-place predicates respectively, \underline{then}  $T_6$ is 5-unstable but 3-stable.
\end{example}

\begin{PROOF}{\ref{y45}}
This proof almost always uses only $\theta = \cf(\theta) > \aleph_0$;
we shall mention when not.

\noindent
0) See the proof of \ref{y41}.

\noindent
1) Note that
\mn
\begin{enumerate}
\item[$(*)_1$]  $(a) \quad$ if $(C_1,C_2)$ is a cut of $\theta \times
  (\theta +1)^*$, \underline{then}  the cofinality of $(C_1,C_2)$ is 

\hskip25pt  one of the following:
  $(0,1),(1,\theta),(1,\partial),(1,1),(\partial,1),(\theta,0)$ 

\hskip25pt with $\partial = \cf(\partial) < \theta$
\sn
\item[${{}}$]  $(b) \quad$ every one of those cofinalities appear.
\end{enumerate}
\mn
[Why?  By inspection.]
\mn
\begin{enumerate}
\item[$(*)_2$]   if $N$ is a model of $T_2$ and $(C_1,C_2)$ is a cut
  of $N$ \underline{then}  the cofinality of $(C_1,C_2)$ is one of the following:
$(0,1),(1,\lambda_1),(1,\partial),(1,1),(\partial,1),(\lambda_2,0)$
  with $\partial = \cf(\partial) < \theta,\lambda_1 = \cf(\lambda_1) \ge
  \theta$ and $\lambda_2 = \cf(\lambda_2) \ge \theta$.
\end{enumerate}
\mn
[Why?  Follows from $(*)_3$ which is proved below.]
\mn
\begin{enumerate}
\item[$(*)_3$]
\begin{enumerate}
\item[(a)]   let $\varphi_1(x,y)$ say: $x<y$ and there
is no $z \in (x,y]$ such that $P(z)$
\sn
\item[(b)]   let $\varphi_2(x,y) = \varphi_1(x,y) \vee
 \varphi_1(y,x) \vee x = y$
\sn
\item[(c)]  if $N \models T_2$ \underline{then}  $\varphi_2$
 defines an equivalence relation on $N$, each 
equivalence class A is $\bbL_{\theta,\theta}$-equivalent to $(\theta +1)^*$ 
\, $(\bbL_{\aleph_1,\aleph_1}$ suffice) hence $N \rest A$ is
anti-well (linearly) ordered, with a 
first element and last element and omitting the first element 
 of co-initiality $\ge \theta$
\sn
\item[(d)]  if $N \models T_2$ \underline{then}  the linear order
$\varphi^N_2$ is $\bbL_{\theta,\theta}$-elementarily equivalent 
to $\theta$.
\end{enumerate}
\end{enumerate}
\mn
[Why?  Should be clear.]

By $(*)_3$, Clause (b) of \ref{y45}(1) holds.  Now Clause (a) of
\ref{y45}(1) follows by checking Definition \ref{y2}.

\noindent
2) Similarly replacing $(\theta +1)^*$ by $\theta^*$.

\noindent
3) Let $\tau = \{<\}$, $M = ({}^{\theta >} 2,\triangleleft)$ a
$\tau$-model so ${<^M} = {\trianglelefteq} \rest {}^{\theta >}2$.
Clause (b) should be clear and anyhow we use just $\Rightarrow$.  
For Clause (a), $T_4$ being 4-unstable
 holds for the formula $\varphi = \varphi(x,y) = (y < x)$ by the 
definition of 4-unstable in \ref{y2}(2).  As being ``5-stable" is
   easier, we shall just prove ``$T_4$ is 3-stable".

For this we prove the following, using $\theta$ 
is a compact cardinal; clearly this
suffices; the $\varphi,\psi$ below are not related to Definition \ref{y2}(4):
\mn
\begin{enumerate}
    \item[$\boxplus$]  Assume $M \models T_4$ and $\delta_1,\delta_2$ are ordinals of cofinality $\ge \theta$, but 
    $\cf(\delta_1) \ne \cf(\delta_2)$ and $J = (\{1\} \times \delta_1) \cup (\{2\} \times \delta_2)$ ordered by 
    $$\alpha_1 < \beta_1 < \delta_1 \wedge \alpha_2 < \beta_2 < \delta_2 \Rightarrow (1,\alpha_1) < (1,\beta_1) < (2,\beta_2) < (2,\alpha_2)$$ 
    and $\varphi = \varphi(\bar x_{[\eps]},\bar y_{[\zeta]}) \in \bbL_{\theta,\theta}(\tau_M)$, $\bar a_s \in {}^\eps\! M$, $\bar b_s \in {}^\zeta\! M$ for $s \in J$ and $M \models \varphi[\bar a_s,\bar b_t]^{\iif(s<t)}$.  \underline{Then} for some 
    $\psi(\bar x,\bar z) \in \bbL_{\theta,\theta}(\tau_M)$ and 
    $\bar c$ from ${}^{\lh(\bar z)}\!M$ we have:
\sn
\begin{enumerate}
\item[$(a)$]  $\delta_1 = \sup\{\alpha_1 < \delta_1 : M \models
``\psi[\bar a_{(1,\alpha_1)},\bar c]"\}$
\sn
\item[$(b)$]  $\delta_2 = \sup\{\alpha_2 < \delta_2:M \models
``\neg\psi[\bar a_{(2,\alpha_2)},\bar c]"\}$.
\end{enumerate}
\end{enumerate}
\mn
Why?  For $\ell=1,2$ let $D_\ell$ be a $\theta$-complete ultrafilter
  on $\delta_\ell$ such that $\alpha < \delta_\ell \Rightarrow
[\alpha,\delta_\ell) \in D_\ell$.  As in \ref{y5}(6), without loss of generality $\bar a_s =
\bar b_s$ and by clause (b) of \ref{y45}(3), 
$M = (\clT,\triangleleft)$ where $\clT,\alpha$ are as there.

Let $\clT^+ = \clT \cup \{\eta \in {}^{\alpha >}2:\lh(\eta)$ is a
limit ordinal and $\beta < \lh(\eta) \Rightarrow \eta \rest \beta
\in \clT\}$, clearly $\eta \in \clT^+ \setminus \clT \Rightarrow
\cf(\lh(\eta)) \ge \theta$ using $T_3 =
\Th_{\bbL_{\theta,\theta}}(M)$.  For $s \in J$ let 
$\bar a_s = \LL a_{s,i}:i < \zeta\RR$ and for each 
$i < \zeta$ we choose $\eta^1_i,\eta^2_i \in \clT^+$ such that:
\mn
\begin{enumerate}
\item[$\bullet$]  $\eta^\ell_i = \bigcup\{\nu \in \clT:\{\alpha <
  \delta_\ell:\nu \trianglelefteq a_{(\ell,\alpha),i}\} \in D_\ell\}$.
\end{enumerate}
\mn
Let $u_\ell = \{\eps < \zeta:\{\alpha <
\delta_\ell:a_{(\ell,\alpha),\eps} = \eta^\ell_\eps\} 
\in D_\ell\}$ clearly
\mn
\begin{enumerate}
\item[$(*)_1$]  $\eps \in u_\ell \Rightarrow
  \eta^\ell_\eps \in \clT$
\sn
\item[$(*)_2$]  $u_\ell \ne \zeta$.
\end{enumerate}
\mn
[Why?  By $s,t \in J \Rightarrow M \models \varphi[\bar
a_s,\bar b_t]^{\iif(s<t)}$, see the statement of
$\boxplus$ hence $s \ne t \Rightarrow
\bar a_s \ne \bar a_t$ but $u_\ell = \zeta \Rightarrow
\bigwedge\limits_{\alpha,\beta < \delta_\ell} \bar a_{(\ell,\alpha)}
  = \bar a_{(\ell,\beta)}$.]

Now we prove $\boxplus$ by cases.
\medskip

\noindent
\underline{Case 1}:  $\eps \in u_1$ but $\eps \notin u_2
\vee (\eps \in u_2 \wedge \eta^1_\eps \ne
\eta^2_\eps)$.

Let $\psi(\bar x_{[\zeta]},\bar c) = (x_{[\eps]} = \eta^1_\eps)$
and check.
\medskip

\noindent
\underline{Case 2}:  $\eps \in u_2$ but $\eps \notin u_1
\vee (\eps \in u_1 \wedge \eta^1_\eps \ne \eta^2_\eps)$.

Let $\psi(\bar x_{[\zeta]},\bar c) = (x_{[\eps]} \ne
\eta^2_\eps)$ and check.
\medskip

\noindent
\underline{Case 3}:  $\eps < \zeta$, $\eps \notin u_1$, 
$\eps \notin u_2$ but $\eta^1_\eps \ne \eta^2_\eps$.

By symmetry without loss of generality $\lh(\eta^1_\eps) > \ell
g(\eta^2_\eps)$, let 
$\nu \in \clT$ be such that $\nu \triangleleft \eta^1_\eps$
but $\nu \ntrianglelefteq \eta^2_\eps$, clearly exists and 
let $\psi(\bar x_\zeta,\bar c) = (\nu \triangleleft x_\eps)$ and check.
\medskip

\noindent
\underline{Case 4}:  $\eps < \zeta$, $\eps \notin u_1 \cup
u_2$, $\eta^1_\eps = \eta^2_\eps$ but for
some $\nu \triangleleft \eta^1_\eps$ we have\\
$\delta_1 = \sup\{\alpha < \delta_1:\nu \triangleleft
a_{(\upharpoonleft,\alpha),\eps}\}$.

Let $\psi(\bar x_\zeta,\bar c) = (\nu \ntrianglelefteq
x_\eps)$.
\medskip

\noindent
\underline{Case 5}:  Like Case 4, for $\delta_2$.

Similarly.

Now if none of the cases above holds, then by $(*)_2$ there is
$\eps < \zeta$ such that $\eps \notin u_1$; by `not Case
2,' $\eps \notin u_2$; by `not Case 3,' $\eta^1_\eps =
\eta^2_\eps$; by `not Case 4,' $\cf(\lh(\eta^1_\eps)) =
\cf(\delta_1)$, and by `not Case 5,' $\cf(\lh(\eta^2_\eps)) =
\cf(\delta_2)$.  Together necessarily $\cf(\delta_1) =
\cf(\delta_2)$, contradicting an assumption.  

So $\boxplus$ holds indeed.  (We may wonder what we can do 
without assuming ``$\theta$ a compact cardinal"; in short, 
if $\partial < \theta \wedge \alpha < 
\cf(\delta_\ell) \Rightarrow |\alpha|^\partial <
\cf(\delta_\ell)$, we can use the $\Delta$-system lemma; otherwise use
\cite[\S7]{Sh:620} which gives a weaker relative of the
$\Delta$-system lemma for, e.g. $\lambda = \mu^+,\mu > 2^{\cf(\mu)}$.)

\sn
4) Easy.

\sn
5) Like the proof of part (3), noting that $<_{\lex}$ is definable in $M$.
\end{PROOF}

\begin{definition}\label{y50}
For a linear order $I$ and $\sigma < \theta$ we define $M_{I,\sigma}$
as the following model:
\mn
\begin{enumerate}
    \item  The universe is  $\{\eta : \eta$ a sequence of length 
    $< \sigma,\ \eta(2i) \in \bbQ,\ \eta(2i+1) \in I\}$.
\sn
    \item  $<^M$ is the natural lexicographic order.
\end{enumerate}
\end{definition}

\begin{example}\label{y56}
1) There is a complete $T \subseteq \bbL_{\theta,\theta}(\{<\})$ which
is definably unstable, 1-unstable but ``3-stable and 4-stable".

\sn
2) We can add ``$T$ has $\theta$-$\ncp$", see Definition \ref{a17} below.
\end{example}

\begin{PROOF}{\ref{y56}}
1) Let $\tau = \{<\}$ and for any cardinality $\lambda$ we define a
$\tau$-model $M_\lambda$ by:
\mn
\begin{enumerate}
    \item  $s \in M_\lambda$ \underline{iff} for some 
    $\alpha = \alpha(s) < \lambda$, $s$ is a function from $\alpha$ to $\{0,1\}$
such that the set $\{\beta < \alpha:s(\beta)=1\}$ is finite$\}$
\sn
    \item  $M_\lambda \models ``s < t"$ \underline{iff}  $s \triangleleft t$.
\end{enumerate}
\mn
Let $T = \Th_{\bbL_{\theta,\theta}}(M_\lambda)$. 

Now
\mn
\begin{enumerate}
    \item[$(*)$]  if $M$ is a model of $T$ \underline{then}  for some cardinal $\lambda$ and $M'$ we have:
\sn
    \begin{enumerate}
        \item  $M'$ is isomorphic to $M$
\sn
        \item  $M' \subseteq M_\lambda$
\sn
        \item  $|M'|$ is closed under initial segments
\sn
        \item  if $\eta \in M'$ and $\gamma < \lambda$ then $\eta \caret \LL (0)_\gamma\RR \in M'$.
    \end{enumerate}
\end{enumerate}
\mn
The rest should be clear.

\noindent
2) As above use the linear order of \ref{y40} instead of $\theta$.
\end{PROOF}

\noindent
We now sum up the implications among the generalizations of stable.
\begin{conclusion}
\label{y62}
1) For $T$ a complete $\bbL_{\theta,\theta}$-theory the following
implications hold:
\begin{enumerate}
\item[$(a)$]  5-unstable $\Rightarrow$ 4-unstable $\Rightarrow T$ is
  unstable $\Rightarrow T$ is $\lambda$-unstable for some $\lambda =
  \lambda^{< \theta} + \theta + \lambda^{|T|} \Rightarrow$ definably
  unstable $\Rightarrow$ 2-unstable $\Leftrightarrow$ 1-unstable.
\sn
    \item[$(b)$]  3-unstable $\Rightarrow$ definably unstable $\Rightarrow$ 2-unstable $\Leftrightarrow$ 1-unstable.
\end{enumerate}
\mn
2) The results in part (1) are best possible, i.e. all implications not
appearing there fail for some such $T$.
\end{conclusion}

\begin{PROOF}{\ref{y62}}
1) \underline{Clause (a)}:
\mn
\begin{enumerate}
\item[$\bullet_1$]  ``$T$ is 5-unstable implies $T$ is 4-unstable".
\end{enumerate}
\mn
[Why?  By \ref{y5}(1)(a) $\Rightarrow$ (b).]
\mn
\begin{enumerate}
\item[$\bullet_2$]  ``4-unstable implies $T$ is unstable" .
\end{enumerate}
\mn
[Why?  By \ref{y5}(1)(b) $\Rightarrow$ (c).]
\mn
\begin{enumerate}
\item[$\bullet_3$]  ``$T$ implies $\lambda$-unstable for some $\lambda
  = \lambda^{< \theta} + \lambda^{|T|}$".
\end{enumerate}
\mn
[Why?  By \ref{y5}(1)(c) $\Rightarrow$ (e).]
\mn
\begin{enumerate}
\item[$\bullet_4$]  ``$\lambda$-unstable for some $\lambda =
  \lambda^{< \theta} + \lambda^{|T|}$ implies definably unstable".
\end{enumerate}
\mn
[Why?  By \ref{y5}(a)(b) $\Rightarrow$ (e).]
\mn
\begin{enumerate}
\item[$\bullet_5$]  ``definably unstable implies 2-unstable".
\end{enumerate}
\mn
[Why?  By \ref{y5}(1)(g) $\Rightarrow$ (i).]
\mn
\begin{enumerate}
\item[$\bullet_6$]  ``2-unstable is equivalent to 1-unstable".
\end{enumerate}
\mn
[Why?  By \ref{y5}(1)(i) $\Rightarrow$ (j).]
\medskip

\noindent
\underline{Clause (b)}:
\mn
\begin{enumerate}
\item[$\bullet_1$]  ``3-unstable implies definably unstable".
\end{enumerate}
\mn
So we are done.

[Why?  By \ref{y5}(2), the second phrase.  The other implications hold by
  clause (a).]
\mn
2) Note that:
\mn
\begin{enumerate}
\item[$\bullet_1$]  ``1-unstable does not imply definably unstable".
\end{enumerate}
\mn
[Why?  This holds by \ref{y41}.]
\mn
\begin{enumerate}
\item[$\bullet_2$]  ``3-unstable does not imply stable.
\end{enumerate}
\mn
[Why?  This holds by \ref{y41}(2).]
\mn
\begin{enumerate}
\item[$\bullet_3$]  ``4-unstable does not imply 3-unstable".
\end{enumerate}
\mn
[Why?  This holds by \ref{y45}(3).]
\mn
\begin{enumerate}
\item[$\bullet_4$]  ``4-unstable does not imply definably 5-unstable".
\end{enumerate}
\mn
[Why?  This holds by \ref{y45}(3).]

So we are done.
\end{PROOF}
\newpage

\section {Saturation of ultrapowers} \label{2}
\bigskip

We define versions of notions of saturation and deal with basic properties.

Note that unlike the first order case, two
$(\lambda,\lambda,\bbL_{\theta,\theta})$-saturated
models of cardinality $\lambda$ are not necessarily isomorphic, see
Definition \ref{a5} and examples in \ref{a6}.  We consider calling the
notion in \ref{a5}, compact instead of saturated, but the word compact
has been over used.

\begin{context}  
\label{a0}
$\theta$ a compact cardinal.
\end{context}

\begin{definition}  
\label{a5}
1) We say $M$ is fully $(\lambda,\partial,L)$-saturated\footnote{Maybe ``compact'' would be more suitable, but late changes are dangerous.}
(we may omit the fully; where $L \subseteq \cL(\tau_M)$ and $\cL$ is a logic; we may
write $\cL$ if $L = \cL(\tau_M)$, the default value is 
$\cL = \bbL_{\theta,\theta})$
 \underline{when} : if $\Gamma$ is a set of $< \lambda$ formulas from $L$ 
with parameters from $M$ with $< 1 + \partial$ free variables, and
$\Gamma$ is $(< \theta)$-satisfiable in $M$, \underline{then} 
$\Gamma$ is realized in $M$.  

\noindent
2) We say ``locally" \underline{when}  using one $\varphi = \varphi(\bar x,\bar y) \in
\cL$, i.e. all members of $\Gamma$ have the form $\varphi(\bar x,\bar
b)$, that is:
\mn
\begin{enumerate}
\item[$(a)$]  if $\partial \le \theta$, then we consider a set of
  formulas of the form $\{\varphi(\bar x_{[\eps]},\bar
  a_\alpha):\alpha < \alpha_*\}$ where $\eps < \partial,\alpha_*
 < \lambda$ (so $\lh(\bar x)=\eps$)
\sn
\item[$(b)$]  if $\partial > \theta$ letting 
$j_* = \lh(\bar x)$, we consider a set of
 formulas of the form $\{\varphi(\LL x_{\eps(i,\alpha)}:i <
 j_\alpha \RR,\bar a_\alpha):\alpha < \alpha_*\}$ where
$\{\eps(i,\alpha):i < j_\alpha,\alpha < \alpha_*\} \subseteq j_*$.
\end{enumerate}
\mn
3) In the full case omitting $\partial$ means $\partial = \lambda$ and in
the local case omitting $\partial$ means $\partial = \theta$; writing
``$\le \partial$" means $\partial^+$.
Omitting $L$ means $\bbL_{\theta,\theta}$ and omitting $\lambda$ means
 $\lambda = \|M\|$.

\noindent
4) Assume $\eps$ is an ordinal $< \theta$ and $\Delta$ is a set of 
formulas of the form $\varphi(\bar x_{[\eps]},\bar y)$.  We say $M$ is
$(\lambda,\Delta)$-saturated when: $\Gamma$ is realized in $M$
whenever $\Gamma$ is a set of $< \lambda$ formulas of the form
$\varphi(\bar x_{[\eps]},\bar a),\bar a \subseteq M$, which is
$(< \theta)$-satisfiable in $M$.  May write
$(\lambda,\theta,\Delta)$-saturated abusing notation.
\end{definition}

\noindent
As said above, this notion does not have the most desirable properties
it has in the first order case as:
\begin{claim}  
\label{a6}
Let $\tau = \{<\}$, $<$ a two-place predicate.

\sn
1) If $T = \Th_{\bbL_{\theta,\theta}}(\theta,<)$, \underline{then}  no model of
$T$ is $(\theta^+,1,\bbL_{\theta,\theta}(\tau))$-saturated.

\sn
2) There is a complete $T \subseteq \bbL_{\theta,\theta}(\tau)$ such
that: $\tau = \tau_T$ is finite and
if $\mu = \mu^{< \kappa}$, $\kappa = \cf(\kappa) \ge \theta$ 
(so possibly $\mu = \kappa$) \underline{then} 
$T$ has non-isomorphic   
$(\kappa,\kappa,\bbL_{\theta,\theta}(\tau))$-saturated models of
cardinality $\mu$ (but a unique smooth one --- see the proof).

\sn
3) In part (2), if $\mu$ is strong limit singular then:
\mn
\begin{enumerate}
\item[$(A)$]  if $\mu$ is of cofinality $\ge \theta$
\underline{then}  $T$ has non-isomorphic special models of cardinality $\mu$; where:
\sn
\begin{enumerate}
\item[$\bullet$]   $M$ is called special when $M$ is the union of the
   $\prec_{\bbL_{\theta,\theta}}$-increasing sequence $\olsi M = 
\LL M_\alpha:\alpha < \cf(\mu)\RR$ such that 
$\|M_\alpha\| < \mu$ and $M_{\alpha +1}$ is
   $(\|M_\alpha\|^+,\|M_\alpha\|^+,\bbL_{\theta,\theta}(\tau))$-saturated
\end{enumerate}
\sn
\item[$(B)$]  if $\mu$ has cofinality $\in [\aleph_1,\theta)$ then $T$
  has $> \mu$ special models of cardinality $\mu$ pairwise
  non-isomorphic; but unique if we demand ``$M$ is smooth" (see in the
  proof)
\sn
\item[$(C)$]  if $\mu$ has cofinality $\aleph_0$ \underline{then}  $T$ has a
  special model of cardinality $\mu$ and this model is unique up to
  isomorphism.
\end{enumerate}
\end{claim}

\begin{remark}
\label{a7}
1) The claim above tells us that saturation does not behave as in the
   first order case, neither concerning existence nor concerning
   uniqueness.

\sn
2) So in part \ref{a6}(2), the counterexample is when $\mu = \kappa$;
note that there are such $\mu$-s: any successor 
of strong limit singular cardinal
which is $\ge \theta$ by \cite{So74}.

\sn
3) Concerning \ref{a6}(3) note that we regain uniqueness if we demand
smoothness; see \cite[2.15=L88r-2.10,2.17=L88r-2.11.1]{Sh:88r}.

\sn
4) Concerning \ref{a6}(3)(c), recall that Chang proved that for such $\mu$,
if two models are $\bbL_{\mu^+\!,\mu}$-equivalent then they are
isomorphic.

\sn
5) Let $L = \bbL_{\theta,\theta}(\tau_M)$.  Why in first order
logic in \ref{a5} we use only $\partial = 1$ and here not?  If $(\forall
\alpha < \lambda)[|\alpha|^{< \theta} < \lambda]$ then the cases
$\partial = 1$ and $\partial =2$ are equivalent \underline{but} for
$\partial = \aleph_1$, a type $p=p(\bar x_{[\omega]})$ may not be
realized though the model is $(\lambda,\partial,L)$-saturated for every
finite $\partial$, unlike first order logic.
\end{remark}

\begin{PROOF}{\ref{a6}}
1) Any model of $T$ is isomorphic to $M = (\delta,<)$ for some ordinal
   $\delta$ of cofinality $\ge \theta$.  Hence it is enough for such
   $\delta$ to prove that $M = (\delta,<)$
satisfies the desired conclusion.  If $\delta = \theta$ the model $M$
omits the type $\{\alpha < x:\alpha < \theta\}$ and if $\delta > \theta$
   then $M$ omits $\{\alpha < x \wedge x < \theta:\alpha < \theta\}$.

\sn
2) Let $\tau = \{<\},<$ a two-place predicate; toward defining a
theory $T$ we first let $\gk = (K,\le_{\gk})$ be defined by:
\mn
\begin{enumerate}
\item[$(*)_1$]  
\begin{enumerate}
\item[(a)]  $K$ is the class of $\tau$-models $M$
  which are trees in the model theoretic sense, i.e. satisfies:
\sn
\begin{itemize}
\item  $x < y \rightarrow x \ne y$
\sn
\item   $(x < y \wedge y < z) \rightarrow x < z$
\sn
\item   $(x < z \wedge y < z) \rightarrow (x < y \vee y < x \vee y = x)$
\end{itemize}
\sn
\item[(b)]  $\le_{\gk}$ is the following two-place
  relation on $K_1:M \le_{\gk} N$ \underline{iff} 
\sn
\begin{enumerate}
\item[($\alpha$)]  $M \subseteq N$
\sn
\item[($\beta$)]   if $\LL a_n:n < \omega\RR$ is
increasing with no upper bound in $M$, \underline{then}  it has no upper bound in $N$.
\end{enumerate}
\end{enumerate}
\mn
Now observe
\mn
\begin{enumerate}
\item[$(*)_2$]  $\gk$ is a weak a.e.c., in the sense that:
\begin{enumerate}
\item[(A)] 
\begin{enumerate}
\item[(a)]  $K$ and $\le_{\gk}$ are closed under isomorphisms
\sn
\item[(b)]  $\le_{\gk}$ is a partial order and $M \in K
  \Rightarrow M \le_{\gk} M$
\sn
\item[(c)]   if $\LL M_i:i < \delta\RR$ is
  $\le_{\gk}$-increasing then $M_\delta \defeq \bigcup\limits_{i < \delta}
  M_i \in K$ and $i < \delta \Rightarrow M_i \le_{\gk} M_\delta$
\sn
\item[(d)]  if $\LL M_i:i \le \delta\RR$ is
  $\le_{\gk}$-increasing then $\bigcup\limits_{i < \delta} M_i
  \le_{\gk} M_\delta$ provided that $\cf(\delta) \ne \aleph_0$
\sn
\item[(e)]   if $M_1 \subseteq M_2$ are $\le_{\gk} N$
  \underline{then}  $M_1 \le_{\gk} M_2$
\sn
\item[(f)]   LST: if $\lambda = \lambda^{\aleph_0}$ \underline{then}  
the LST-property holds up to $\lambda$
\end{enumerate}
\end{enumerate}
\sn
\item[(B)]
\begin{enumerate}
\item[(a)]  $\gk$ satisfies the amalgamation property, in fact, essentially
  disjoint union suffice, i.e. if $M_0 
\subseteq M_1,M_0 \subseteq M_2$ are from $K$ and $M_1 \cap M_2 =
M_0$, \underline{then}  $M_3 = M_1 \cupdot M_2$ does $\le_{\gk}$-extend $M_\ell$
for $\ell=0,1,2$.
\newline
Note that to say $M_3 \defeq M_1 \cupdot M_2$ means 
$M_3$ has universe $|M_1| \cup |M_2|$ and $<^{M_*}$ is defined by
 $a_1 <^{M_3} a_2$ \underline{iff}  at least one of the following holds:
\sn
\begin{enumerate}
\item[($\alpha$)]  $a_1 <^{M_1} a_2$
\sn
\item[($\beta$)]  $a_1 <^{M_2} a_2$
\sn
\item[($\gamma$)]  $a_1 \in M_1 \setminus M_0$ and
  $a_2 \in M_2 \setminus M_0$ and for some $b \in M_0$
\sn
\begin{itemize}
\item  $a_1 \le^{M_1} b <^{M_2} a_2$
\end{itemize}
\sn
\item[($\delta$)]   as in ($\gamma$) but we interchange $M_1,M_2$
\sn
\item[($\eps$)]   $a_1 \in M_1 \setminus M_0$ and $a_2 \in M_1
  \setminus M_0$ and the sets $\{b \in M_0 : a_1 \le^{M_1} b\}$, $\{b \in
  M_0:a_2 \le^{M_1} b\}$ are equal and non-empty (recalling $M_\ell$
  is a tree)
    \end{enumerate}
\sn
    \item[(b)]   similarly $\gk$ has the JEP, even as the disjoint union
\sn
\item[(c)]  (skewed amalgamation) if $M_0 \subseteq
  M_1$ and $M_0 \le_{\gk} M_2$ all from $K$
and $M_1 \cap M_2 = M_0$ \underline{then}  $M_3 = M_1 \cup M_2$ defined as in (B)(a)
 above satisfies $M_2 \subseteq M_3$ and $M_1 \le_{\gk} M_3$
\sn
\item[(d)]   if $A \subseteq M \in K,A \ne \varnothing$
  then $M \rest A \in K$ (but possibly $M \rest A \nleq_{\gk} M$).
\end{enumerate}
\end{enumerate}
\end{enumerate}
\mn
[Why?  For clause (B)(c), clearly $\ell \le 3 \Rightarrow M_\ell \in
K$ and $\ell < 3 \Rightarrow M_\ell \subseteq M_3$.  For proving $M_1
\le_{\gk} M_3$ let $\bar a = \LL a_n:n < \omega\RR$ be
$<^{M_1}$-increasing and $c \in M_3 \setminus M_1$ be an upper bound
(for $<^{M_3}$) of $\{a_n:n < \omega\}$.  So one of the five cases 
in (B)(a) holds for infinitely many pairs $(a_n,c)$, so without loss of generality it holds for all $(a_n,c)$.

If clause ($\alpha$) - then $c \in M_1$ and we are done, and if clause
($\beta$) then $a_n \in M_0$ and use $M_0 \le_{\gk} M_2$.  If clause
($\gamma$), then there is $b_n \in M_0$ such that $a_n \le^{M_1} b_n
\le^{M_2} c$, so $b_n \in M_0,\{b_n:n < \omega\}$ linearly ordered, by
Ramsey theorem (as $M_1$ is a tree) without loss of generality $\bar b =
\LL b_n:n < \omega\RR$ is monotone.  If $\bar b$ is increasing,
then it is increasing in $M_1$ and clearly has no upper bound in $M_1$
(as it will be an upper bound of $\bar a$), hence in $M_0$ but it has one in
$M_2$, contradicting $M_0 \le_{\gk} M_2$.  If $\bar b$ is (monotone
and) not increasing then it is $\le$-decreasing hence $b_0 \in M_0 \subseteq
M_1$ is an upper bound of $\bar a$, contradiction.  

Next, if we use Clause ($\delta$), the proof is easier:
$\bigwedge\limits_{n} a_n \in M_2$ hence $\bigwedge\limits_{n} a_n \in
M_1 \cap M_2 = M_0$ and $c \in M_3 \setminus M_1 = M_2 \setminus
M_0$ so use $M_0 \le_{\gk} M_1$.

Lastly, if clause ($\eps$), then there is $b \in M_0$ above all the
$a_n$-s so we finish as earlier.

So we are done proving $(*)_2$.]

In particular
\mn
\begin{enumerate}
\item[$(*)_3$]  if $\LL M_i:i < \delta\RR$ is
  $\le_{\gk}$-increasing \underline{then}  $\bigcup\limits_{i < \delta} M_i \in K$
  does $\le_{\gk}$-extend $M_i$ for $i < \delta$.
\end{enumerate}
\mn
Next for $\kappa \ge \theta$ and let
\mn
\begin{enumerate}
    \item[$(*)_4$]  $K^{\ec}_\kappa = \{M \in K: \text{if } 
    M \le_{\gk} N,\ A \subseteq M$ has cardinality $< \kappa$ and $\bar a \in {}^{\kappa>}\!N$ then some $\bar b \in {}^{\lh(\bar a)}\!M$ realizes $\tp_{\qf}(\bar a,A,N)\}$.
\end{enumerate}
\mn
Clearly
\mn
\begin{enumerate}
\item[$(*)_5$] 
\begin{enumerate}
\item[(a)]   if $M_1 \in K$ has cardinality $\le \mu =
\mu^{< \kappa}$ \underline{then}  some $M_2 \in K^{\ec}_\kappa$ 
 has cardinality $\mu$ and $\le_{\gk}$-extends $M_1$
\sn
\item[(b)]   any $M \in K^{\ec}_\kappa$ has elimination
  of quantifiers in $\bbL_{\theta,\theta}$ up to
$x < y$, $x=y$ and $\varphi_*(\bar x_{[w]}) = (\exists y)
(\bigwedge\limits_{n} x_n < y)$; also $M$ is
$(\kappa,\kappa,\bbL_{\theta,\theta})$-saturated, recalling
$\kappa \ge \theta$
\sn
\item[(c)]   any $M_1,M_2 \in K^{\ec}_\kappa$ are
$\bbL_{\theta,\theta}$-equivalent and even
$\bbL_{\infty,\theta}$-equivalent
\sn
\item[(d)]  $K^{\ec}_{\kappa_2} \subseteq
  K^{\ec}_{\kappa_1}$ when $\theta \le \kappa_1 \le \kappa_2$.
\end{enumerate}
\end{enumerate}
\mn
Hence we define $T$ as follows: (it is well defined by $(*)_5$(c))
\mn
\begin{enumerate}
\item[$(*)_6$]  $T = \Th_{\bbL_{\theta,\theta}}(M)$ whenever 
$M \in  
K^{\ec}_\theta$.
\end{enumerate}
\mn
So
\mn
\begin{enumerate}
    \item[$(*)_7$]  $T$ is a complete $\bbL_{\theta,\theta}$-theory, $\tau_T = \{<\}$ and if $\kappa \ge \theta$, 
    $\mu = \mu^{< \kappa}$
    \underline{then}  $T$ has a $(\kappa,\kappa,\bbL_{\theta,\theta})$-saturated model of cardinality $\mu$ (even extending
    any pre-given $M \in \Mod_T$ of cardinality $\le \mu$).
\end{enumerate}
\mn
Lastly
\mn
\begin{enumerate}
\item[$(*)_8$]  if $\mu = \mu^{< \kappa},\kappa \ge \theta$ \underline{then} 
there are $> \mu$ pairwise non-isomorphic
  $(\kappa,\kappa,\bbL_{\theta,\theta})$-saturated models of $T$ of
  cardinality $\mu$.
\end{enumerate}
\mn
Why?  First, 

\mn
\underline{\textbf{Case 1}}:  assume $\mu$ is regular uncountable.

For $M \in K$ with universe $\lambda$ let 
$$\smooth_0(M) = \{\delta < \mu : \cf(\delta) = \aleph_0 \text{ and }
M \rest \delta \le_{\gk} M\}$$ 
and for any $M \in K$ of cardinality $\lambda$ let 
$\smooth(M) = \smooth_0(N)/\cD_\mu$ for any
$N$ isomorphic to $M$ with universe $\lambda$, recalling 
$\cD_\mu$ is the club filter on $\mu$.  

This makes sense because:
\mn
\begin{itemize}
\item  If $M_1,M_2 \in K$ have universe $\lambda$ then $\smooth_0(M_1)
  = \smooth_0(M_2) \mod \cD_\mu$.
\end{itemize}
\mn
We say such $M$ is smooth when $\smooth(M) = \lambda/\cD_\lambda$.

Easily for any $S \subseteq \{\delta < \lambda:\cf(\delta) =
\lambda\}$ there is $M = M_S \in \Mod_T$ of cardinality $\mu$ such that
$\smooth(M) = S/\cD_\mu$ and even $M_S \in K^{\ec}_\kappa$.  So if
$S_1,S_2 \subseteq \lambda$ and $S_1 \setminus S_2$ is stationary then
$M_{S_1} \not\cong M_{S_2}$, so by $(*)_5(c)$ we are done.

Note
\mn
\begin{enumerate}
\item[$\boxplus_1$]  If $\mu = \mu^{<\mu} > \aleph_0$ \underline{then}  up to
isomorphism there is one and only one smooth $M \in K^{\ec}_\mu$
  which is $(\mu,\mu,\bbL_{\theta,\theta})$-saturated of cardinality
  $\mu$; where
\sn
\item[$\boxplus_2$]  $M \in K$ of cardinality $\mu = \cf(\mu)$ is smooth
  when $\smooth(M)=\varnothing/\cD_\mu$.
\end{enumerate}
\mn
Details on $\boxplus_2$ see $(*)_9$ - $(*)_{11}$ in the end of the proof.

Second, next 

\mn
\textbf{\underline{Case 2}}:  Assume $\mu$ is singular of cofinality 
$\ge \aleph_1$.

For special models in our context the hope was to show that any two special model
are $\bbL_{\infty,\theta}$-equivalent.

Let $\bar\kappa = \LL \kappa_i : i < \cf(\mu)\RR$ be increasing
with limit $\mu$ such that $\kappa_i > \theta$ and
$\lambda_i = 2^{\kappa_i} < \kappa_{i+1}$.

So we can consider:
\mn
\begin{enumerate}
    \item[$\boxplus_3$]   $K^{\sep}_{\bar\kappa} = \big\{\bigcup\{M_i : i < \cf(\mu)\} : M_i \in K^{\ec}_{\kappa_i}$ is $\kappa^+_i$-saturated of cardinality $\lambda_i,\ \le_{\gk}$-increasing with $i\big\}$.
\end{enumerate}
\mn
Now
\mn
\begin{enumerate}
    \item[$\oplus$]
    \begin{enumerate}
        \item[(a)]  Any $M \in K^{\sep}_{\bar\kappa}$ is special and $K^{\sep}_{\bar\kappa} \ne \varnothing$. Moreover, if 
        $M_1 \in K$ has cardinality $\le \mu$ \underline{then}  
        there is $N \in K^{\sep}_{\bar\kappa}$ such that 
        $M \le_{\gk} N$.
\sn
        \item[(b)]  Any two models from $K^{\sep}_{\bar\kappa}$ are $\bbL_{\infty,\theta}$-equivalent.
\sn
        \item[(c)]  There are non-isomorphic 
        $M_1,M_2 \in K^{\sep}_{\bar\kappa}$.
    \end{enumerate}
\end{enumerate}
\mn
Why does $\oplus$ hold?
\medskip

\noindent
\underline{Clause (a)}:  The existence of $N \in K^{\ec}_\mu$ as well
as ``any $M \in K^{\sep}_{\bar\kappa}$ is special" are
obvious by the definitions.  For the second demand (density) assume
$M \in K$ has cardinality $\le \mu$, without loss of generality of cardinality $\mu$.
Let $|M|$ be $\bigcup\limits_{i < \kappa} A_i$ with 
$|A_i| = \lambda_i$.

We choose $M_i$ by induction on $i \le \kappa$ such that:
\mn
\begin{enumerate}
    \item[$\oplus_1$]
    \begin{enumerate}
        \item[(a)]  $M_i \subseteq M$ has cardinality [???]
\sn
        \item[(b)]  $\LL M_j : j \le i\RR$ is $\le_{\gk}$-increasing.
\sn
        \item[(c)]  $M_i \le_{\gk} M$
\sn
        \item[(d)]  If $i = j+1$ then $A_j \subseteq M_i$.
    \end{enumerate}
\end{enumerate}
\mn
Next we choose $N_i$ by induction on $i \le \kappa$ such that:
\mn
\begin{enumerate}
    \item[$\oplus_2$]
    \begin{enumerate}
        \item[(a)]  $N_i \in K$ is $\kappa_i$-saturated of cardinality $\lambda_i$.
\sn
        \item[(b)]  $\LL N_j : j \le i\RR$ is $\le_{\gk}$-increasing.
\sn
        \item[(c)]  $M_i \le_{\gk} N_i$
\sn
        \item[(d)]  $N_i \cap M = M_i$.
    \end{enumerate}
\end{enumerate}
\mn
Why can we carry the induction?  For $i=0$ obviously, by the $\JEP$
and the density of $\kappa^+_i$-saturated in cardinality $\lambda_i$.
For $i=j+1$ recalling $\gk$ has amalgamation ($\LST$ and as above).
For limit $i$ of cofinality $> \aleph_0$ - similarly.  

Lastly, for $i$ of
cofinality $\aleph_0$ the proof is as in $(*)_2(B)(c)$.
\medskip

\noindent
\underline{Clause (b)}:  Is obvious when $\cf(\mu) \ge \theta$.

But even without this assumption we can prove a stronger result:
\mn
\begin{enumerate}
\item[$\oplus_3$]
\begin{enumerate}
\item[(b)$^+$]  if $M_\ell \in K^{\sep}_{\bar\kappa}$ for $\ell=1,2$
  and $\kappa < \mu$ then $M_1,M_2$ are
  $\bbL_{\infty,\kappa}$-equivalent.
\end{enumerate}
\end{enumerate}
\mn
Why?  Without loss of generality $\kappa = \lambda^+_0 \ge \cf(\mu)$ and $\olsi M_\ell =
\LL M_{\ell,i}:i < \cf(\mu)\RR$ witness $M_\ell \in
K^{\sep}_{\bar\kappa}$.

Let $\cA_\ell$ be the set of $A$ such that:
\mn
\begin{enumerate}
\item[($\alpha$)]  $A \subseteq M_\ell,|A| \le \lambda_0$
\sn
\item[($\beta$)]  if $a \in A \setminus M_i,i < \cf(\mu)$ and
  $B^\ell_{a,i} = \{b \in M_{1,i}:b <_{M_1} a\}$ has cofinality $\le 
\lambda_0$ then $B^\ell_{a,i} \cap A$ is cofinal in $M_\ell$
\sn
\item[($\gamma$)]  if $a_n \le^{M_\ell} a_{n+1} \le^{M_\ell} b$ 
and $a_n \in A \cap M_i$ for $n < \omega,b \in M_j$,
  \underline{then}  there is such $b$ in $A \cap M_j$.
\end{enumerate}
\mn
Let $\cF_0$ be the set of $f$ such that:
\mn
\begin{itemize}
    \item  For some $A_1 \in \cA_1$ and $A_2 \in \cA_2$, $f$ is an isomorphism from $M_1 \rest A_1$ onto $M_2 \rest A_2$ preserving the property in ($\beta$) above.
\end{itemize}
\mn
Now $\cF$ witness ``$M_1,M_2$ are $\bbL_{\infty,\kappa}$-equivalent.
We leave the checking to the reader.

What about
\underline{Clause (c)}: ``Two non-isomorphic ones"?  We give three ways to
do this.
\medskip

\noindent
\underline{The First Way}:

We can get $2^\mu$ pairwise non-isomorphic
$(\kappa,\kappa,\bbL_{\theta,\theta})$-equivalent models which are
special and even in $K^{\sep}_{\bar\kappa}$
when $\mu$ is strong limit singular.  A way to do it is to work 
as in \cite{Sh:511} where we construct ``complicated" 
sequences of subtrees of ${}^{\sigma\ge}\lambda$ and use 
them to construct, e.g. Boolean Algebras.  We do not elaborate, but
shall give details in the other ways.
\medskip

\noindent
\underline{A Second Way}:  

Giving in some respect a stronger version, when $\mu$
is strong limit of cofinality $\kappa > \aleph_0$ is as follows.  Let
$\LL \lambda_i:i < \kappa\RR$ be increasing continuous with
limit $\lambda_\kappa = \mu$, $\lambda_{i+1} = (\lambda_{i+1})^{\lambda_i}$, $\lambda_0 =
(\lambda_0)^{\aleph_0}$ and $S_0,S_1 \subseteq S^\kappa_{\aleph_0}$ be
stationary disjoint and $\eps \in S_1 \Rightarrow
\lambda_{\eps +1} = 2^{\lambda_\eps}$.  We choose $\bar
M_\eps$ by induction on $\eps \le \kappa$ such that:
\mn
\begin{enumerate}
    \item[$(*)_{8.1}$]
    \begin{enumerate}
        \item[(a)]  $\olsi M_\eps = \LL M_\eta : \eta \in {}^{(\lambda_\eps)}2\RR$
\sn
        \item[(b)]  $\LL M_{\eta \rest \lambda_\zeta} : \zeta \le \eps \RR$ is $\subseteq$-increasing continuous.
\sn
        \item[(c)]  $M_\eta \in K$ has universe $\lambda_{\lh(\eta)}$.
\sn
        \item[(d)]  $M_\eta \in K^{\ec}_{\lambda_\eps}$ if $\eta \in
  {}^{\lambda_{\eps +1}}2$
\sn
        \item[(e)]  $M_{\eta \rest (\zeta +1)} \le_{\gk} M_\eta$ 
        for $\eta \in {}^{\lambda_\eps}2$ and $\zeta < \eps$.
\sn
        \item[(f)]  If $\eta \ne \nu \in {}^{(\lambda_\eps)}2$, 
        $\eps = \zeta +1$, $\zeta \in S_1$, and 
        $f \in \cF_{\eta,\nu}$ (see below) \underline{then}  
        for some $\rho \in \lim(I_f)$ (see below) we have: there is $a \in M_\eta$ such that $(\forall n)[\rho(n) < a]$ but for no $b \in M_\nu$ do we have 
        $M_\nu \models (\forall n)[f(\rho(n)) < b]$, where
        \begin{enumerate}
            \item[$\oplus$]  $\cF_{\eta,\nu}$ is the set of functions $f$ such that
            \begin{itemize}
                \item  $\dom(f)$ is a subtree of 
                ${}^{\omega >}(\lambda_\eps)$ with $\lim(\dom(f))$ of cardinality $2^{\lambda_\eps}$.
\sn
                \item  $\rho \in \dom(f) \Rightarrow M_\eta 
                \models ``\LL \rho(\ell):\ell < \lh(\rho)\RR$ is increasing"
\sn
                \item  For every $\xi < \zeta$, for all but 
                $< \lambda_\eps$ members $\rho$ of $\dom(f)$, we have $\Rang(\rho) \nsubseteq \lambda_\xi$.
\sn
                \item  If $\rho \caret \LL \alpha \RR,\varrho \caret \LL \beta \RR \in \Dom(f)$ are $\triangleleft$-incomparable then $M_\eta \models ``\rho \caret \LL \alpha \RR,\rho \caret \LL \beta \RR$ are incomparable".
\sn
                \item  $f$ is one to one.
            \end{itemize}
        \end{enumerate}
    \end{enumerate}
\end{enumerate}
\mn
Now
\mn
\begin{enumerate}
    \item[$(*)_{8.2}$]  we can carry the induction.
\end{enumerate}
\mn
[Why?  For $\eps =0$ trivially, and $\eps$ limit use union; for
$\eps = \zeta +1$, $\zeta \notin S_1$ use $(*)_5$(a) and for $\eps =
\zeta +1$, $\zeta \in S_1$ by cardinality consideration we can take care
of clause (f) and then use $(*)_5$(a) to take care of clause (d).]
\mn
\begin{enumerate}
\item[$(*)_{8.3}$]  if $\eta \in {}^\mu 2$ then $M_\eta$ is a special
  model of $T$.
\end{enumerate}
\mn
[Why?  By $(*)_{8.1}$(b),(c),(d).]
\mn
\begin{enumerate}
    \item[$(*)_{8.4}$]  If $\eta \ne \nu \in {}^\mu 2$ then 
    $M_\eta \in K^{\ec}_\mu$ is not $\le_{\gk}$-embeddable 
    into $M_\nu$.
\end{enumerate}
\mn
[Why?  By \cite[Claim 2.4,pg.111]{Sh:136}; see more in 
Rubin-Shelah \cite{Sh:117} and \cite[Ch.XI]{Sh:f}.]
\medskip

\noindent
\underline{Third Way}:  Giving $\mu^+$ non-isomorphic models is 
by the simple black box of \cite[\S1,1.5=L4.5A,pg.3]{Sh:309}, but we 
elaborate\footnote{Can we get $2^\mu$ ones?  In this particular case, 
yes, but we shall not elaborate; we can use
\cite[1.9=L4.6,pg.5]{Sh:309}.} giving a self contained proof.  
Let $\LL M_i:i < \mu\RR$ be a sequence of members of 
$K^{\ec}_\kappa$ (so models of $T$ each of cardinality $\mu$) and we shall find a model from $K$ of cardinality $\mu$ not $\le_{\gk}$-embeddable into any $M_i$, this clearly suffices by $\oplus$(a), the density.

We define a model $M \in K$ as follows:
\mn
\begin{enumerate}
    \item[(a)]  Its set of elements is the set of $\eta$-s such that
\sn
    \begin{enumerate}
        \item[($\alpha$)]  $\eta$ is a sequence of length $\le \omega$.
\sn
        \item[($\beta$)]  $\eta(0) \in \mu$ if $\lh(\eta) > 0$.
\sn
        \item[($\gamma$)]  $\eta(1+n) \in M_{\eta(0)}$ when 
        $1+n < \lh(\eta)$.
\sn
        \item[($\delta$)]  $M_{\eta(0)} \models ``\eta(1+n) < \eta(1+n+1)"$ when $1 + n + 1 < \lh(\eta)$.
\sn
        \item[($\eps$)]  If $\lh(\eta) = \omega$ then 
        $M_{\eta(0)} \models ``\neg(\exists x) \big[\bigwedge\limits_{n} \eta(1+n) < x \big]"$
    \end{enumerate}
\sn
    \item[(b)]  The order $<^M$ is $\triangleleft$, `being an initial segment.'
\end{enumerate}
\mn
Let $N \in K^{\ec}_\kappa$ be such that $M \le_{\gk} N$ and $N$ has
cardinality $\mu$.  Now indeed $i < \mu = N$ is not
$\le_{\gk}$-embeddable into $M_i$ as in \cite[\S1,1.5=L4.5A]{Sh:309};
in details toward contradiction assume $f$ is an
isomorphism from $N$ onto $M_i$.  Define $\eta_n \in N$ of length
$n+1$ by induction on $n$ as follows: if $n=0$ then 
$\eta_n = \LL i \RR \in N$ so $\eta_n(0)=i$. If $\eta_n$ has been
defined then we let
$\eta_{n+1} = \eta_n \caret \LL f(a_n)\RR$; it is well defined as
$a_n \in N$ hence $f(\eta_n) \in M_i$ and clearly $\eta_n
\triangleleft \eta_{n+1}$ hence $M_i \models f(\eta_n) < f(\eta_{n+1})$.

Now we ask: does the $<^{M_i}$-increasing sequence $\LL
f(\eta_n):n < \omega \RR$ have an upper bound in $M_i$?  If $a$ is
such an upper bound, $f^{-1}(a)$ is above $\{\eta_n:n < \omega\}$ so
necessarily is the sequence 
$\bigcup\limits_{n} \eta_n$ which does not belong to
$N$.  If there is no such $a,\eta = \bigcup\limits_{n} \eta_n \in N$
and $f(\eta)$ satisfies the demand, contradiction, so we are done
proving $(*)_8$.]

Why are we done proving part (3)?   Clauses (A),(B) -- the existence 
of $2^\mu$ pairwise non-isomorphic
special models from $K^{\ec}_\theta$ of cardinality $\lambda$ is
proved in ``the second way" of the proof of $(*)_8$ in part (1).  The
uniqueness of the smooth special model is just like Lemma
\cite[2.18=L88r-2.11,pg.18]{Sh:88r} and see Definition
\cite[2.15=L88r-2.10]{Sh:88r}, but see $(*)_{10}$ below.
\medskip

\noindent
\underline{Proof of $\boxplus_2$}:  Easy as above because here 
smoothness holds automatically as quoted above but we elaborate:
\mn
\begin{enumerate}
\item[$(*)_9$]  if $\lambda = \lambda^{< \lambda} > \aleph_0$ and
  $\alpha < \lambda \Rightarrow |\alpha|^{\aleph_0} < \lambda$ and
  $M_1,M_2$ are smooth $\le_{\gk}$-saturated $\lambda$-saturated models of
  cardinality $\lambda$, \underline{then}  $M_1,M_2$ are isomorphic.
\end{enumerate}
\mn
Why?  For $\ell=1,2$ let $\LL M_{\ell,\alpha}:\alpha <
\lambda\RR$ be $\le_{\gk}$-increasing continuous with union
$M_\ell$ such that $\alpha < \lambda \Rightarrow \|M_{\ell,\alpha}\| <
\lambda$; possible because $\alpha < \lambda \Rightarrow
|\alpha|^{\aleph_0} < \lambda$.]

Now we choose
$f_\eps,\alpha_{1,\eps},\alpha_{2,\eps},N_{1,\eps},N_{2,\eps}$ by
induction on $\eps < \lambda$ such that:
\mn
\begin{enumerate}
\item[$(*)_{10}$]
\begin{enumerate}
\item[(a)]  $N_{\ell,\eps} \le_{\gk} M_\ell$ has cardinality $<
  \lambda$
\sn
\item[(b)]  $f_\eps$ is an isomorphism from $N_{1,\eps}$ onto
$N_{2,\eps}$
\sn
\item[(c)]  $\alpha_{\ell,\eps} = \alpha(\ell,\eps)$ is increasing
  with $\eps$ for $\ell = 1,2$
\sn
\item[(d)]  if $\zeta < \eps,\ell = 1,2$ then
  $M_{\ell,\alpha(\ell,\zeta)} \subseteq N_{\ell,\eps} \subseteq
  M_{\alpha(\ell,\eps)}$.
\end{enumerate}
\end{enumerate}
\mn
The rest should be clear.
\mn
\begin{enumerate}
\item[$(*)_{11}$]  We have $M_1 \cong M_2$ when for $\ell=1,2$:
\sn
\begin{enumerate}
\item[(a)]  $M_\ell \in K$ is of cardinality $\mu$
\sn
\item[(b)]  $M_\ell = \bigcup\limits_{\ell < \kappa} M_{\ell,i}$
\sn
\item[(c)]  $\LL M_{\ell,i}:i < \kappa\RR$ is
  $\le_{\gk}$-increasing continuous
\sn
\item[(d)]  $M_{\ell,i+1} \in K$ is $\|M_{\ell,i}\|^+$-saturated.
\end{enumerate}
\end{enumerate}
\mn
Why true?   Similar to the proof above.  Note that if $\kappa =
\aleph_0$, then the ``continuous" in clause (c) is redundant.

\noindent
3) Clauses (A),(B) of \ref{a5}(3) were proved inside the proof of part
(2) and Clause (C) follows from $L_{\infty,\mu}$-equivalence.
\end{PROOF}

\begin{claim}  
\label{a8}
1) If $D \in \uf_\theta(I)$ is $(\lambda,\theta)$-regular
and $M_1,M_2$ are $\bbL_{\theta,\theta}$-equivalent and $\tau(M) =
\tau$ has cardinality $\le \lambda$ \underline{then}  $M^I_1/D,M^I_2/D$ are 
$\bbL_{\lambda^+\!,\lambda^+}$-equivalent, moreover
$\bbL_{\infty,\lambda^+\!,\lambda^+}$-equivalent (so one is 
$(\lambda^+,\lambda^+,\bbL_{\theta,\theta})$-saturated iff the other is).

\noindent
2) Similarly for $D \in \fil_\theta(I)$ which is $(\lambda,\theta)$-regular.
\end{claim}

\begin{remark}\label{a8d}
Recall that 
$$\bbL_{\chi,\mu,\gamma}(\tau) = \{\varphi(\bar x) \in
\bbL_{\chi,\mu}(\tau):\varphi(\bar x) \text{ has quantifier depth} <
\gamma\}$$ and $\bbL_{\infty,\lambda^+}(\tau) =
\bigcup\{\bbL_{\chi,\lambda^+}(\tau):\chi$ a cardinal $> \lambda\}$ and
$$\bbL_{\lambda^+\!,\lambda^+}(\tau) = \bigcup\{\bbL_{\lambda^+\!,\lambda^+\!,\gamma}:\gamma < \lambda^+\}.$$
Note that unlike the first order case we cannot demand
$\bbL_{\infty,\lambda^+}$-equivalence. 
\end{remark}

\begin{PROOF}{\ref{a8}}
1) Let $\gamma < \lambda^+$.
As $D$ is $(\lambda,\theta)$-regular there is a sequence 
$\LL (u_s,v_s,\Delta_s) : s \in I \RR$ such that 
$v_s \in [\gamma]^{< \theta}$, $u_s \in [\lambda]^{<\theta}$, 
$\Delta_s$ a set of $< \theta$ formulas of
$\bbL_{\theta,\theta}(\tau_T)$ and 
$$\alpha < \gamma \wedge \beta < \lambda \wedge \varphi(\bar x) \in \bbL_{\theta,\theta}(\tau_T) \Rightarrow \{s : \alpha \in v_s,\ \beta \in u_s, \text{ and }\varphi(\bar x) \in \Delta_s\} \in D.$$
For $s \in I$ let $\Game_s$ be the game
$\Game_{\Delta_s,u_s,v_s}(M_1,M_2)$; see Definition \ref{x2}.  As
$M_1,M_2$ are $\bbL_{\theta,\theta}$-equivalent by \ref{x4} the
protagonist wins this game $\Game_s$ which means that it 
has a winning strategy {\bf st}$_s$.  
Let $N_\ell = M^I_\ell/D$, and it suffices to find a
strategy {\bf st} for the protagonist in the game
$\Game_{\bbL_{\theta,\theta},\lambda,\gamma}$.  The strategy is
obvious (see proof in \cite[1.3=Ld11]{Sh:1101}) but we give details.

We say $\bfs$ is a reasonable state when it consists of:
\sn
\begin{enumerate}
    \item[(a)]  $\gamma_{\bfs} < \gamma$, $n_{\bfs} < \omega$
\sn
    \item[(b)]  A member $A$ of $D$.
\sn
    \item[(c)]  A set $J$ of cardinality $< \theta$.
\sn
    \item[(d)]  $f^\ell_\alpha \in M^I_\ell$ for $\ell \in \{1,2\}$,
    $\alpha < \lambda$.
\sn
    \item[(e)]  If $s \in A$ \underline{then} $\gamma_{\bfs} \in v_s$
    and $(n_{\bfs},g_{\bfs,s})$ is a winning state for the
    isomorphism player in the game $\Game_{\Delta_s,u_s,v_\alpha}$,
    where the partial function $g_{\bfs,s}$ is 
    $\big\{(f^1_\alpha(s),f^2_\alpha(s)):\alpha \in u_s \big\}$, so necessarily of cardinality $\le |u_s| < \theta$.
\end{enumerate}
\mn
2) The same proof as part (1) using only $\Delta$-s which are sets
of $<\theta$ atomic formulas of $\bbL_{\theta,\theta}(\tau_T)$.
\end{PROOF}

\begin{definition}  \label{a10}
1) Assume $\bar \mu = (\mu_1,\mu_2)$ but if $\mu_1 = \mu$, 
$\mu_2 = \theta$ we may write $\mu$; and 
$\lambda \ge \mu_1 \ge \mu_2 \ge \theta$.
We define a two-place relation 
$\blacktriangleleft_{\lambda,\bar\mu,\theta}$ on the
class of complete theories $T$ (in $\bbL_{\theta,\theta}$, of course) of 
cardinality $\le \lambda$.  We say $T_1 
\blacktriangleleft_{\lambda,\bar\mu,\theta} T_2$ \underline{iff}  for every $D \in 
\ruf_\theta(\lambda)$ and models $M_1,M_2$ of $T_1,T_2$,
   respectively we have: if $M^\lambda_2/D$ is locally
$(\mu^+_1,\mu^+_2,\bbL_{\theta,\theta})$-saturated 
\underline{then}  so is $M^\lambda_1/D$.

\sn
2) We say fully or write $\blacktriangleleft^{\ful}_{\lambda,\bar\mu,\theta}$,
   \underline{when}  we deal with full saturation.  We may omit $\bar\mu$ when
$\lambda = \mu_1,\mu_2 = \theta$.  We define
$\blacktriangleleft_{\lambda,\bar\mu,\theta},
\blacktriangleleft^{\ful}_{\lambda,\bar\mu,\theta}$ analogously. 
\end{definition}

\begin{remark}
\label{a11}
1) Note that $\blacktriangleleft$ is a quasi-order and not a partial
order, in particular, is \underline{not} a strict order.

\noindent
2) The relation of $\blacktriangleleft$ 
here to the classical one of Keisler is quite
close.  Keisler uses ``$D$ a regular ultrafilter on $\lambda$".  The
demand of regular is natural for several reasons.  The most relevant
is that using it Keisler proves that $\lambda^+$-saturation of
$M^\lambda/D$ depends only on the first order theory of $M$.  By
request we use a different symbol.

Naturally, we demand here $(\lambda,\theta)$-regularity because to
preserve the $\bbL_{\theta,\theta}$-theory we need the ultrafilter to
be $\theta$-complete, so the strongest possible regularity is for
$(\lambda,\theta)$.  Also the choice of saturation is natural.
\end{remark}

\noindent
We now turn to generalizing $\trianglelefteq^*$.
\begin{definition}  
\label{a9}
Assume $\bar\mu = (\mu_1,\mu_2)$, $\bar \chi = (\chi_1,\chi_2)$ 
and $\lambda \ge \theta$, $\mu_1 \ge 
\mu_2 \ge \theta$; if $\mu_1 = \mu$, $\mu_2 = \theta$ we may write
$\mu$ instead of $\bar\mu$; similarly for $\bar\chi$; if $\bar\chi =
(\mu,\theta)$ then we may omit $\bar\chi$.

\noindent
1) We say $T$ is [locally/fully] $(\lambda,\bar\mu,\theta)$-minimal
\underline{when}  for every complete $T_0 \supseteq T$ with $\tau(T_0)
\setminus \tau(T)$ of cardinality $\le \lambda$, for some $T_1$, we
have:
\mn
\begin{enumerate}
    \item[(a)]  $T_1 \supseteq T_0$ is a complete theory in $\bbL_{\theta,\theta}(\tau_{T_1})$.
\sn
    \item[(b)]  $T_1$ has no model of cardinality $< \theta$.
\sn
    \item[(c)]  $\tau(T_0) \subseteq \tau(T_1)$ and 
    $|\tau(T_1) \setminus \tau(T_0)| \le \lambda$.
\sn
    \item[(d)]  If $M_1$ is a model of $T_1$ of cardinality 
    $> \mu_2$ then $M_1 \rest \tau_T$ is [locally/fully] 
    $(\mu^+_1,\mu^+_2,\bbL_{\theta,\theta})$-saturated.
\end{enumerate}

\mn
2) For complete $T_1,T_2$ with no model of cardinality $< \theta$, we 
say $T_1 \blacktriangleleft^*_{\lambda,\bar\mu,\bar\chi,\theta} T_2$ \underline{when} 
for every complete $T^+_1 \supseteq T_1$ such that $|\tau(T^+_1)
\setminus \tau(T_1)| \le \lambda$ for some $T_3,\tau_2$ we have:
\mn
\begin{enumerate}
\item[(a)]  $T_3$ is a complete theory in
  $\bbL_{\theta,\theta}(\tau(T_3))$ 
\sn
\item[(b)]  $|\tau(T_3) \setminus \tau(T^+_1)| \le \lambda$ and 
$\tau(T_1) \subseteq \tau(T^+_1) \subseteq \tau(T_3)$
\sn
\item[(c)]  $T^+_1 \subseteq T_3$
\sn
\item[(d)]  $\tau_2 \subseteq \tau(T_3)$ and 
$T_3 \rest \tau'_2$ is isomorphic to $T_2$ over $\tau(T_1)$,
(if $\tau(T^+_1) \cap \tau(T_2) = \varnothing$ we can demand 
$T^+_1 \cup T_2 \subseteq T_3$; so the isomorphism above maps
$\tau'_2$ onto $\tau(T_2)$, $T_3 \rest \tau_2$ onto $T_2$, preserving
the number of places and being a predicate/function symbol) and is
  the identity on $\tau(T_1)$
\sn
\item[(e)]  if $M_3$ is a model of $T_3$ and $M_3 \rest \tau_2$ is
  locally $(\mu^+_1,\mu^+_2)$-saturated then $M_3 \rest \tau(T_1)$ is locally
$(\chi^+_1,\chi^+_2)$-saturated.
\end{enumerate}
\mn
3) We define $T_1 \blacktriangleleft^{*,\ful}_{\lambda,\bar\mu,\theta}
T_2$ is as in part (2) omitting the ``locally".

\noindent
4) In part (2), if we omit $\bar\mu,\bar\chi$
we mean $\|M_3\|$, i.e. $T_1 \blacktriangleleft^*_{\lambda,\theta} T_2$
means as above but we replace clause (e) in part (2) by:
\mn
\begin{enumerate}
    \item[(e)$'$]  if $M_3$ is a model of $T_3$ and $M_3 \rest \tau'_2$ is locally $(\|M_3\|,\|M_2\|)$-saturated \underline{then}  $M_3 \rest \tau'_1$ is locally $(\|M_3\|,\|M_3\|)$-saturated.
\end{enumerate}
\end{definition}

\begin{remark}\label{a9d}
0) We may note that $\trianglelefteq^*$ 
is defined similarly in the first order case.

\sn
1) Why the $T_0$ in \ref{a9}(1) and $T^+_1$ in \ref{a9}(2) 
in the definition?  Because otherwise it is not clear
why those relations are partial orders because $\bbL_{\theta,\theta}$ fail
the Robinson lemma, i.e. if $T_\ell \subseteq
\bbL_{\theta,\theta}(\tau_\ell)$ is complete for $\ell=1,2$ and
$\tau_0 = \tau_1 \cap \tau_2,T_1 \cap \bbL_{\theta,\theta}(\tau_0) =
T_2 \cap \bbL_{\theta,\theta}(\tau_0)$ \underline{then}  $T_1 \cup T_2$ does not
necessarily have a model); see \cite{BF}.

\sn
2) We may be worried that this will cause $\neg(T_1
\blacktrianglelefteq^*_{\lambda,\bar\mu,\bar\chi,\theta} T_2)$ because of 
trivial reasons, i.e. because for some $T^+_1 \supseteq T_2$ there is
no $T_3$ satisfying clauses (a)-(d) of Definition \ref{a9}(2).  But
this is not the case because
\mn
\begin{enumerate}
\item[$\boxplus$]  if $T_\ell \subseteq
  \bbL_{\theta,\theta}(\tau_\ell)$ has a model of cardinality $\ge
  \theta$ for $\ell=1,2$ and $\tau_1 \cap \tau_2 = \varnothing$ \underline{then} 
  $T_1 \cup T_2$ has a model of cardinality $\ge \theta$.
\end{enumerate}
\mn
[Why?  Because by the compactness for $\bbL_{\theta,\theta}$ and the
downward LST property if $\lambda = \lambda^{< \theta} + |T_\ell|$
then $T_\ell$ has a model of cardinality $\lambda$.]

\sn
3) For $\bbL^1_\kappa$ it holds; see \S3.
\end{remark}

\begin{conclusion}  
\label{a12}
1) $\blacktriangleleft^*_{\lambda,\bar\mu,\theta},
\blacktriangleleft_{\lambda,\bar\mu,\theta}$ are partial orders
(as are the full versions).

\sn
2) In Definition \ref{a10} the choice of $M_1,M_2$ does not matter.

\sn
3) If $T_1 \blacktriangleleft^*_{\lambda,\bar\mu,\theta} T_2$ \underline{then}  $T_1
\blacktriangleleft_{\lambda,\bar\mu,\theta} T_2$; also for the full
versions.
\end{conclusion}

\begin{PROOF}{\ref{a12}}  
1) Easy.

\sn
2) By \ref{a8}.

\sn
3) By part (2).
\end{PROOF}

\begin{claim}  
\label{a14}
1) $\Th_{\bbL_{\theta,\theta}}((\theta,<))$ is a
$\blacktriangleleft^*_{\lambda,\bar\mu,\theta}$-maximal and a
$\blacktriangleleft_{\lambda,\bar\mu,\theta}$-maximal theory (so $\bar\chi
= (\mu,\theta)$, see beginning of Definition \ref{a9}).

\sn
2) $\Th_{\bbL_{\theta,\theta}}(\theta,=)$ is a 
$\blacktriangleleft^*_{\lambda,\bar\mu,\theta}$-minimal and
$\triangleleft_{\lambda,\mu,\theta}$-minimal theory.

\sn
3) $T$ is $(\lambda,\bar\mu,\theta)$-minimal, (see Definition
\ref{a9}(1)) iff $T$ is 
$\blacktriangleleft^*_{\lambda,\bar\mu,\theta}$-minimal.
\end{claim}

\begin{PROOF}{\ref{a14}}
1) Easy:  we never get even local saturation, recalling \ref{a9d}(2).

\sn
2) Easy: even the (full)
 $(\lambda^+,\lambda^+,\bbL_{\theta,\theta})$-saturated means just ``of
 cardinality $\ge \lambda^+$".

\sn
3) Easy, too, just read the definitions.
\end{PROOF}
\newpage

\section {The $\ncp$ and local minimality} \label{3}

\begin{definition}  
\label{a17}
1) We say $T$ has the $\theta$-n.c.p. \underline{when}  it fails the
$\theta$-c.p.  We say $T$ has the $\theta$-c.p. \underline{when} :
  some $\varphi = \varphi(\bar x_{[\eps]},
\bar y_{[\zeta]}) \in \bbL_{\theta,\theta}(\tau_T)$ so
 $\eps,\zeta < \theta$ is a witness of $\theta$-c.p., that is, 
for every $\partial < \theta$ there 
are a model $M$ of $T$ and $\Gamma$ such that:
\mn
\begin{enumerate}
\item[$(*)_{M,\Gamma,\theta,\partial}$]  $\bullet \quad \Gamma \subseteq 
\{\varphi(\bar x_{[\eps]},\bar b):\bar b \in {}^\zeta\! M\}$
\sn
\item[${{}}$]  $\bullet \quad |\Gamma| < \theta$
\sn
\item[${{}}$]  $\bullet \quad \Gamma$ is $(< \partial)$-satisfiable in $M$
\sn
\item[${{}}$]  $\bullet \quad \Gamma$ is not satisfiable in $M$.
\end{enumerate}
\mn
2) For $\eps < \theta$, if $\Delta \subseteq 
\Phi_{T,\eps} \defeq \{\varphi(\bar x_{[\eps]},\bar
   y_\varphi):\varphi \in \bbL_{\theta,\theta}(\tau_T)\}$ is of
   cardinailty $< \theta$ we define the $\spec(\Delta,T)$ as the
   set of cardinals $\partial < \theta$ such that $\partial \ge 2$ and
for some model $M$ of $T$ and 
sequence $\LL \varphi_\alpha(\bar x_{[\eps]},\bar
   y_{\varphi_\alpha}):\alpha < \partial\RR$ of members of
$\Delta$ and $\bar a_\alpha \subseteq M$ of length 
$\lh(\bar y_{\varphi_\alpha})$ for $\alpha < \partial$, the set 
$\{\varphi_\alpha(\bar x_{[\eps]},\bar a_\alpha):
\alpha < \partial\}$ is not realized in $M$ but any subset of
cardinality $< \partial$ is realized.

\sn
3) For $\varphi = \varphi(\bar x_{[\eps]},\bar y_{[\zeta]}) \in
\Phi_{T,\eps}$ let $\spec(\varphi,T) = \spec(\{\varphi\},T\}$.

\sn
4) We may replace $\Delta$ by a sequence listing its members (even
   with repetitions).
\end{definition}

\begin{observation}  
\label{a19}
1) $T$ has $\theta$-c.p. \underline{iff}  for some $\varphi,\spec(\varphi,T)$ is
   unbounded in $\theta$ \underline{iff}  for some $\eps < \theta$ and
   $\Delta \subseteq \Phi_{T,\eps}$ of cardinality $< \theta$
   the set $\spec(\Delta,T)$ is unbounded in $\theta$.

\sn
2) In the definition of ``the theory $T$ has the $\theta$-c.p.", of
``$S = \spec(\varphi,T)$" and of ``$S = \spec(\Delta,T)$" see
Definition \ref{a17}, the model $M$ does not matter; of course, for
$T$ a complete $\bbL_{\theta,\theta}$-theory.

\sn
3) If $\eps < \theta$ and 
$\Delta \subseteq \Phi_{T,\eps}$ has cardinality 
$< \theta$ \underline{then}  for some $\psi = \psi(\bar
   x_{[\eps]},\bar y_\psi)$ we have:
\mn
\begin{enumerate}
    \item[$(a)$]  $\spec(\Delta,T) \subseteq \spec(\psi,T)$; \underline{moreover}, they are equal.
\sn
    \item[$(b)$]  If $M \models T$ then $$\{\varnothing\} \cup \{\varphi(M,\bar a):\varphi(\bar x_{[\eps]},\bar y) \in \Delta \text{ and } \bar a \in {}^{\lh(\bar y)}\!M\} = \{\psi(M,\bar a):\bar a \in {}^{\lh(\bar y)}\!M\}$$ (well, assuming
$\|M\| > 1$).
\end{enumerate}
\end{observation}

\begin{PROOF}{\ref{a19}}  
1) Obviously, the second assertion implies the first and the third trivially
implies the first by part (3) so we are left with proving ``the first
implies the second".

For $\partial < \theta$, let $M,\Gamma$ be as in \ref{a17}(1) for
   $\partial$, so necessarily $|\Gamma| \ge \partial$, let $\Gamma_1
   \subseteq \Gamma$ be of minimal cardinality such that $\Gamma_1$ is
   not realized in $M$.  So $\partial \le |\Gamma_1| \in
   \spec(\varphi,T)$.

\sn
2) Read Definition \ref{a17}.

\sn
3) Use definition by cases as in \cite{Sh:c}, (see
\cite[Ch.II,\S(2.1),pg.29]{Sh:c} 
and \S2 here; it suffices to assume the theory $T$ 
has no model with just one element).  That is, let 
$\LL \varphi_i(\bar x_{[\eps]},\bar y_i) : i < i_*\RR$ list $\Delta$, 
$\zeta = \sup\{\lh(\bar y_i) : i < i_*\}$ so $\zeta < \theta$
and let $$\psi = \psi(\bar x_{[\eps]},\bar y_{[\zeta + i_*+1]})
= \bigwedge\limits_{i < i_*} \Big[ \big(y_{\zeta +i_*} = y_{\zeta +i} \wedge
\bigwedge\limits_{j < i} y_{\zeta + i_*} \ne y_{\zeta +j} \big)
\rightarrow \varphi(\bar x_{[\eps]},\bar y \rest \zeta_i)\Big].$$
Now check. 
\end{PROOF}

\noindent
For first order $T$, $\aleph_0$-$\cp = \fcp$ follows from unstability
(by \cite[Ch.II,\S2]{Sh:a} = \cite[Ch.II,\S2]{Sh:c}), but not so
here.  The most interesting part of \ref{a22} is \ref{a22}(4) as we
have many non-implications.

\begin{claim}  
\label{a22}
1) There is a 5-unstable $T$ with $\spec(\bbL(\tau_T),T) = \aleph_0$ 
which is 3-unstable (see Definition \ref{a17}(2); yes, here we use
$\Delta =$ the set of first order formulas).

\sn
2) There is a 1-unstable, definably stable $T$ which has 
the $\theta$-$\cp$.

\sn
3) Assume $M = (\lambda,E^M)$, $E^M$ an equivalence relation on $\lambda$
and $\lambda \ge \theta$, $T = \Th_{\bbL_{\theta,\theta}}(M)$, \underline{then} 
$T$ is 1-stable; and $T$ has the $\theta$-c.p. iff 
$$\theta = \sup\!\big\{(a/E^M):a \in M \text{ and } \theta > |a/E^M|\big\}.$$

\sn
4) If $T$ is $\theta$-n.c.p. and is 1-unstable, then it is definably stable.
\end{claim}

\begin{PROOF}{\ref{a22}}  
1) Let $T$ be the theory of $I$ for any dense linear order $I$
which is $\theta$-saturated (in the first order sense) with neither first nor
last member.  This is the $T_5$ of \ref{y45}(4).

\sn
2) $T_0 = \Th((\theta,<))$ which by \ref{y45}(1) is 1-unstable,
definably stable; by inspection $\spec(\varphi,T) = \Card \cap \theta$
when $\varphi(x,y_0,y_1) = (x < y_1 \wedge x \ne y_0)$ so $T_0$ has the
$\theta$-c.p. 

\sn
3) Easy, too.

\sn
4) So we are assuming $T$ has the $\theta$-n.c.p. and is 1-unstable.
As $T$ is 1-unstable there is $\varphi(\bar x_{[\eps]},\bar
y_{[\eps]}) \in \bbL(\tau_T)$ witnessing it, hence we can choose:
\mn
\begin{enumerate}
\item[$(*)_1$]
\begin{enumerate}
\item[(a)]  a model $M$ of $T$ and $\bar a_\alpha \in {}^\eps\! M$ such
that
\sn
\item[(b)]   $M \models \varphi[\bar a_\alpha,\bar
  a_\beta]^{\iif(\alpha < \beta)}$ for $\alpha < \beta < \theta$
\sn
\item[(c)]  without loss of generality $M$ and $T$ has cardinality $\theta$
\sn
\item[(d)]  $\varphi(\bar x_{[\eps]},\bar y_{[\eps]} \vdash \neg
  \varphi(\bar y_{[\eps]},\bar x_{[\eps]})$.
\end{enumerate}
\end{enumerate}
\mn
By $\theta$ being a compact cardinal and $M \in \Mod_T$, every $p \in \bfS_\varphi(M)$ being definable because $T$ is definably stable, we can
find:
\mn
\begin{enumerate}
    \item[$(*)_2$]  $\psi = \psi(\bar y_{[\zeta]},\bar z_{[\xi]}) \in \bbL(\tau_T)$ such that: if $M \models T$ and $p \in \bfS_\varphi(M)$ then for some $\bar c \in {}^\xi M$ we have: if $\bar
b \in {}^\zeta\! M$ then $\varphi(\bar x_{[\eps]},\bar b) \in p$ \underline{iff} $M \models \psi[\bar b,\bar c]$
\sn
\item[$(*)_3$]
\begin{enumerate}
\item[(a)]  $\Delta = \{\varphi(\bar x_{[\eps]},\bar
  y_{[\eps]}),\varphi^\perp(\bar x_{[\eps]},\bar y_{[\eps]})\}$ see
  Definition \ref{y3}(2)
\sn
\item[(b)]  let $\partial = \partial_\Delta$ be $< \theta$ but $>
\sup[\spec(\Delta_\ell,\tau)]$ for $\ell=1,2$, see Definition
\ref{a17}(2).
\end{enumerate}
\end{enumerate}
\mn
Let
\mn
\begin{enumerate}
\item[$(*)_4$]  $\LL \bar c_\xi:\xi < \theta\RR$ list ${}^\xi
  \mu$ each appearing $\theta$-times
\sn
\item[$(*)_5$]  let $S = \{\delta < \theta:\cf(\delta) > \partial\}$.
\end{enumerate}
\mn
Now fix $\delta \in S$ for a while, we choose $\bar b_{\delta,\alpha}$
by induction on $\alpha < \theta$ such that:
\mn
\begin{enumerate}
\item[$(*)_6$]
\begin{enumerate}
\item[(a)]  $\bar b_{\delta,\alpha} \in {}^\eps\! M$
\sn
\item[(b)]  $M \models \varphi[\bar b_\beta,\bar b_{\delta,\alpha}]$
  for $\beta < \delta$
\sn
\item[(c)]  $M \models \varphi[\bar b_{\delta,\alpha},
\bar b_{\delta,\beta}]$ for $\beta < \alpha$
\sn
\item[(d)]  if possible (under (a)+(b)+(c)) then we have $M \models 
\psi[\bar b_{\delta,\alpha},\bar c^*_\alpha]$.
\end{enumerate}
\end{enumerate}
\mn
We can carry the induction, because for $\bar b$ to satisfy clauses
(a),(b),(c) it has to realize a $\Delta$-type $p_{\delta,\alpha}$ and
every member is satisfied by $\bar a_\beta$ for $\beta < \alpha$ large
enough, so recalling $\cf(\delta) > \partial$ and the choice of
$\partial$, we can carry the induction indeed; where
$p_{\delta,\theta} = \{\varphi(\bar a_\alpha,\bar x),\varphi(\bar
x,a_{\delta,\beta}):\alpha < \delta,\beta < \theta\}$ is a type in
$M$.  Hence there is $q_\delta \in \bfS(M)$ extending it.a

Now by the choice of $\psi$, there is $\bar d_\delta \in {}^\xi M$
such that:
\mn
\begin{itemize}
\item  $\bar b \in {}^\zeta\! M \Rightarrow [M \models \psi[\bar b,\bar
  c_\delta]$ iff $\varphi(\bar x,\bar b) \in p_\delta]$.
\end{itemize}
\mn
Clearly there is $\alpha(\delta) < \theta$ such that $\bar
c_{\alpha(\delta)} = \bar d_\delta$ hence
\mn
\begin{itemize}
\item  $r_\delta = p_{\delta,\alpha(\delta)}(\bar x_{[\eps]}) \cup
  \{\neg \psi(\bar x_{[\eps]},\bar c_{\alpha(\delta)})\}$ is
  contradictory, but of course
\sn
\item  every subset of $r_\delta$ with $< \cf(\delta)$ members is
  realized.
\end{itemize}
\mn
So $r_\delta$ contradicts ``$T$ has the $\theta$-n.c.p.
\end{PROOF}

\noindent
More generally
\begin{claim}  
\label{a23}
Assume $T = \Th_{\bbL_{\theta,\theta}}(M),M$ a 
$\theta$-saturated model (in the first order sense)
with $\Th_{\bbL}(M)$, the first order theory of $M$, being unstable
(e.g. random graph).

\sn
1) $T$ is 5-unstable.

\sn
2) $T$ has $\theta$-$\ncp$ provided that $\theta =
  \sup\{\theta':\theta' < \theta$ is a compact cardinal$\}$.

\sn
3) $T$ has the $\theta$-$\cp$ \underline{when}  $(a)$ and $(b) \vee (b)' \vee (b)''$ where:
\mn
\begin{enumerate}
\item[$(a)$]  the first order theory $\Th_{\bbL}(M)$ has the
  independence property (hence is unstable)
\sn
\item[$(b)$]  for some $\kappa < \theta$ we have
$\theta = \sup\{\mu$: there is a graph $G$ on $\mu$ such
  that $\chr(G) > \kappa$ but $A \in [\mu]^{< \mu} \Rightarrow \chr(G \rest A)
  \le \kappa\}$
\end{enumerate}
\mn
(maybe $(b)',(b)''$ are more transparent)
\mn
\begin{enumerate}
\item[$(b)'$]  $\theta = \sup\{\mu:\mu = \cf(\mu) < \theta$ and some
  stationary $S \subseteq S^\mu_{\aleph_0}$ does not reflect$\}$ or
just
\sn
\item[$(b)''$]  like (b) 
replacing $\aleph_0$ by some regular $\kappa < \theta$.
\end{enumerate}
\mn
4) $T$ has the $\theta$-$\cp$ \underline{when}  (a) and (b) $\vee$ (b)$'$ where:
\mn
\begin{enumerate}
\item[$(a)$]  the first order theory $\Th_{\bbL}(M)$ has the
  strict order property (hence is unstable)
\sn
\item[$(b)$]  for some regular $\kappa < \theta$ we have $\theta =
\sup\{\mu^{< \kappa}:\mu = \cf(\mu)$ and $I^\kappa/D$ has a
  $(\mu,\mu)$-cut for some ultrafilter $D$ on
  $\kappa$ and $\theta$-saturated dense linear order $I\}$, we can fix
$D$ and $I$; see Golshani-Shelah \cite[Th.3.3]{Sh:1075}
\end{enumerate}
\mn
(maybe more transparently)
\mn
\begin{enumerate}
\item[(b)$'$]  for some regular $\kappa < \theta$ we have $\theta =
\sup\{\mu^{< \kappa}:\mu$ is a successor cardinal, $\mu = \mu^{< \mu}
> \kappa^+$ and there are a
stationary $S \subseteq S^\mu_\kappa$ and $\bar C = \LL
C_\delta:\delta < \mu$ limit$\RR$ such that $C_\delta$ is a closed
unbounded subset of $\delta$ disjoint to $S$ and $\delta_1 \in
C_{\delta_2} \Rightarrow C_{\delta_1} = C_{\delta_2} \cap \delta_1\}$.
\end{enumerate}
\mn
5) $T$ has the $\theta$-$\ncp$ if $\Th_{\bbL}(M)$ is stable.
\end{claim}

\begin{remark}
\label{a23d}
1) Recall that a first order $T_0$ is unstable iff it has the
independence property or the strict order property, hence part
(3),(4),(5) of \ref{a23} covers all complete first order $T$.

\noindent
2) The statements in \ref{a23}(3)$(b)'$, \ref{a23}(4)$(b)'$ are consistent
by a relative of Laver indestructibility; see,
e.g. \cite[1.3=La7]{Sh:945}.

Note that \cite[Th.3.3]{Sh:1075} use conditions weaker than
\ref{a23}(4)(b)$'$, because by \cite{Sh:922} the assumptions on $\mu$
and $\kappa$ implies $\diamondsuit_S$.
\end{remark}

\begin{PROOF}{\ref{a23}}
1) Let $\varphi(\bar x,\bar y) \in \bbL(\tau_T)$ be a first order
   formula which has the order property for $T$.  Easily it witnesses
   that $T$ is 5-unstable.

\noindent
2) Easy, but we shall elaborate.

So let $\varphi = \varphi(\bar x,\bar y) \in
\bbL_{\theta,\theta}(\tau_T)$ be a formula and we shall prove that
$\spec(\varphi,T)$ is bounded in $\theta$.  As $\theta$ is strongly
inaccessible there is $\sigma < \theta$ such that $\varphi \in
\bbL_{\sigma,\sigma}(\tau_T)$ so $\lh(\bar x) + \lh(\bar y) <
\sigma$.  By the assumption without loss of generality $\sigma$
is a compact cardinal.  Now for every cardinal $\partial \in
[\sigma,\theta)$ and $\tau_M$-model $N$ consider the statement
\mn
\begin{enumerate}
\item[$(*)^+_{N,\varphi,\partial}$]  if $\bar b_i \in {}^{\lh(\bar
  y)}\!N$ for $i < \partial$ and every subset of $p(\bar x) \defeq  
\{\varphi(\bar x,\bar b_i):i < \partial\}$ of cardinality $<
\partial$ is realized in $N$ \underline{then}  $p(\bar x)$ is realized in $N$.
\end{enumerate}
\mn
Now first it suffices to prove $(*)^+_{M,\varphi,\partial}$ for every
such $\partial$ because this statement can be phrased as a sentence
$\psi_{\varphi,\partial}$ in $\bbL_{\theta,\theta}(\tau_T)$ and it
means $\partial \notin \spec(\varphi,T)$.

Second, assume the antecedent of $(*)^+_{M,\varphi,\partial}$ so $\LL \bar
b_i:i < \partial \RR$ are as above, let $B = \bigcup\{\bar b_i:i <
\partial\}$ hence $p$ is a $(< \sigma)$-satisfiable 
$\lh(\bar x)$-type in $M$ over $B,B \subseteq M,|B| = \partial$.
Hence there is an
$\bbL_{\sigma,\sigma}(\tau_T)$-complete type $q(\bar x)$ in $\bfS^{\lh(\bar x)}_{\bbL_{\sigma,\sigma}(\tau_T)}(M)$ extending it;
the existence of $q(\bar x)$ is the point at which we use ``$\sigma$
is a compact cardinal".

Let $q'(\bar x)$ be the set of first order formulas in $q(\bar x)$ so
clearly $q'(\bar x) \in \bfS^{\lh(\bar x)}_{\bbL}(M)$; as
$M$ is $\theta$-saturated clearly some 
$\bar a \in {}^{\lh(\bar x)}\!M$ realizes $q'(\bar x) \rest B$.  We are done because in $M$ every
$\bbL_{\sigma,\sigma}(\tau_T)$ formula is equivalent to a Boolean
combination of first order formulas.  In other words, without loss of generality $M$
has elimination of quantifiers for first order formulas; and it
follows that it has elimination of quantifiers also for
$\bbL_{\sigma,\sigma}(\tau_T)$; so we are done.

\noindent
3) Trivially $(b)' \Rightarrow (b)''$ and by 
\cite[1.2=La6]{Sh:1006} we have 
$(b)'' \Rightarrow (b)$ so we can assume (a) + (b).

Let $\varphi(\bar x_{[m]},\bar y_{[n]}) \in \bbL(\tau_T)$ be a
first-order formula with the independence property for
$\Th_{\bbL}(M)$.   Define $\psi(\bar x_{[\kappa]},
\bar y^0_{[n]},\bar y^1_{[n]}) \in
\bbL_{\kappa^+,\aleph_0}(\tau_T)$ or pedantically $\in
\bbL_{\kappa^+,\kappa^+}(\tau_T)$ as saying:
\mn
\begin{enumerate}
\item[$(*)_1$]  for each $\ell \in \{0,1\}$ there is a 
unique $i_\ell < \kappa$ such that 
$\varphi(\bar x_{[m i_\ell,m(i_\ell+1))},\bar y^\ell_{[n]})$ 
and moreover $i_0 \ne i_1$.
\end{enumerate}
\mn
We claim $\sup(\spec_\psi(T)) = \theta$. By clause (b), for some
unbounded $\Theta \subseteq \Card \cap \theta$ for every $\mu \in
\Theta$ there is a graph $G_\mu$ with set of nodes $\mu$ such that
$\chr(G_\mu) > \kappa$ but $u \in [\mu]^{< \mu}$ implies $\chr(G_\mu
\rest u) \le \kappa$.  Since $\varphi$ has the independence property and
$M$ is (first-order) saturated, we can find $\LL \bar b_i:i <
\mu\RR$ with $\bar b_i \in {}^n M$ such that for every $\bar{\bft} \in 
{}^\mu 2$ there is $\bar a \in {}^m M$ with $\bigwedge\limits_{i
  < \mu} \varphi^M[\bar a,\bar b_i]^{\iif(\bft(i))}$.

Now let:
\mn
\begin{enumerate}
\item[$(*)_2$]  $\Gamma_\mu = 
\{\psi(\bar x,\bar b_i,\bar b_j):i < j<\mu \text{ and }
(i,j) \in \edge(G_\mu)\}$.
\end{enumerate}
\mn
Easily
\mn
\begin{enumerate}
\item[$(*)_3$]   $\Gamma_\mu \text{ demonstrates } \mu \in \spec_\psi(T)$.
\end{enumerate}
\mn
Let $I$ be as there and let $D$ be a uniform ultrafilter on $\kappa$ such that
$\Theta$ is unbounded in $\theta$ where

\[
\Theta = \{\mu:\mu = \mu^{<\kappa} \text{ and in } I^\kappa/D 
\text{ there is a } (\mu,\mu)\text{-cut}\}.
\]

\mn
Let $\mu \in \Theta$; let the first order formula $\varphi =
\varphi(\bar x_{[n]},\bar x_{[m]})$ exemplify that $\Th_{\bbL}(M)$ has
the strict order property.  For notational simplicity assume $n=1=m$.
We choose $a_s \in {}^m M$ for $s \in I$ such that $M \models (\forall
x)(\varphi(x,a_s) \rightarrow \varphi(x,a_2))$ \underline{iff}  $s <_I t$.  

By the choice of $\mu$, there are $f^1_\alpha,f^2_\alpha \in {}^\kappa
I$ such that in $I^\kappa/D$ we have $\alpha < \beta < \mu \Rightarrow
f^1_\alpha/D < f^1_\beta/D < f^2_\beta/D < f^2_\alpha/D$, but
$I^\kappa/D$ omits the type $p = \{f^1_\alpha/D < x <
f^2_\alpha/D:\alpha < \mu\}$.  By \cite[Lemma 2.1]{Sh:1075} if $J$
is the completion of $I$ then also $J^\kappa/D$ omits the type $p$.

Let $\psi(\bar x_{[\kappa]},\bar y_{[\kappa]},\bar z_{[\kappa]})$ be
the formula $\bigvee\limits_{A \in D} \, \bigwedge\limits_{i \in
  A}(\varphi(x_i,z_i \wedge \neg \varphi(x_i,y_i)$.

We define $\bar b^\ell_\alpha = \LL b^\ell_{\alpha,\eps}:\eps <
\kappa\RR$ for $\alpha < \mu,\ell \in \{1,2\}$ by
$b^\ell_{\alpha,\eps} = a_{f^\ell_\alpha(\eps)} \in M$.

Now let $\Gamma_\mu = \{\psi(\bar x,\bar b^1_\alpha,\bar
b^2_\alpha):\alpha < \mu\}$ and the rest should be clear.

\noindent
4) Clause (b)$'$ implies clause (b) is proved in 
Golshani-Shelah \cite[Th.3.3]{Sh:1075}.  So we can assume (a) + (b)
and the proof is similar to the proof of part (2).

\noindent
5) Without loss of generality $\tau(T)$ has cardinality $< \theta$.  Assume
$\eps < \theta,\varphi(\bar x_{[\eps]},\bar y) \in
\bbL_{\theta,\theta}(\tau_T)$, let $\zeta = \lh(\bar y)$ 
and $\Gamma = \{\varphi(\bar x_{[\eps]},\bar a_\alpha):
\alpha < \alpha_* < \theta\}$ is a set of
$\bbL_{\theta,\theta}$-formulas with parameters from $M$.
Without loss of generality $\LL \bar a_\alpha:\alpha < \alpha_*\RR$ is with no
repetitions, we let $\kappa = (|T| + |\zeta|)^{|T|+|\eps|}$.  

We shall use freely:
\mn
\begin{enumerate}
\item[$(*)$]  if $\alpha < \alpha'_*$ and $\bar b',\bar b'' \in
  {}^\delta M$ realize the same first order type over $\bar a_\alpha$
  then $M \models \varphi[\bar b',\bar a_\alpha] \equiv \varphi[\bar
  b'',\bar a_\alpha]$.
\end{enumerate}
\mn
We shall assume $\Gamma$ is $(\le 2^\kappa)$-satisfiable in $M$ and 
prove that it is satisfiable in $M$; this easily suffices.  
Let $A = \bigcup\{\bar a_\alpha:\alpha < \alpha_*\}$ and we try by
induction on $i < \kappa^+$ to choose $M_i \prec_{\bbL} M$ 
of cardinality $\le 2^\kappa$, increasing
   continuous with $i$ such that: if $p(\bar x_{[\eps]}) 
\in \bfS^\eps_{\bbL}(M_i \cup A)$ does not fork over $M_i$ then for
some $\alpha < \alpha_*,\bar a_\alpha \subseteq M_{i+1}$ and
   $p(\bar x_{[\eps]}) \nvdash \varphi_\alpha(\bar
x_{[\eps]},\bar a_\alpha)$.  If we are stuck in $i$,
i.e. $M_i$ is well defined but we cannot choose $M_{i+1}$, then as
$[p_1,p_2 \in \bfS^\eps_{\bbL}(M_i \cup A)$ does not fork over
$M_i \Rightarrow (p_1 = p_2 \Leftrightarrow p_1 \rest M_i = p_2 \rest
M_i)]$ and $\bfS^\eps_{\bbL}(M_i)$ has cardinality $(\sup_n|\bfS^n_{\bbL}(M_i)|)^{|\eps|} \le (2^\kappa)^{|\eps|} =
2^\kappa$, clearly for some $p(\bar x) \in 
\bfS^\eps_{\bbL}(M_i \cup A)$ not forking over $M_i$
there is no such $\alpha$, but
$p(\bar x)$ is realized in $M$ hence so is $\Gamma$.

What if we succeed to carry the induction?  Choose $\bar b$ which
realizes $\Gamma' = \{\varphi(\bar x_{[\eps]},\bar
a_\alpha):\bar a_\alpha \subseteq M_i$ for some $i < \kappa^+\}$, now
$\{\alpha < \alpha_*:\bar a_\alpha \subseteq M_{\kappa^+}\} \le
\|M_{\kappa^+}\|^{|\zeta|} \le 2^\kappa$, hence $\Gamma'$ indeed is
realized in $M$ say by $\bar b \in {}^\eps\! M$ and let 
$q \in \bfS^\eps_{\bbL}(M_{\kappa^+} \cup A)$ extend
$\tp_{\bbL}(\bar b,M_{\kappa^+},M)$ and does not fork over
$M_{\kappa^+}$.  Without loss of generality $\bar b$ realizes $q$ in $M$ using a
partial automorphism of $M$.

Now for every $i < \kappa^+$, by the induction $\tp_{\bbL}(\bar b,M_\kappa
\cup A)$ is not a non-forking extension of $\tp(\bar b,M_i) = p$ hence
also $\tp(\bar b,M_\kappa)$ is not.  Contradiction to 
``$\Th_{\bbL}(M)$ is stable".
\end{PROOF}

\begin{claim}  
\label{a25}
The model $N = M^I/D$ is not $(\chi^+,\theta,\bbL_{\theta,\theta})$-saturated
(even locally, and even just for $\varphi$-types) \underline{when} :
\mn
\begin{enumerate}
\item[$(a)$]  $D \in \uf_\theta(I)$
\sn
\item[$(b)$]  $\varphi(\bar x_{[\eps]},\bar y_{[\zeta]})$ witnesses $T$
  has the $\theta$-c.p.
\sn
\item[$(c)$]  $\chi = \lcr_\theta(\spec(\varphi,T),D)$ see \ref{w8}(3),
  equivalently letting 
\newline
$(J,<_J,P^J) = (\theta,<,\spec(\varphi,T))^I/D$ we have
\newline
$\chi = \min\{|\{s:s <_J t\}|:t \in P^J$, but $(\exists^{\ge
\theta}s)(s <_J t)\}$.
\end{enumerate}
\end{claim}

\begin{PROOF}{\ref{a25}}  
Straightforward or see the proof of \ref{a46} below.
\end{PROOF}

\begin{remark}
In \ref{a26}, \ref{a28} + more below the distinction $T,T_1$ is not
necessary.  But it is natural in the way we shall quote them; that is we
consider properties of $T$ and ask for $T_1 \supseteq T$ large enough
such that ``$M \models T_1 \Rightarrow M \rest \tau_T$ satisfies ..."
\end{remark}

\begin{definition}
\label{a26}
We say that $(\varphi,M,\bar{\bfa},\bar{\bfb})$ strongly
$\chi$-witnesses or $(M,\bar{\bfa},\bar{\bfb})$ strongly
$(\chi,\varphi)$-witness that $T$ is 1-unstable \underline{when}  for some $T_1
\supseteq T$: (if $\chi = \theta$ we may omit it)
\mn
\begin{enumerate}
\item[$\circledast_1$]  
\begin{enumerate}
\item[(a)]  $M$ is a model of $T_1$
\sn
\item[(b)]  $\varphi = \varphi(\bar x_{[\eps]},\bar
  y_{[\zeta]}) \in \bbL_{\theta,\theta}(\tau(T_1))$
\sn
\item[(c)]
\begin{enumerate}
\item[($\alpha$)]  $\bar a_\alpha \in {}^\eps\! M,
\bar b^1_\beta \in {}^\zeta\! M$ for $\alpha,\beta < \chi$ are
such that $M \models 
\varphi[\bar a_\alpha,\bar b_\beta]^{\iif(\alpha < \beta)}$
\sn
\item[($\beta$)]  $\bar{\bfa} = \LL \bar a_\alpha:
\alpha < \chi\RR$ and $\bar{\bfb} = \LL 
\bar b_\alpha:\alpha < \chi\RR$
\end{enumerate}
\sn
\item[(d)]  for every $\bar a \in {}^\eps\! M$ for some truth 
value $\bft$ for every $\beta < \chi$ large enough
we have $M \models \varphi[\bar a,\bar b_\beta]^{\iif(\bft)}$
\sn
\item[(e)]   for every $\bar b \in {}^\zeta\! M$
for some truth value $\bft$ for every $\alpha < \chi$ large
enough we have $M \models \varphi[\bar a_\alpha,\bar b]^{\iif(\bft)}$.
\end{enumerate}
\end{enumerate}
\end{definition}

\begin{remark}
\label{a26d}
Definition \ref{a26} is a case of ``$\LL \bar a^1_\alpha \caret
\bar b^1_\alpha:\alpha < \chi\RR$ is convergent", see
\cite[\S2,Def 2.1=L300a-2.1,pg.25]{Sh:300a}.
\end{remark}

\begin{observation}
\label{a27}
1) Assume the triple $(M,\bar{\bfa},\bar{\bfb})$ strongly
   $(\chi,\varphi)$-witnesses that $T$ is 1-unstable and $\chi =
\cf(\chi) \ge \theta$.  If $\lambda = \lambda^{< \theta} + |\tau_T|$ and
$\sigma = \cf(\sigma) \in [\theta,\lambda]$, \underline{then}  there is a triple
$(M',\bar{\bfa}',\bar{\bfb}')$ which strongly 
$(\sigma,\varphi)$-witness $T$ is 1-unstable and $\|M'\| = \lambda$.
We can add $\|M\| \le \lambda \Rightarrow M \prec_{\bbL_{\theta,\theta}} M'$
 and $\chi > \lambda \Rightarrow M' \prec_{\bbL_{\theta,\theta}} M$.

\noindent
2) If for every $\tau' \subseteq \tau(T)$ of cardinality $< \theta$
such that $\varphi \in \bbL_{\theta,\theta}(\tau')$ there is a
strong $(\chi,\varphi)$-witness for $T \cap
   \bbL_{\theta,\theta}(\tau)$ being 1-unstable for some $\chi =
   \cf(\chi) \ge \theta$ \underline{then}  there is a strong
   $(\chi,\varphi)$-witness for $T$ being 1-unstable for every $\chi =
   \cf(\chi) \ge \theta$.
\end{observation}

\begin{PROOF}{\ref{a27}}
1) First let $D \in \ruf_\theta(\lambda)$ and so by \ref{x31}(3)
for some $\chi_1 = \cf(\chi_1) \in [\lambda^+,2^\lambda)$ and 
$\bar{\bfa}',\bar{\bfb}'$, we have $(M^I/D,\bar{\bfa}',\bar{\bfb}')$
   strongly $(\chi_1,\varphi)$ witness $T$ is 1-unstable.  Now apply
   the downward LST argument.

\noindent
2) Easy, too.
\end{PROOF}

\begin{observation}
\label{a27c}
For any model $M$ satisfying $\|M\| = \|M\|^{< \theta}$ 
there is an expansion $M^*_1$ by the new
function symbols $F_\xi(\xi < \theta),F_\xi$ being $\xi$-place such
that $M' \equiv_{\bbL_{\theta,\theta}} M \Rightarrow \|M'\| = 
\|M'\|^{<\theta}$.
\end{observation}

\begin{PROOF}{\ref{a27c}}
Choose $F^{M_2}_\xi:{}^\xi M_2 \rightarrow M$ which is one-to-one.
\end{PROOF}

\begin{claim}
\label{a27g}
Assume $T \subseteq \bbL_{\theta,\theta}(T_1)$ is complete 1-unstable
theory as witnessed by $\varphi(\bar x,\bar y)$.

For any theory $T_1 \supseteq T$ and regular $\chi \ge \theta$ there
are $M,\bar{\bfa},\bar{\bfb}$ as in Definition \ref{a26} with $M
\in \Mod_{T_1}$.  
\end{claim}

\begin{PROOF}{\ref{a27g}}
Let $\lh(\bar x) = \eps < \theta,\lh(\bar y) = \zeta <
\theta$.

Let $P,<$ be new predicates, i.e. $\notin \tau(T_1)$ with $\eps
+ \zeta,\eps + \zeta + \eps + \zeta$ places
respectively and let $F_\xi$ be a new $\xi$-place function symbol.  

Let $T_2$ be the set of
$\bbL_{\theta,\theta}(\tau_{T_1} \cup \{P,<,F_\xi:
\xi < \theta\})$-sentences such that for any $\tau(T)$-model $M_2$ we have:
 $M_2 \models T_2$ iff
\mn
\begin{enumerate}
\item[$(*)_1$]  
\begin{enumerate}
\item[(a)]  $M_2 \models T_1$
\sn
\item[(b)]  $<^{M_2}$ linearly ordered $P^{M_2}$, of
  cofinality $\ge \theta_1$ for any $\theta_1 < \theta$
\sn
\item[(c)]   if $\bar a_1 \caret \bar b_1 \in P^{M_2},\bar
  a_2 \caret \bar b_2 \in P^{M_2},\bar a_\ell \in
  {}^\eps(M_2),\bar b_\ell \in {}^\zeta(M_2)$ for $\ell=1,2$ and
$\bar a_1 \caret \bar b_1 <^{M_2} \bar a_2 
\caret \bar b_2$ \underline{then}  $M_2 \models 
\varphi(\bar a_1,\bar b_2) \wedge \neg \varphi(\bar a_2,\bar b_1)$
\sn
\item[(d)]  for every $\bar a' \in {}^\eps(M_2)$ for some 
truth value $\bft$, for every $\bar a \caret \bar b \in P^{M_2}$
which is $<^{M_2}$-large enough 
(and $(\lh(\bar a),\lh(\bar b)) = (\eps,\zeta)$, of course)
we have $M_2 \models \varphi[\bar a',\bar b]^{\iif(\bft)}$
\sn
\item[(e)]    for every $\bar b' \in {}^\zeta(M_2)$ 
for some truth value $\bft$, for every $\bar a \caret \bar b 
\in P^{M_2}$ which is $<^{M_2}$-large enough, we have
$M_2 \models \varphi[\bar a,\bar b']^{\iif(\bft)}$.
\end{enumerate}
\end{enumerate}
\mn
Now
\mn
\begin{enumerate}
\item[$(*)_2$]  $T_2$ is an $\bbL_{\theta,\theta}$-theory.
\end{enumerate}
\mn
Why?  For this it suffices to prove that every $T'_2 \subseteq T_2$
  of cardinality $< \theta$ has a model, so without loss of generality $|\tau_{T_1}| < \theta$
  and let $M_1 \models T_1$.  As $T$ is complete 1-unstable as
witnessed by $\varphi$ for every $\gamma < \theta$ there are
$\LL (\bar a^\gamma_i,\bar b^\gamma_i):i < \gamma\RR$ in
$M_1$ as in Definition \ref{y2}(1), i.e. 
$M_1 \models \varphi[\bar a^\gamma_i,\bar b^\gamma_j]^{\iif(i<j)}$ 
for $i,j < \gamma$.  

By compactness of $\bbL_{\theta,\theta}$ possibly changing $M_1$ we
have $\LL (\bar a_i,\bar b_i):i < \theta\RR$ as above.  By the
LST argument without loss of generality $\|M_1\| = \theta$, hence
$|{}^\eps(M_1)| + |{}^\zeta(M_1)| = \theta$.

Let $\LL \bar c_\alpha:\alpha < \theta\RR$ list
${}^\eps|(M_1)|$ and $\LL \bar d_\alpha:\alpha <
\theta\RR$ list ${}^\zeta|(M_1)|$. 

We define $f:[\theta]^3 \rightarrow \{0,1\}$ by:
\mn
\begin{enumerate}
\item[$(*)_3$]  if $\alpha < \beta < \gamma < \theta$ then
  $f(\{\alpha,\beta,\gamma\}) = 1$ \underline{iff}  $j < \alpha \Rightarrow M_1
\models ``\varphi[\bar c_j,\bar b_\beta] \equiv \varphi[\bar c_j,\bar
  b_\gamma]"$ and $j < \alpha \Rightarrow M_1 \models ``\varphi[\bar
  a_\beta,\bar d_j] \equiv \varphi[\bar a_\gamma,\bar d_j]"$.
\end{enumerate}
\mn
But $\theta$ is, of course, weakly compact so $f$ is constant on
$[\cU]^3$ for some $\cU \in [\theta]^\theta$; easily necessarily $f$
is constantly 1.

We now define $M_2$ expanding $M_1$ by

\[
P^{M_2} = \{\bar a_\alpha \caret \bar b_\alpha:\alpha \in \cU\}
\]

\[
<^{M_2} = \{\bar a_\alpha \caret \bar b_\alpha \caret \bar a_\beta
\caret \bar b_\beta:\alpha < \beta \text{ are from } \cU\}.
\]

\mn
Easily $M_2 \models T'_2$ hence we are done proving $(*)_2$.
\mn
\begin{enumerate}
\item[$(*)_4$]  for every $\lambda$ there is a model $M_2$ of $T_2$
  such that $\cf(P^{M_2},<^{M_2}) \ge \lambda^+$.
\end{enumerate}
\mn
[Why?  Let $M_2 \models T_2,D \in \ruf_{\chi,\theta}(\lambda)$ then
  $(M_2)^\lambda/D$ is as required by \ref{x31}(3).]
\mn
\begin{enumerate}
\item[$(*)_5$]  for every regular $\chi \ge \theta$ and $\lambda =
  \lambda^{< \theta} + |T_1| + \chi$ there is a model $M_2$ of $T_2$
of cardinality $\lambda$ such that $\cf(P^{M_2},<^{M_2}) = \chi$.
\end{enumerate}
\mn

[Why?  By $(*)_4$ applied with $((\chi + \lambda + \theta)^{<
  \theta})^+$ here standing for $\lambda$ there and then use the LST argument.]

To finish note that
\mn
\begin{enumerate}
\item[$(*)_6$]  if $M_2 \models T_2$ and $\LL (\bar a_\alpha \char
  94 \bar b_\alpha):\alpha < \chi\RR$ is $<^{M_2}$-increasing
  cofinal in $P^{M_2}$ and $(\lh(\bar a_\alpha),
\lh(\bar b_\alpha)) = 
(\eps,\zeta)$ \underline{then}  $(\varphi,M_2,\LL \bar
  a_\alpha:\alpha < \chi\RR,\LL \bar b_\alpha:\alpha <
  \chi\RR)$ is as required in Definition \ref{a26}.
\end{enumerate}
\mn
[Why?  Read the Definition of $T_2$.]
\end{PROOF}

\begin{remark}
\label{a27k}
1) We can strengthen the conclusion of \ref{a27g} to
\mn
\begin{enumerate}
\item[$(*)$]  for every $\bar d \in {}^{\theta >}\mu$ the sequence
  $\LL \tp_{\bbL_{\theta,\theta}(\tau)}(\bar a^1_\alpha \caret
  \bar a^2_\alpha,\Rang(\bar d),M):\alpha < \chi\RR$ is eventually
  constant.
\end{enumerate}
\mn
How?  In $(*)_3$ we can change somewhat the demand:
\mn
\begin{enumerate}
\item[$(*)'_3$]  for $\alpha < \beta < \gamma < \theta$ then
  $f(\{\alpha,\beta,\gamma\})=1$ \underline{iff}  for every $j < \alpha$ and
  formula $\vartheta(\bar x_{[\eps + \zeta]},\bar y_{[\eps +
    \zeta]})(\tau(T'_2))$ we have $M_1 \models \vartheta[a^1_\beta \char
  94 \bar a^2_\beta,\bar c_j] \Leftrightarrow M_1 \models \vartheta[\bar
  a^1_\gamma,\bar a^2_\gamma,\bar c_j]$.
\end{enumerate}
\mn
We similarly change $(*)_1(c) +(d)$.

\noindent
2) Clearly if $T \vdash ``(P,<)$ is a linear order of cofinality $\ge
\partial"$ for every $\partial < \theta$ and $\lambda = \lambda^{<
  \theta} + |T| \ge \kappa = \cf(\kappa) \ge \theta$, then $T$ has a
model $N$ of cardinality $\lambda$ such that $\cf(P^N,<^N) = \kappa$.
This is proved inside the proof of \ref{a27g} and holds by \ref{x31}(3).
\end{remark}

\begin{claim}  
\label{a28}
If (A) then (B) where:
\mn
\begin{enumerate}
\item[(A)]
\begin{enumerate}
\item[(a)]  $T$ is a complete $\bbL_{\theta,\theta}(\tau_T)$-theory
\sn
\item[(b)]  $T$ is 1-unstable as witnessed by
$\varphi(\bar x_{[\eps]},y_{[\zeta]})$ and let
$\psi = \psi(\bar x_{[\zeta]},\bar y_{[\eps]}) = 
\varphi(\bar y_{[\eps]},\bar x_{[\zeta]})$
\sn
\item[(c)]  $T_1 \supseteq T$ is a complete
  $\bbL_{\theta,\theta}(\tau_1)$-theory and $|\tau(T_1) 
\setminus \tau(T)| \le \lambda$
\sn
\item[(d)]  $\bfx$ is a non-trivial $(\theta,\theta)$-$\luft$
\sn
\item[(e)]  $\chi = \cf(\uflp_{\bfx}(\theta,<))$
  hence $\chi = \chi^{< \theta}$, see \ref{x17} - \ref{x18}
\end{enumerate}
\sn
\item[(B)]  for some $M_1 \models T_1$ the model 
$\uflp_{\bfx}(M_1)$ is not $(\chi^+,\{\varphi\})$-saturated or not
$(\chi^+,\{\psi\})$-saturated, see Definition \ref{a5}(4).
\end{enumerate}
\end{claim}

\begin{PROOF}{\ref{a28}}  
\smallskip

\noindent
\underline{Case 1}:  $|T_1| \le \theta$.

Choose $D_* \in \ruf_{\chi,\theta}(\chi)$ hence $D_*$ is uniform.  
Let $(M,\LL
\bar a^1_\alpha:\alpha < \theta\RR,\LL \bar b^1_\alpha:\alpha
< \theta\RR)$ be a strong $\varphi$-witness for $T$ being
1-unstable, see Definition \ref{a26}, exists by Claim \ref{a27g}.

Let $M^+ = (M,P^{M^+},<^{M^+})$ where $P^{M^+} = \{\bar a^1_\alpha
\caret \bar b^1_\alpha:\alpha < \theta\}$ and $<^{M^+} = \{(\bar a^1_\alpha
\caret \bar b^1_\alpha,\bar a^1_\beta \caret \bar b^1_\beta):
\alpha < \beta < \theta\}$ and let $N^+ = \lup_{\bfx}(M^+)$ hence
clearly $N^+ = (\lup_{\bfx}(M),P^{N^+},<^{N^+})$ and $N = \lup_\bfx(M)$.  By clause (A)(e) of the claim, clearly 
$(P^{N^+},<^{N^+})$ is a linear order of cofinality $\chi$ so we can
choose an increasing cofinal sequence $\LL \bar a^3_\alpha \char
94 \bar b^3_\alpha:\alpha < \chi\RR$ in $(P^{N^+},<^{N^+})$, and 
by \ref{x7}
\mn
\begin{enumerate}
\item[$(*)_1$]  if $\bar a \in {}^\eps|N^+|$ and $\bar b \in
  {}^\zeta|N^+|$ \underline{then}  for some truth values $\bft(1),\bft(2)$ 
for every $\alpha < \chi$ large enough $N^+ \models
 ``\varphi[\bar a,\bar b^3_\alpha]^{\iif(\bft(1))} \wedge \varphi[\bar
  a^3_\alpha,\bar b]^{\iif(\bft(2))}"$; of course this is a property of $N$.
\end{enumerate}
\mn
We try to choose $(N_\alpha,\bar b^4_\alpha)$ by induction on
$\alpha < \chi$ such that:
\mn
\begin{enumerate}
\item[$(*)^2_\alpha$]  
\begin{enumerate}
\item[(a)]  $N_\alpha \prec_{\bbL_{\theta,\theta}} N$ has cardinality $\chi$
\sn
\item[(b)]   if $\beta < \alpha$ \underline{then}  $\bar a^3_\beta
  \caret \bar a^4_\beta \subseteq N_\alpha \subseteq N$
\sn
\item[(c)]   if $\beta < \alpha$ \underline{then}  $N_\beta \cup
  \bar b^3_\beta \caret \bar b^4_\beta \subseteq N_\alpha$
\sn
\item[(d)]  $\bar b^4_\alpha \in {}^\zeta\! N$ is from $N^+$ satisfies: 
\sn
\begin{itemize}
\item  for every $\bar a \in
  {}^\eps(N_\alpha + \bar a^4_\alpha)$ 
we have $N \models \varphi[\bar a,\bar b^4_\alpha]$ iff 
$\{\beta < \chi:N \models \varphi[\bar a,\bar b^3_\beta]\} 
\in D_*$, equivalently
\sn
\item   $\bar b^4_\alpha$ realizes $\{\varphi(\bar a,
\bar y_{[\zeta]})^{\iif(\bft)}:\bar a \in {}^\eps(N_\alpha + 
\bar a^4_\alpha)$ and 
$\{\beta < \chi:N \models \varphi(\bar a,\bar b^3_\beta)^{\iif(\bft)}\} \in D_*$ and $\bft \in \{0,1\}\}$.
\end{itemize}
\end{enumerate}
\end{enumerate}
\mn
If we are stuck at $\alpha$ \underline{then}  obviously we can choose
$N_\alpha$ as required in clauses (a),(b),(c) of $(*)^2_\alpha$ 
hence there is no 
$\bar b^4_\alpha$ as required in $(*)^2_\alpha(d)$ hence $N$ is not
$(\chi^+,\theta,\{\psi\})$-saturated, (as otherwise $N_\alpha$ easily
exists).  Now as $N =$ l.u.p.$_{\bfx}(M)$ the desired
conclusion (B) holds for $M_1=M$.  So we can assume
that we succeed to carry the induction so $M_3 \defeq \cup \{N_\alpha:\alpha <
\chi\}$ is $\prec_{\bbL_{\theta,\theta}} N$.  
Now the pair $(M_3,\LL (\bar a^3_\alpha,\bar b^3_\alpha,\bar
b^4_\alpha):\alpha < \chi\RR)$, recalling that (by \ref{x32}) 
necessarily $\chi = \chi^{< \theta}$, satisfies
$\boxplus^\chi_{M_3,\LL(\bar a^3_\alpha,\bar b^3_\alpha,
\bar b^4_\alpha):\alpha < \chi\RR}$, where for a linear order $I$ 
and model $M_*$ we let
\medskip

$\boxplus^I_{M_*,\LL(\bar a^3_s,\bar b^3_s,
\bar b^4_s):s \in I\RR} \qquad (a) \quad M_*$ is a model of $T_1$
\smallskip

\hskip100pt  $(b) \quad \bar b^3_s,\bar b^4_s \in {}^\zeta(M_*)$ 
and $\bar a^3_s \in {}^\eps(M_*)$
\smallskip

\hskip100pt $(c) \quad$ if $\bar a \in {}^\eps(M_*)$ then
for some truth value

\hskip125pt  $\bft$ we have for every $s \in I$ large enough 

\hskip125pt  $M_* \models 
\varphi[\bar a,\bar b^3_s]^{\iif(\bft)}
  \wedge \varphi[\bar a,\bar b^4_s]^{\iif(\bft)}$
\smallskip

\hskip100pt  $(d) \quad M_* \models ``\varphi[\bar a^3_s,\bar b^4_t]"$  
for $s,t \in I$
\smallskip

\hskip100pt  $(e) \quad$ if $s,t < \chi$ then $M_* 
\models ``\varphi[\bar a^3_s,\bar b^3_t]"$ \underline{iff}  $s < t$.
\medskip

\noindent
Why?  For clause (c) let $\alpha < \chi$ be such that $\bar a \in
  {}^\eps(N_\alpha)$.  Now for all $\beta \in [\alpha,\chi)$
recall clause $(*)^2_\beta(d)$ and
  $(*)_1$.  For clause (d), by $\circledast_1(c)(\alpha)$ of
\ref{a26} we have $\alpha_1 < \beta_1 \Rightarrow N \models
  \varphi[\bar a^1_{\alpha_1},\bar b^1_{\beta_1}]$, hence by the choice of
$\LL \bar a^3_\gamma \caret \bar b^3_\gamma:\gamma < \chi\RR$
  we have $\gamma \in (\alpha,\chi) \Rightarrow N \models \varphi[\bar
  a^3_\alpha,\bar a^3_\gamma]$ so by $(*)^2_\alpha(d)$ we have $N
 \models \varphi[\bar a^3_\alpha,\bar b^4_\beta]$ as required in (d).

As for clause (e) by $\circledast_1(c)(\alpha)$ of \ref{a26} we have
$\beta,\alpha < \chi \Rightarrow N \models 
\varphi[\bar a^1_\alpha,\bar b^1_\beta]^{\iif(\alpha<\beta)}$ hence
by the choice of $\LL \bar a^3_\gamma \caret \bar
b^3_\gamma:\gamma < \chi\RR$ we have $\alpha,\beta < \chi
\Rightarrow N \models \varphi[\bar a^3_\alpha,\bar b^3_\beta]^{\iif(\alpha <
  \beta)}$.  So the pair $(M_3,\LL (\bar a^3_\alpha,\bar
b^3_\alpha,\bar b^4_\alpha):\alpha < \chi\RR$ is as promised.

As $|\tau_{T_1}| \le \theta$ by the case assumption, by the downward
LST theorem there are $M_4 
\prec_{\bbL_{\theta,\theta}} M_3$ of cardinality $\theta$ and an 
increasing sequence $\LL \alpha(i):\eps < \theta\RR$
 of ordinals $< \chi$ such that
$(M_4,\LL (b^3_{\alpha(\eps)},\bar a^3_{\alpha(\eps)},\bar
b^4_{\alpha(\eps)}):\eps < \theta \RR$ satisfies 
$\boxplus^\chi_{M_4,\LL (\bar a^3_{\alpha(\eps)},
\bar b^3_{\alpha(\eps)},\bar b^4_{\alpha(\eps)}):
\eps < \chi\RR}$.

Now it is easy to see that $\lup_{\bfx}(M_4)$ is not locally
$(\chi^+,\theta,\{\varphi\})$-saturated, a detailed proof is
included in the proof of Case 2.
\bigskip

\noindent
\underline{Case 2}:  $|T_1| > \theta$

Let $\tau_2 = \tau(T_1) \cup \{P,<,F_i,G_j,H_j):i < \eps,j <
\zeta\}$ where the union is disjoint, and $P,<$ are unary and binary
predicates respectively and $F_i,G_j,H_j$ are unary function symbols.

Let $T_2$ be the set of $\bbL_{\theta,\theta}(\tau_2)$-sentences such
that
\mn
\begin{enumerate}
\item[$(*)_3$]  for a $\tau_2$-model $M_2$ we have $M_2 \models T_2$ \underline{iff} 
\mn
\begin{enumerate}
\item[$(a)$]  $M_2 \models T_1$
\sn
\item[$(b)$]  $(P^{M_2},<^{M_2})$ is a linear order of cofinality $>
\partial$ for every $\partial < \theta$
\sn
\item[$(c)$]  $I = (P^{M_2},<^{M_2}),M'_3 = M_2 \rest
  \tau(T_1),\bar{\bfa} = \LL(\bar a^3_t,
\bar b^3_t,\bar b^4_t):t \in P^{M_2}\RR$
satisfies $\boxplus^I_{M_2,\bar{\bfa}}$ where we let
\sn
\begin{enumerate}
\item[${{}}$]  $\bullet \quad \bar a^3_t = \LL F^{M_2}_i(t):i <
  \eps\RR$ 
\sn
\item[${{}}$]  $\bullet \quad \bar b^3_t = \LL G^{M_2}_j(t):j <
  \zeta\RR$ 
\sn
\item[${{}}$]  $\bullet \quad \bar b^4_t = \LL H^{M_2}_j(t):j <
  \zeta\RR$.
\end{enumerate}
\end{enumerate}
\end{enumerate}
\mn
By Case 1 applied to $T_1 \cap \bbL_{\theta,\theta}(\tau')$ for any
$\tau' \subseteq \tau_T$ of cardinality $\le \theta$ such that
$\varphi(\bar x,\bar y) \in \bbL_{\theta,\theta}(\tau')$, hence clearly
$T_2$ is a theory.

By the proof of \ref{a27g}, for every $\lambda = \lambda^{< \theta} +
|T_1| \ge \kappa = \cf(\kappa) \ge \theta$, the theory $T_2$ has a
model $N = N_{\lambda,\kappa}$ of cardinality $\lambda$ such that
$\cf(P^N,<^N) = \kappa$, see \ref{a27k}(2), \ref{x31}(3).  
Applying this to the case
$\kappa = \theta$, consider $N^* = \lup_{\bfx}(N_{\lambda,\theta})$, so
$(P^{N^*},<^{N^*})$ has cofinality $\chi$, so let $\LL
t_\eps = t(\eps):\eps < \chi \RR$ be
increasing and cofinal in it and for $t \in P^{N_{\lambda,\theta}}$
let $\bar a^3_t = \LL F^{N_*}_i(t):i <
\eps\RR,\bar b^3_t = \LL G^{N_*}_j(t):j
< \zeta\RR,\bar b^4_t = \LL H^{N_*}_j(t:j <
\zeta)$, so the statement $\boxplus = \boxplus^\chi_{N_*,\bar{\bfa}_1}$ where $\bar{\bfa}_1 = \LL
(\bar a^3_{t(\xi)},\bar b^3_{t(\xi)},\bar b^4_{t(\xi)}):\xi <
  \chi\RR$ clearly holds.

Now for every $\bar a \in {}^\eps(N^*)$ by $(*)_3(c)$ and clause
(c) of $\boxplus$ clearly for some ordinal 
$\eps(\bar a) < \chi$ and truth value $\bft(\bar a)$ we have
\mn
\begin{enumerate}
\item[$(*)_5$]  if $\eps(\bar a) \le \xi < \chi$ then $N_* \models 
``\varphi[a,\bar b^3_{t(\xi)}]^{\iif(\bft(\bar a))} \wedge \varphi[\bar a,\bar b^4_{t(\xi)}]^{\iif(\bft(\bar a))}"$.
\end{enumerate}
\mn
For $\alpha \le \chi$ let $p_\alpha = \{\varphi(\bar x,\bar b^4_{t(\xi)}),\neg
\varphi(\bar x,\bar b^3_{t(\xi)}):\xi < \alpha\}$.  Now by $(*)_3(c)$ and
clauses (d),(e) of $\boxplus$ the sequence $\bar a^3_{t(\alpha)}$
realizes $p_\alpha$ in $N_*$ when $\alpha < \chi$
hence $p_\chi$, the increasing union of $\LL p_\alpha:\alpha < 
\chi\RR$ is $(< \chi)$-satisfiable in $N_*$.  However, by
$(*)_5$ no $\bar a \in {}^\eps(N_*)$ realizes $p_\chi$, 
so $p_\chi$ exemplifies $N_* = \uflp(M_4)$ is not 
$(\chi^+,\varphi(\bar x,\bar y))$-saturated so we have gotten
the desired conclusion.
\end{PROOF}

\begin{theorem}  
\label{a30}
Assume $T$ is a complete theory (in $\bbL_{\theta,\theta})$,
has $\theta$-n.c.p. and is definably stable and $\lambda = \lambda^{<\theta}$.

\noindent
1) $T$ is locally $\blacktriangleleft_{\lambda,\theta}$-minimal.

\noindent
2) If $D \in \ruf_{\lambda,\theta}(I)$
and $M \models T$ \underline{then}  $M^I/D$ is locally 
$(\lambda^+,\theta,\bbL_{\theta,\theta})$-saturated.
\end{theorem}

\begin{remark}
\label{a3k}
Note Theorem \ref{a30} deals with local 
$\blacktriangleleft_\lambda$-minimality, whereas \ref{a32} below deals with 
local $\blacktriangleleft^*_\lambda$-minimality and Claim \ref{a28} deals
with non-$\blacktriangleleft^*_\chi$-minimality.
\end{remark}

\begin{PROOF}{\ref{a30}}  
1) By part (2).

\noindent
2) Without loss of generality $|\tau_T| \le \theta$.

Let $\varphi(\bar x,\bar y) \in \bbL_{\theta,\theta}$ and $\partial =
\partial_\varphi < \theta$ witness $\varphi(\bar x,\bar y)$ fail the
$\theta$-c.p. and let $\eps = \lh(\bar x),\zeta = \ell
g(\bar y)$ and $N = M^I/D$, where $D \in \ruf_\theta(\lambda)$
and $M$ is a model of $T$ and $p(\bar x) = p_0(\bar x)$ is a positive
$\varphi$-type
in $N$ of cardinality $\le \lambda$, so $p(\bar x) \subseteq
\{\varphi(\bar x,\bar b):\bar b \in {}^{\lh(\bar y)}\!N\}$
 is $(< \theta)$-satisfiable in $N$.

As $\theta$ is a compact cardinal there is
$p_1(\bar x) \in \bfS^\eps_\varphi(N)$ extending $p_0(x)$.
By Definition \ref{y3} there are $\psi(\bar y,\bar z) \in
\bbL_{\theta,\theta}(\tau_T)$ and $\bar c \in {}^{\lh(\bar z)}\!N$
which define $p_1(\bar x)$.  Let $\bar c_s \in {}^{\lh(\bar z)}\!M$
for $s \in I$ be such that $\bar c = \LL \bar c_s:s \in
I\RR/D$ and for $s \in I$ let $\Gamma_s =
\{\varphi(x,\bar b)^{\iif(\bft)}:M \models
``\psi[\bar b,\bar c_s]^{\iif(\bft)}"$ and $\bft \in \{0,1\}\}$.

Let $I_\partial = \{s \in I:\Gamma_s$ is $(< \partial)$-satisfiable 
in $M_s$, that is if $\bar b_\alpha \in {}^\zeta(M_s)$ and $M_s \models 
\psi[\bar b_\alpha,\bar c_s]^{\iif(\bft(\alpha))}$ for $\alpha <
\partial$ then $M \models \exists \bar x \bigwedge\limits_{\alpha <
\partial} \varphi(\bar x,\bar b_\alpha)^{\iif(\bft(\alpha))}\}$; so by
\ref{x7} necessarily $I_\partial \in D$.

By the choice of $\partial$ and of $I_\partial$ for every $s \in
I_\partial$ the set $\Gamma^+_{\bfs} = \{\varphi(\bar
x,\bar b):M \models ``\psi[b,\bar c_s]"\}$ is $(< \theta)$-satisfiable
in $M_s$.

Let $\chi$ be large enough such that $M \in \clH(\chi)$ and let $\gB =
(\clH(\chi),\in,M)^I/D$.  As $s \in I \Rightarrow \Gamma^+_s \in
\clH(\chi)$ we have $\Gamma^+ \defeq \LL \Gamma^+_s:s \in I\RR/D 
\in \gB$ and $\gB \models ``\Gamma^+$ is a 
$(< \bfj(\theta))$-satisfiable over $M$" where 
$\bfj:\clH(\chi) \rightarrow
\gB$ is the canonical embedding.
Let $\Gamma' = \{\varphi(\bar x,\bar a):\gB \models ``\varphi(\bar
x,\bar a) \in \Gamma"\}$.
Hence to prove $p_0(\bar x)$ is realized it suffices to show
\mn
\begin{enumerate}
\item[$\bullet$]  there is $w \in \gB$ such that 
$\varphi(\bar x,\bar b)^{\iif(\bft)} \in p_0(x) 
\Rightarrow \gB \models ``\bar b \in w$ and $|w| < \bfj(\theta)"$.
\end{enumerate}
\mn
By \ref{x8}(2) this holds.
\end{PROOF}

\begin{theorem}
\label{a32}
Assume the complete $T \subseteq \bbL_{\theta,\theta}$ has
$\theta$-$\ncp$ and is 1-stable hence (by \ref{y5}) definably stable and $T_0
\supseteq T$ is a complete $\bbL_{\theta,\theta}$-theory.  \underline{Then}  for some
$\bbL_{\theta,\theta}$-theory $T_1 \supseteq T_0$ of cardinality $(|T|
+ \theta)^{< \theta}$, we have:
\mn
\begin{enumerate}
\item[$\bullet$]  if $M_1$ is a model of $T_1$, letting $\lambda$
be its cardinality, \underline{then}  $M' \rest \tau_T$ is locally
$(\lambda,\theta,\bbL_{\theta,\theta})$-saturated and $\lambda =
\lambda^{< \theta} \subseteq |T|$.
\end{enumerate}
\end{theorem}

\begin{remark}
\label{a32d}
Instead of ``$T$ is 1-stable" to prove $M_1$ is locally
$(\lambda,\theta,\Delta)$-saturated it is enough to assume
\mn
\begin{enumerate}
\item[$(a)$]  $\Delta \subseteq \bbL_{\theta,\theta}(\tau_T)$ has
cardinality $< \theta$
\sn
\item[$(b)$]  if $\varphi_1(\bar x,\bar y) \in \Delta$ then some
$\psi_{\varphi_1}(\bar y,\bar z)$ is as in the definition of definably
stable
\sn
\item[$(c)$]  $\Delta$ is closed under redividing the variables and
permuting variables
\sn
\item[$(d)$]  each $\varphi_1(\bar x,\bar y) \in \Delta$ is 1-stable
  in $T$.
\end{enumerate}
\end{remark}

\begin{PROOF}{\ref{a32}}
For any $\varphi(\bar x,\bar y) \in \bbL_{\theta,\theta}(\tau_T)$ let
$\psi_\varphi(\bar y,\bar z_\varphi)$ be as in Definition \ref{y3} of 
definably stable
for $\varphi$ and $T$, see Definition \ref{y3}(1) recalling $T$ is
definably stable by \ref{y5}(1).   For $\gamma < \theta$ 
let $\vartheta_{\varphi,\gamma}(\bar z_\varphi)$ be the formula saying 
that $(\forall \ldots \bar y_i 
\ldots)_{i < \gamma}(\bigwedge\limits_{i < \gamma} 
\psi(\bar y_i,\bar z) \rightarrow \exists x 
\bigwedge\limits_{i < \gamma} \varphi(\bar x,\bar y))$ and let
$\vartheta_\varphi(\bar z_\varphi) =
\vartheta_{\varphi,\partial_\varphi}(\bar z_\varphi)$.

Let $\Delta_\varphi \subseteq \{\varphi,\neg \varphi\}$
and let $\varphi^{[*]}(\bar x,\bar y_*)$ be as in 
\ref{a19}(3) for $\Delta$ and let $\theta_\varphi <
\theta$ be large enough and for $\Delta \subseteq
\bbL_{\theta,\theta}(\tau_T)$ be of cardinality $< \theta$, let 
$\theta_\Delta < \theta$ be large enough.

Now
\mn
\begin{enumerate}
\item[$(*)_1$]  let $T_2$ be the set of sentences in
$\bbL_{\theta,\theta}(\tau_2)$ where $\tau_2$ implicitly defined below
such that $M_2 \models T_2$ \underline{iff} :
\sn
\begin{enumerate}
\item[(a)]  $M_2 \models T_0$
\sn
\item[(b)]  $<^{M_2}$ is a well ordering of $|M_2|$ of cofinality
$\ge \theta$
\sn
\item[(c)]  if $\varphi = \varphi(\bar x,\bar y) \in
\bbL_{\theta,\theta}(\tau_T)$ and $\bar c \in \vartheta_\varphi(M_2)$
and $d \in M_2$ \underline{then}  $\bar a^{\varphi,M_2}_{\bar c,d} \defeq \LL
F_{\varphi,i}(d,\bar c):i < \lh(\bar z_\varphi)\RR$ realizes
$p^{\varphi,M_2}_{\bar c,d} \defeq \{\varphi(x,\bar b):\bar b \in
{}^\zeta(M_2)$ and $i < \lh(\bar b) \Rightarrow b_i < d$ and $M_2
\models \psi_\varphi[\bar b,\bar c]\}$
\sn
\item[(d)]  $P^{M_2}$ is a closed unbounded set of $d$-s such that: if
$\Delta \subseteq \bbL_{\theta,\theta}(\tau_{T_2})$ has 
cardinality $< \theta$ and $\partial = \partial_\Delta < \theta$ is
large enough and $\cf(\{d':d' <^{M_1} d\},<^{M_1}) 
\ge \theta_\Delta$ then $M^{<d}_2 \defeq M_2 \rest \{d':d' < d^{M_2}\} 
\prec_\Delta M_2$
\sn
\item[(e)]  $a \mapsto \LL G^{M_2}_\eps(a):\eps <
\zeta\RR$ is a function from $M_2$ onto ${}^\zeta(M_2)$ for each
$\zeta < \theta$.
\end{enumerate}
\end{enumerate}
\mn
Now
\begin{enumerate}
\item[$(*)_2$]  $T_2$ is a theory.
\end{enumerate}
\mn
[Why?  Choose $\chi = \chi^{< \theta} \ge |T_2|$, let $M_0 \models
T_0$ be a $(\chi^+,\{\varphi\})$-saturated model (or just a locally
$(\chi^+,\theta,\bbL_{\theta,\theta}(\tau_T))$-saturated model); 
exists by \ref{a30}
+ L.S.T.  Choose $\LL M^2_\alpha:\alpha < \chi^+\RR$ a
$\prec_{\bbL_{\theta,\theta}}$-increasing sequence of
$\prec_{\bbL_{\theta,\theta}}$-submodels of $M_0$, 
each of cardinality $\chi$ increasing fast enough, 
i.e. choose $M^2_\alpha$ by induction on $\alpha$.   The rest should be clear.]
\mn
\begin{enumerate}
\item[$(*)_3$]  let $\tau_3 = \tau_2 \cup \{Q,F\},Q$ a unary predicate,
  $F$ a unary function symbol and
$T_3 \subseteq \bbL_{\theta,\theta}(\tau_3)$ is a set of sentences
such that a $\tau_3$-model $M_3$ satisfies $T_3$ \underline{iff} :
\sn
\begin{enumerate}
\item[(a)]  $M_3 \models T_2$
\sn
\item[(b)]  $Q^{M_3} \subseteq P^{M_3}$ is $<^{M_3}$-unbounded
\sn
\item[(c)]  $F^{M_3}$ maps $Q^{M_3}$ onto $|M_3|$ hence $Q^{M_3}$
is of cardinality $\|M_3\|$
\sn
\item[(d)]  if $d \in M_3$ and $\bar c \in {}^{\lh(\bar z)}(M^{<d}_3)$  
\underline{then}  $\LL e \in M_3:e$ satisfies $M_3 \models ``d < e \wedge
Q(e)"\RR$ is 2-indiscernible (even $n$-indiscernible for every $n$)
over $\bar c$ in $M_3 \rest \tau_2$
\end{enumerate}
\sn
\item[$(*)_4$]  $T_3$ is a theory.
\end{enumerate}
\mn
[Why?  Easy, e.g. it is enough to consider
$(\Delta,2)$-indiscernibility and for this imitate the proof of
\ref{a27g}.]
\mn
\begin{enumerate}
\item[$(*)_5$]  assuming $\varphi = \varphi(\bar x,\bar y) \in
  \bbL_{\theta,\theta}(\tau_T)$ for some cardinal 
$\partial^1_\varphi < \theta$, if $M_3 
\models T_3,\bar c \in \vartheta_\varphi(M_3)$ and 
$\bar b \in {}^{\lh(\bar y)}(M_3)$ 
\underline{then}  for some $A = A^{\varphi,M_3}_{\bar c,\bar b}
\subseteq P^{M_3}_{\Delta_\varphi}$ of cardinality $<
  \partial^1_\Delta$ we have:
\sn
\begin{enumerate}
\item[$\bullet$]  if $d_1,d_2 \in P^M$ and
$(\forall d \in A)(d_1 \le d \equiv d_2 \le d)$ then $M_3 \models
``\varphi[\bar a^{M_3,\varphi}_{\bar c,d_1},\bar b] \equiv
\varphi[\bar a^{M_3,\varphi}_{\bar c,d_2},\bar b]"$.
\end{enumerate}
\end{enumerate}
\mn
[Why?  Straightforward because $T$ is definably stable and $<^{M_3}$
is a linear well ordering but we give details.  Let
$\partial^1_\varphi < \theta$ be large enough.

Suppose $M_3 \models T_3$ hence $(|M_3|,<^{M_3})$ is a well ordering.
Without loss of generality $|M_3|$ is an ordinal $\alpha_*$ and $<^{M_3}$ is the usual
order so $\cf(\alpha_*) \ge \theta$.  Suppose $\bar c \in
\vartheta_\varphi(M_3)$ and $\bar b \in {}^{\lh(\bar y)}(|M_3|)$
and we shall prove that there is $A = A^{\varphi,M_1}_{\bar c,\bar b}
\subseteq P^{M_2}_{\Delta_\varphi}$ as required.

Toward this we choose by induction on $n$ a set $A_n$ such that:
\mn
\begin{enumerate}
\item[$(*)_{5.1}$]  $(a) \quad A_n \subseteq P^{M_3}$
has cardinality $\le \partial^1_\varphi$
\sn
\item[${{}}$]  $(b) \quad m < n \Rightarrow A_m \subseteq A_n$ and
$A_0 = \{\min\{\alpha \in P^{M_3}:\bar b \subseteq M^{< \alpha}_3\}\}$
\sn
\item[${{}}$]  $(c) \quad$ if $\alpha \in A_n$ and $\cf(M^{<\alpha}_3
\cap P^{M_3}) \ge \theta_{\Delta_\varphi}$, \underline{then}  there are $\psi_*,
\bar c_\alpha$ such that

\hskip25pt   (letting $\psi_{\varphi[*]} = 
\psi(\bar y_{[*]},\bar z_*)$: we have
\sn
\begin{enumerate}
\item[${{}}$]  $(\alpha) \quad \bar c_\alpha \in {}^{\lh(\bar
z_*)}(M^{< \alpha}_3)$
\sn
\item[${{}}$]  $(\beta) \quad$ if $\bar a \in (M^{< \alpha}_3)$ \underline{then} 
 $M_3 \models \varphi[\bar a,\bar b]$ \underline{iff}  $M_3 \models
\psi_*[\bar a,\bar c_\alpha]$
\sn
\item[${{}}$]  $(\gamma) \quad \bar c_\alpha \subseteq M^{< \beta}_3$
  for some $\beta < \alpha$ which belongs to $A_{n+1}$
\end{enumerate}
\sn
\item[${{}}$]  $(d) \quad$ if $\alpha \in A_n$ and $\cf(M^{< \alpha}_1
\cap P^{M_3}_{\Delta_\varphi},<^{M_3}) < \theta_{\Delta_\varphi}$ then

\hskip25pt  $(A_{n+1} \cap M^{<\alpha}_3 \cap P^{M_3})$ is cofinal in 
$(P^{M_3},<^{M_3})$.
\end{enumerate}
\mn
Recall $(P^{M_3}_{\Delta_\varphi},<)$ is a well order of cofinality
$\ge \theta$.

Now let $A = \bigcup\limits_{n} A_n$ and we shall prove $\bullet$ of $(*)_5$; 
suppose $d_1,d_2 \in P^{M_3} \setminus A$ and
$(\forall d \in A)(d < d_1 \equiv d < d_2)$.  If $\bar b \subseteq
{}^{\min(d_1,d_2)>}(M_3)$ then $d_1,d_2$ are $<^{M_3}$-above the unique
member of $A_0$, hence clearly $M_3 \models ``\varphi[\bar a^{M_3}_{\bar
c,d_1},\bar b] \equiv \varphi[\bar a^{M_3}_{\bar c,d_2},\bar b]"$ as
required.

If not, let $d'' \in A \subseteq P^{M_3}$ be minimal such that $d_1 < d''$
(equivalently $d_2 < d''$).  Now $d''$ cannot be the first, a successor or of
cofinality $< \theta$ in $(P^{M_3},<^{M_3})$ hence $(M^{< d''}_3 \cap
P^{M_3})$ has cofinality $\ge \theta_{\Delta_\varphi}$ 
(see $(*)_{5.1}(d)$ and use $(*)_{5.1}(c)$).  Let $\alpha = d''$ and
$\beta = \sup(A \cap \alpha)$, by $(*)_{5.1}(c)(\gamma)$ we have $\bar
c_\alpha \subseteq M^{< \beta}_3$ so by $(*)_{5.1}(c)(\beta)$ again
$M_3 \models ``\varphi[\bar a^{M_3}_{\bar c,d_1},\bar b] \equiv
\varphi[\bar a^{M_3}_{\bar c,d_2},\bar b]"$.  So we are done proving $(*)_5$.] 
\mn
\begin{enumerate}
\item[$(*)_6$]  if $\varphi = \varphi(\bar x,\bar y) \in
  \bbL_{\theta,\theta}(\tau_T)$, for $\partial^2_\varphi < \theta$
  large enough,  if $M_3 \models T_3,
\bar c \in \vartheta_i(M_3),\bar b \in {}^{\lh(\bar y)}(M_3)$
\underline{then}  for some $B \subseteq Q^{M_3}$ of cardinality 
$< \partial^2_\varphi$ and for some truth value $\bft$ we have
\sn 
\item[${{}}$]  $\bullet \quad$ if $\alpha \in Q^{M_3} \setminus B$
then $M_3 \models ``\varphi[\bar a^{M_3}_{\bar c,d},\bar b]^{\iif(\bft)}"$.
\end{enumerate}
\mn
[Why?  As otherwise we get contradiction to $\varphi$ is 1-stable.  
In details, let $M_3,\bar b$ be a counterexample; let
$\partial_2 < \theta$ be large enough and $\kappa =
\cf(|M_3|,<^{M_3})$ let $\kappa \ge \theta$; and let $\LL d_i:i <
\kappa\RR$ be $<^{M_3}$-increasing cofinal and $d_i \in Q^{M_3}$.

Now $\bar b \in {}^\zeta(M_3)$ hence there is $d_* \in Q^M$ such that $\bar
b \subseteq M^{<d_*}_3$; so for some truth value, $d_* \le^{M_3} d
\Rightarrow M_3 \models ``\varphi[\bar a^{M_3}_{\bar c,d},\bar
b]^{\iif(\bft)}"$.

Let $A^{\varphi,M_3}_{M_3,\bar c,\bar b}$ be as in $(*)_5$ and $E = 
E_{M_3,\bar c,\bar b} = \{(d_1,d_2):d_1,d_2 \in Q^{M_3}$ and
$(\forall d \in A^{\varphi,M_3}_{M_3,\bar c,\bar b})(d < d_1
\equiv d < d_2 \wedge d = d_1 \equiv d = d_2)\}$ is an equivalence
relation and let $A^+_{M_3,\bar c,\bar b} = \{d \in
Q^M:d/E_{M_3,\bar c,\bar b}$ has $\le \partial_2$ members$\}$.  Now if $d
\in Q^{M_3} \setminus A^+_{M_3,\bar c,\bar b} 
\Rightarrow M_3 \models ``\varphi[\bar
a^{M_3}_{\bar c,d},\bar b]^{\iif(\bft)}"$, we are done, otherwise
let $d^*$ be a counterexample.  Let $d^*_1 = \min(d^*/E)$ and $d^*_2 
\in (A_{M_3,\bar c,\bar b} \setminus M^{<d^*}_3)$ and let $d^*_3 = d_*$.

Now $M_3$ satisfies
\mn
\begin{enumerate}
\item[$(*)_{6.1}$]
\begin{enumerate}
\item[(a)]  $M_3 \models ``d^*_1 < d^*_2 < d^*_3
  \wedge Q(d^*_1) \wedge Q(d^*_2) \wedge Q(d^*_3)$
\sn
\item[(b)]  for some $\bar b' \in {}^\eps(M_3)$ we have
$M_3 \models (\forall t) \in [d^*_1 < t < d^*_2 \wedge P(t)
  \rightarrow \varphi(\LL F_i(t),\bar c):
i < \eps,\bar b']^{\iif(\neg \bft)}]$ and 
$M^*_3 \models (\forall t)[d^*_3 < t \wedge P(t) \rightarrow 
\varphi(\LL F_i(t),\bar c):i < \eps\RR,
\bar b)^{\iif(\bft)}]$.
\end{enumerate}
\end{enumerate}
\mn
By the demand on $Q^{M_3}$
\mn
\begin{enumerate}
\item[$\bullet$]  for every $d'_1 < d'_2 < d'_3$ from $Q^{M_3}$ for
  some $\bar b' \in {}^\zeta(M_3)$ we have $M_3 \models 
(\forall t)[d'_1 < t < d'_2 \wedge P(t) \rightarrow
\varphi(\LL F_i(t,\bar c):i < \eps \RR,
\bar b')^{\iif(\neg \bft)}]$ and $M^*_3 \models (\forall t)[d'_3 
< t \wedge P(t) \rightarrow \varphi(\LL F_i(t,\bar c):
i < \eps),\bar b')^{\iif(\bft)}]$.
\end{enumerate}
\mn
\relax From this clearly $T$ has the order property, contradiction, so
$(*)_6$ holds indeed.]   Now the required
saturation follows.  That is, assume $\bar c \in
\vartheta(M_3),p_{\bar c} = \{\varphi(\bar x,\bar b):M \models
\psi[\bar b,\bar c]\}$, so a type of cardinality $\le \|M\|^{|\ell
  g(\bar x)|}$ but $\|M\| =
\|M\|^{< \theta}$ by \ref{x32}, and every $\varphi(\bar x,\bar b) \in
p_{\bar c}$ is realized by every $\bar a^{M_3}_{\bar c,d}$ for every
$d \in Q^{M_3}$ except possibly $\le \partial_2$ many.  As $|Q^M| = \|M\|$
by $(*)_5(c)$, we are done.
\end{PROOF}

\noindent
We can now sum up, giving full characterization of two versions of
local minimality.  Note that at last we state the main results
\ref{a33}, \ref{a33g}.
\begin{conclusion}  
\label{a33}
Assume $T$ is a complete $\bbL(\tau_T)$-theory.

Assume $\lambda = \lambda^{< \theta} \ge 2^\theta + |T|$, \underline{then}  $T$
is locally $(\lambda,\theta)$-minimal \underline{iff}  $T$ is 1-stable with
$\theta$-n.c.p.
\end{conclusion}

\begin{PROOF}{\ref{a33}}  

\noindent
\underline{Case 1}:  $T$ has the $\theta$-c.p.

Let $T_1 \supseteq T$.
Let $D_1 \in \ruf_\theta(\lambda)$ and $D_2$ be an e.g. normal ultrafilter
on $\theta$ and so $D = D_1 \times D_2 \in \ruf_\theta(\lambda \times
\theta)$.  If $M \models T_1$ then $M^{\lambda \times \theta}/D \cong
(M^\lambda/D_1)^\theta/D_2$; let $M_0 = M,M_1 = M^\lambda_0/D$ and $M_2 =
M^\theta_1/D$, all models of $T_1$. So $M^{\lambda \times \theta}/D$ is
isomorphic to $M^\theta_1/D$ and the latter is not locally
$((2^\theta)^+,\theta,\bbL_{\theta,\theta}(\tau_T))$-saturated by \ref{a25},
(hence not $(\lambda^+,\theta,\bbL_{\theta,\theta})$-saturated).
\medskip

\noindent
\underline{Case 2}:  $T$ is 1-unstable.

Let $T_1 \supseteq T$ and $M \models T_1$ and $M^+$ be a
$\theta$-complete expansion of $M$.

Now apply Claim \ref{a28} to the theory $T_1$ so for some $M_1 \models
T_1$, so for some $(\theta,\theta)$-$\luft$ $\bfx$ we have $\theta =
\cf(\lup_{\bfx}(\theta,<))$ (this exists by \ref{x31}(3)) hence the model 
$\lup_{\bfx}(M)$ is not locally
$(\theta^+,\theta,\bbL_{\theta,\theta}(\tau_T))$-saturated 
so we are done.
\medskip

\noindent
\underline{Case 3}:  $T$ is 1-stable with $\theta$-n.c.p.

Use Theorem \ref{a32}.
\end{PROOF}

\begin{conclusion}
\label{a33g}
Assume $\lambda = \lambda^{< \theta} \ge 2^\theta + |T|$ and $T$ is a
complete $\bbL_{\theta,\theta}(\tau)$-theory of cardinality $\le
\lambda$.  \underline{Then}  $T$ is $\triangleleft_{\lambda,\theta}$-minimal
\underline{iff}  $T$ is definably stable with the $\theta$-$\ncp$ \underline{iff}  $T$ is
1-stable with the $\theta$-n.c.p.
\end{conclusion}

\begin{PROOF}{\ref{a33g}}
The third and second clauses are equivalent by \ref{a22}(4).
The proof splits to cases and is similar to the proof of \ref{a33}.
\medskip

\noindent
\underline{Case 1}:  $T$ has the $\theta$-$\cp$

Exactly as in the proof of \ref{a33}.
\medskip

\noindent
\underline{Case 2}:  $T$ is definably unstable

By Claim \ref{y5}(1), $T$ is 1-unstable.
Again use \ref{a28} but now using $\bfx$ which is simply $D \in
\ruf_\theta(\lambda)$; true \ref{a28} say ``for some $M_1$" but recall
\ref{a8}. 
\medskip

\noindent
\underline{Case 3}: $T$ is definably stable with the $\theta$-$\ncp$

Use \ref{a30}.
\end{PROOF}

\begin{claim}  
\label{a36}
1) If the set $\spec(\varphi(\bar x,\bar y),T)$ includes every regular
 $\partial < \theta$ or just belongs to every normal ultrafilter on
 $\theta$ and $\lambda \ge \theta$ \underline{then}  $T$ 
is $\blacktriangleleft_{\lambda,\theta}$-maximal.

\noindent
1A) Moreover, if $\spec(\varphi(\bar x,\bar y),T)$ belongs to every
normal ultrafilter on $\theta$ and $\lambda \ge 2^\theta$ \underline{then}  for
every theory $T_0 \supseteq T$ of cardinality $\le \lambda$
for some $\bbL_{\theta,\theta}$-theory $T_1$ extending $T_0$ of
cardinality $\lambda$ for every
model $M_1$ of $T_1,M_1 \rest \tau_T$ is not locally
$\theta^+$-saturated; so $T$ is $\blacktriangleleft_{\lambda,\theta}$-maximal.

\noindent
1B) In (1A) we can replace ``$\lambda \ge 2^\theta$" by ``$\lambda \ge \theta$
and $\theta \setminus \spec(\varphi,T)$ is not in the 
$(\lambda,\theta)$-weakly
compact ideal on $\theta$ (see in the proof)".

\noindent
2) There is a model $M_* = (\theta,E^M),E^M$ an equivalence relation
such that $T = \Th_{\bbL_{\theta,\theta}}(M)$ satisfies 
$\spec(xEy,T) = \theta \cap \Card$ hence $T$ is 
$\blacktriangleleft_{\lambda,\theta}$-maximal for every $\lambda$ and even
$\triangleleft^*_{\lambda,\bar\mu,\theta}$-maximal. 

\noindent
3) Assume $\kappa$ is supercompact with the Laver diamond.
There is a sequence of models $\LL M_A:A \subseteq \theta
\RR$ such that:
\mn
\begin{enumerate}
\item[$(a)$]  $M_A = (\theta,E_A)$ for $A \subseteq \theta,E_A$ an
  equivalence relation on $\theta$ 
\newline
such that letting $T_A = \Th(M_A)$
  we have
\sn
\item[$(b)$]  for $\lambda = \lambda^{<\theta},T_A
\blacktriangleleft_{\lambda,\theta} T_B$ \underline{iff}  $A \subseteq B$ \underline{iff} 
  $T_A \trianglelefteq^*_{\lambda,\theta} T_B$
\end{enumerate}
\end{claim}

\begin{PROOF}{\ref{a36}}  
1) By \ref{a25}, because for $\theta$-complete which is
not $\theta^+$-complete\footnote{being $(\lambda,\theta)$-regular is a
  stronger condition} 
ultrafilter on a set $I$ recalling \ref{x8}(3) and
``$\prod\limits_{\alpha < \theta} \alpha/D$ has cardinality $\theta$"
 we know that $\theta \in \{\prod\limits_{s \in I} 
\theta_s/E:\theta_s \in \spec(\varphi(\bar x,\bar y))a\}$.

\noindent
1A) To make the rest of the proof be also a proof of part (1B), 
let $\bbB$ be the Boolean Algebra 
$\cP(\theta)$ and let $\cF = \{f:f \in {}^\theta \theta$
satisfies $f(\alpha) < 1 + \alpha\}$.  Also without loss of generality, $|T| \le \theta$.

Let $M_0$ be a model of $T_0$ such that letting $M = M_0 \rest \tau$
we have $\clH(\theta) \subseteq M,M \rest \clH(\theta)
\prec_{\bbL_{\theta,\theta}} M$.  Let $M_1$ be an expansion of $M$ by
$\le \lambda$ symbols including $P^{M_1} = \clH(\theta),
P^{M_1}_u = u$ for $u \in \bbB,F^M_f \rest \theta = f$ for $f \in \cF$ and
the relations $R_1 = (\in \rest \clH(\theta))$ and $R^{M_1}_2 =
\{(\beta,\partial) \caret \bar a_{\partial,\beta}:\partial \in
\spec(\varphi,T),\beta < \partial\}$, where $\{\varphi(\bar x,\bar
a_{\partial,\beta}):\beta < \partial\}$ exemplified $\partial \in
\spec(\varphi,T)$ in the model $M$.

Lastly, let $T_1 = \Th_{\bbL_{\theta,\theta}}(M_1) \cup \{P_\theta(c)
\wedge (\exists^{\ge \partial} y)(y \in c):\partial < \theta\}$
recalling $\theta \in \bbB$.  The
rest should be clear but we shall give details.

Let $M_2$ be a model of $T_1$, so $(P^{M_2}_\theta,\in^{M_2} 
\rest P^{M_2}_\theta)$
is a linear order which is a well ordering, so without loss of generality $P^{M_2}_\theta =
\alpha_*$ for some ordinal $\alpha_*$ and $\in^{M_2} \rest P^{M_2}$ is the
usual order and $c^{M_2} \in P^{M_2}_\theta = \alpha_*$ is necessarily
$\ge \theta$, so $\theta \in P^{M_2}_\theta$.

Let $D = \{u \in \bbB:M_2 \models P_u(\theta)\}$ so this is an
ultrafilter on the Boolean algebra $\bbB$ 
which is $\theta$-complete and normal (for $\cF$, i.e. $(\forall f \in
\cF)(\exists A \in D)$[$f \rest A$ is constant]).  By the
assumption of the claim, $u_* \defeq \spec(\varphi,T) \in D$, so $M_2
\models ``P_{u_*}(\theta)"$ and let $p_* = \{\varphi(\bar x,\bar
a):\LL \beta,\theta\RR \caret \bar a \in R^{M_2}_2$ for some
$\beta < \theta\}$.

Now
\mn
\begin{enumerate}
\item[$\bullet$]  $p_*(\bar x)$ is not realized in $M_2$, i.e. $M_2
\rest \tau_T$.
\end{enumerate}
\mn
[Why?  Because $M_1$ satisfies the sentence saying this even replacing
  $\theta$ by any member of $P_{\spec(\varphi,T)}$ and $M_2
\models T_2$.]
\mn
\begin{enumerate}
\item[$\bullet$]  if $\partial < \theta$ \underline{then}  every subset of
$p_*$ of cardinality $\le \partial$ is satisfiable in $M_2 \rest \tau_T$.
\end{enumerate}
\mn
[Why?  Similarly.]

\noindent
1B) The proof is as in (1A), but the demand
\mn
\begin{enumerate}
\item[$(*)$]  there is $\bbB \subseteq \cP(\theta)$ of cardinality
  $\lambda$, including $[\theta]^{< \theta}$ but we also have 
$\cF \subseteq \{f \in {}^\theta
\theta:(\forall \alpha < \theta)(f(\alpha) < 1 + \alpha)\}$ of
cardinality $\le \lambda$ satisfying $\alpha < \theta \wedge f \in \cF
\Rightarrow f^{-1}\{\alpha\} \in \bbB$ such that \underline{there is
no} uniform $\theta$-complete ultrafilter $D$ on $\bbB$ such that $f \in \cF
\Rightarrow (\exists \alpha)(f^{-1}\{\alpha\} \in D)$.
\end{enumerate}
\mn
In the proof ``the ultrafilter $D$ is normal for $\cF$" means $f \in
\cF \Rightarrow (\exists \alpha < \theta)(f^{-1}\{\alpha\} \in D)$.
By the way this implies 
$\theta$-complete when $\cF$ is the set of all regressive
$f \in {}^\theta \theta$.  Why?  If $A = \bigcup\limits_{i < \partial}
A_i$, let $f:\theta \rightarrow \theta$ be $f(\alpha)$ is 0 if $\alpha
< \partial$ and if $\min\{i < \partial:\alpha \in A_i\}$ if $\alpha
\ge \partial$.

\noindent
2) E.g. $E^M = \{(\alpha,\beta):\alpha + |\alpha| = \beta + |\beta|\}$
   satisfies the first demand; the first ``hence" follows by (1), the
   second hence by (1B).

\noindent
3) Let $C = \{\mu < \theta : \mu \text{ is strong limit}\}$, let $\LL
   S_i:i < \theta\RR$ be a partition of $C$ to $\theta$ unbounded
   subsets of $C$ such that for each $i$ there is a normal ultrafilter
$D^*_i$ on $\theta$ to which $S_i$ belongs; moreover, for every $\lambda
\ge \theta$ for some normal ultrafilter $D$ on $[\lambda]^{< \theta}$
the set $\{u \in [\lambda]^{< \theta} : u \cap \theta \in S_i\}$ belongs
to $D$.  Well known to exist,  see Kanamori-Magidor \cite{KnMg78}.
For $A \subseteq \theta$,
   let $E_A$ be an equivalence relation on $\theta$ such that
$\{|(\alpha/E_A|:\alpha < \theta\} = \bigcup\{S_i:i \in A\}$. So the
 following claim \ref{a38} will suffice.
\end{PROOF}

\begin{claim}  
\label{a38}
Assume $\theta < \lambda = \lambda^{< \theta}$ and $f_*:\theta \rightarrow
\theta$ satisfies $\alpha < \theta \Rightarrow \alpha < f_*(\alpha)
\in \Card$ and there is a transitive class
$\bfM \supseteq {}^\lambda \bfM$, a model of ZFC including the ordinals
 and an elementary embedding $\bfj$ of $\bfV$ into 
$\bfM$ with critical point $\theta$
such that $(\bfj(f_*))(\theta) = \lambda$.

Let $E$ be a thin enough club of $\theta,S_1 = \Rang(f_* \rest E)$ and
let $S_2 = \{2^\mu:\mu \in S_1\}$.

\underline{Then}  there is $D \in \ruf_\theta(\lambda)$ such that we have:
\mn
\begin{enumerate}
\item[$(a)$]  if $f:\lambda \rightarrow S_1$ \underline{then}  the cardinal
  $\prod\limits_{\alpha < \lambda} f(\alpha)/D$ is $< \theta$ or is $\ge
  \lambda$
\sn
\item[$(b)$]  for some $f:\lambda \rightarrow S_1$ we have
  $\prod\limits_{\alpha < \lambda} f(\alpha)/D$ is $\lambda$
\sn
\item[$(c)$]  if $f:\lambda \rightarrow S_2$ \underline{then}  the cardinality
  $\prod\limits_{\alpha < \lambda} f(\alpha)/D$ is $< \theta$ or is 
$\ge 2^\lambda$
\sn
\item[$(d)$]  for some $f:\lambda \rightarrow S_2$ we have
$\prod\limits_{\alpha < \lambda} f(\alpha)/D$ is $2^\lambda$.
\end{enumerate}
\end{claim}

\begin{PROOF}{\ref{a38}}  
Let $E = \{\mu <\theta:\mu$ strong limit 
and $\Rang(f_* \rest \mu) \subseteq \mu\}$, it is the club of $\theta$,
mentioned in the claim.  Let $S_1 = \{f_*(\mu):\mu \in E\}$ 
and $S_2 = \{2^{f_*(\mu)}:\mu \in S_1\}$.

Let $D$ be the following normal ultrafilter on $I =
 [\lambda]^{<\theta}$

\[
\{\cU \subseteq I:\{\bfj(\alpha):\alpha < \lambda\} \in \bfj(\cU)\}.
\]

\mn
Hence the following set belongs to $D:\{s \in I:s \cap \theta \in E$ and 
$|s| = f_*(s \cap \theta)\}$.  

Clearly $D$ is a $\theta$-complete $(\lambda,\theta)$-regular
ultrafilter on a set $I$, even normal and fine, and the set $I$ has 
cardinality $\lambda^{< \theta} =
\lambda$, so (by renaming)  can serve as $D$ in the claim.

Let $G_s:\cP(s) \rightarrow |\cP(s)|$ be one to one onto for each $s \in I$.

By the normality of $D$, in $(\theta,<)^I/D$, the 
$\theta^\tthh$ element is $f_0/D$ where $f_0:I \rightarrow \theta$ is
defined by $f_0(s) = \min(\theta \setminus s)$.

Now clause (b) holds for the function $f_* \circ f_0$, because 
 $\prod\limits_{s \in I} (f_* \circ f_0)(s),<)$ is 
isomorphic to $(\lambda,<)$ by the choice of $D$, hence 
$f_* \circ f_0/D$ is the 
$\lambda^\tthh$ member of $(\theta,<)^I/D$.  As for clause (a) 
if $g/D \in \theta^I/D,\Rang(g) \subseteq S_1$ and
$g <_D f_* \circ f_0$ then by the normality of 
$D,\prod\limits_{s} g(s)/D$ has cardinality $< \theta$.

Note that $f_* \circ f_0(s) = \min\{\gamma \in S_1:\gamma > \sup(s
\cap \theta)\}$.  

To prove clause (d) let $f_2 \in {}^I \theta$ be $f_2(s) = \min\{\gamma
\in S_2:\gamma > \sup(s \cap \theta)\}$, so $f_2(s) = 2^{f_*(s
  \cap \theta)}$ when $s \cap \theta \in E$ and easily
$\prod\limits_{s \in I} f(s)/D$ is of cardinality $\le \theta^I =
\theta^\lambda = 2^\lambda$.  In fact, it is of cardinality
$2^\lambda$ as exemplified by $\LL f_{\cU}/D:\cU \subseteq
\lambda\RR$ where for 
$\cU \subseteq \lambda$ let $f_{\cU}:I \rightarrow \theta$ be
$f_{\cU}(s) = G_s(\cU \cap s)$.  Also clause (c) follows, similarly to
the proof of clause (a).
\end{PROOF}
\newpage

\section {Global $\cp$ and full minimality} \label{4}

\begin{definition}
\label{a43}
1) Let $T \subseteq \bbL_{\theta,\theta}(\tau_T)$ be complete.  We say
   $T$ has the global $\theta$-c.p. (negation: global 
$\theta$- n.c.p.) \underline{when}  for some 
pair $(\bar\varphi,\bar\partial)$ it has the global
 $(\bar\varphi,\bar\partial)$-c.p., see below.

\noindent
2) $T$ has the global $(\bar\varphi,\bar\partial)$-c.p. \underline{when}  for
   some $S$ and $\eps$:
\mn
\begin{enumerate}
\item[$(a)$]  $S \subseteq \theta$ belongs to some normal ultrafilter
on $\theta$ and is a set of cardinals
\sn
\item[$(b)$]  $\eps < \theta$ and $\bar\varphi = \LL
\varphi_\alpha(\bar x_{[\eps]},\bar y_{\varphi_\alpha}):\alpha
< \theta\RR$ where $\varphi_\alpha \in \bbL_{\theta,\theta}(\tau_T)$
\sn
\item[$(c)$]  $\bar\partial = \LL \partial_\alpha:\alpha \in
S\RR$ and $\partial_\alpha$ is a cardinal $\in [\alpha,\theta)$
\sn
\item[$(d)$]  if $\alpha \in S$ then $\partial_\alpha \in
\spec(\bar\varphi \rest \alpha,T)$, see Definition \ref{a17}(3),(4).
\end{enumerate}
\end{definition}

\begin{observation}
\label{a44}
If $T$ has the $\theta$-$\cp$ \underline{then}  $T$ has the global $\cp$.
\end{observation}

\begin{claim}
\label{a46}
Assume $D$ is a normal ultrafilter on $\theta$ and $T$ has the
global $(\bar\varphi,\bar\partial)$-c.p., $S = \Dom(\bar\partial) \in D$ and
$M$ is a model of $T$ and $\chi = \theta^\theta/D$ or just $\chi = \Pi
\bar\partial/D$.

\noindent
1) $N = M^\theta/D$ is not fully
$(\chi^+,\theta,\bbL_{\theta,\theta})$-saturated.

\noindent
2) If $T_1 \supseteq T$ \underline{then}  for some model $M_1$ of $T_1$, the
model $(M_1 \rest \tau(T))^\theta/D$ is not fully
 $(\chi^+,\theta,\bbL_{\theta,\theta})$-saturated. 
\end{claim}

\begin{PROOF}{\ref{a46}}
1) Let $M \models T$ and for $i \in S$ let $\LL
\varphi_{\xi(i,j)}(\bar x_{[\eps]},\bar a_{i,j}):j <
\partial_i\RR$ witness $\partial_i \in \spec(\bar\varphi \rest
i,T)$ and $j < \partial_i \Rightarrow \xi(i,j) < i$.  
Let $\partial'_\eps$ be
$\partial_\eps$ if $\eps \in S$ and 1 if $\eps
\in \lambda \setminus S$.  We can fix $\bar f = \LL
f_\alpha:\alpha < \chi\RR$ such that $f_\alpha \in
\prod\limits_{\eps < \theta} \partial'_\eps$ and $\bar
f$ is a set of representatives for $\prod\limits_{i < \theta}
\partial'_i/D$.  For each $\alpha < \chi$, as $D$ is a normal
ultrafilter on $\theta$ to which $S$ belongs and 
$i \in S \Rightarrow \xi(i,f_\alpha(i)) <
i$ clearly for some $\zeta(\alpha) < \theta$ 
we have $S_\alpha \defeq \{i < \theta:i \in S$ and
$\xi(i,f_\alpha(i)) = \zeta(\alpha)\} \in D$ and let $\bar a^*_\alpha
\subseteq N$ be of length $\lh(\bar y_{\varphi_{\zeta(\alpha)}})$
such that $\bar a_\alpha = \LL \bar a_{i,f_\alpha(i)}:i \in
S_\alpha\RR/D$ and let $\Gamma = \{\varphi_{\zeta(\alpha)}(\bar
x_{[\eps]},\bar a_\alpha):\alpha < \chi\}$.

Of course,
\mn
\begin{enumerate}
\item[$(*)_0$]   $\Gamma$ has cardinality $\le \chi$
\sn
\item[$(*)_1$]   $\Gamma$ is a set of
$\bbL_{\theta,\theta}(\tau_T)$-formulas with parameters from $N$
\sn
\item[$(*)_2$]  $\Gamma$ is $(<\theta)$-satisfiable $M$.
\end{enumerate}
\mn
[Why?  Let $u \subseteq \chi$ have cardinality $< \theta$, hence
$\zeta(*) = \sup\{\zeta(\alpha):\alpha \in u\}$ is $< \theta$ and let
$S_* = \{i \in S$: if $\alpha \in u$ then $f_\alpha(i) =
\zeta(\alpha)$ and $|u| < i\}$.  Clearly $S_* \in D$ and if $i \in
S_*$ then $\{\varphi_{\zeta(\alpha)}(\bar x_{[\eps]},\bar
a_{i,f_\alpha(i)}):\alpha \in u\} \subseteq \{\varphi_{\xi(i,j)}(\bar
x_{[\eps]},\bar a_{i,j}):j <\partial_i\}$ and\footnote{The $\le
\partial_i$ is for technical reasons, anyhow $\partial_i =
|\partial_i+1|$.} has cardinality $< |i| < \partial_i$ hence is
realized in $M$, so $M \models (\exists \bar x_{[\eps]})
\bigwedge\limits_{\alpha \in u} \varphi_{\zeta(\alpha)}(\bar
x_{[\eps]},\bar a_{i,f_\alpha(i)})$.  Hence $N \models (\exists
\bar x_{[\eps]}) \bigwedge\limits_{\alpha \in u}
\varphi_{\zeta(\alpha)}(\bar x_{[\eps]},\bar a_\alpha)$ so we
are done.]
\mn
\begin{enumerate}
\item[$(*)_3$]  $\Gamma$ is not realized in $N$.
\end{enumerate}
\mn
[Why?  As in the proof of Case 2 of \ref{a28}, without loss of generality $\theta
\subseteq M$.  Let $\tau^* = \tau_T \cup \{P_\zeta,Q,<,R,F:\zeta <
\theta\}$ where $P_\zeta$ is a $(2 + \lh(\bar y_{\varphi_\zeta}))$-place
predicate, $Q$ is unary, $R$ is a $(1 + \eps)$ place predicate
and $F$ a unary function symbol.

For $i \in S$ let $M^+_i =
(M,Q^{M^+_i},P^{M^+_i}_\zeta,<^{M^+},R^{M^+_i},F^{M^+_i})_{\zeta <\theta}$ where
\mn
\begin{enumerate}
\item[$(*)_{3.1}$]  $\bullet \quad Q^{M^+_i} = \partial_i$
\sn
\item[${{}}$]  $\bullet \quad <^{M^+_i}$ the order on $\partial_i$
\sn
\item[${{}}$]  $\bullet \quad P^{M^+_i}_\zeta = \{\LL
\zeta,j\RR \caret \bar a_{i,j}:j < \partial_i$ and $\xi(i,j)=\zeta\}$
\sn
\item[${{}}$]  $\bullet \quad R^{M^+_i} = \{\LL j \RR \caret
\bar b:j < \partial_i$ and $\lh(\bar b) = \eps$ and $M
\models \varphi_{\xi(i,j)}[\bar b,\bar a_{i,j}]\}$
\sn
\item[${{}}$]  $\bullet \quad F^{M^+_i}(j) = \xi(i,j) < i$.
\end{enumerate}
\mn
Let $N^+ = \prod\limits_{i \in S} M^+_i/D$, so $N = N^+ \rest \tau_T$,
let $\bfi = \LL i:i \in S\RR/D \in N^+$ and $\partial =
\LL \partial_i:i \in S\RR/D \in N^+$
\mn
\begin{enumerate}
\item[$(*)_{3.2}$]  in $N^+$ there is no $\bar b \in
{}^\eps(N^+)$ such that for every $j \in Q^{N^+},N^+ \models
``j < \partial \rightarrow R[j,\bar b]"$
\sn
\item[$(*)_{3.3}$]  in $N^+$ if $j \in Q^{N^+}$ and $F^{N^+}(j) =
\zeta <\theta$ then $N^+ \models (\forall \bar
x_{[\eps]})(\forall \bar y)[P_\zeta(j,\zeta,\bar y) \rightarrow
R(j,\bar x_{[\eps]}) \equiv \varphi_\zeta(\bar
x_{[\eps]},\bar y)]$.
\end{enumerate}
\mn
Let
\mn
\begin{enumerate}
\item[$(*)_{3.4}$]   $\Gamma = \{\varphi_\zeta(\bar
x_{[\eps]},\bar a)$: for some $j \in Q^{N^+},\zeta =
F^{N^+}(j)$ we have $N^+ \models ``P_\zeta(j,\zeta,\bar a)"\}$.
\end{enumerate}
\mn
Together
\mn
\begin{enumerate}
\item[$(*)_{3.5}$]   $\Gamma$ is a set of
$\chi$ formulas from $\bbL_{\theta,\theta}(\tau_T)$ with parameters
from $N$ which is $(< \theta)$-satisfiable in $N$ but not 
realize in $N$ so we are done.
\end{enumerate}
\mn
2) Follows by (1).
\end{PROOF}

\begin{discussion}
\label{a49d}
Considering Theorem \ref{a33g}, \ref{a54} it is natural to wonder what
are the implications between ``$T$ has the $\theta$-$\ncp$" and ``$T$
has the global $\theta$-$\ncp$".

By \ref{a49} below the second does not imply the first and by
\ref{a44}, the first implies the second.
\end{discussion}

\begin{claim}
\label{a49}
There are a vocabulary $\tau,|\tau| \le \theta$ and a complete $T \subseteq
\bbL_{\theta,\theta}(\tau)$ which have $\theta$-n.c.p. but has the
global c.p.
\end{claim}

\begin{PROOF}{\ref{a49}}
For $i < \theta$ let 
$\partial_i$ be an infinite cardinal $\in [i,\theta)$.  Let $\tau
  = \{E,P_\zeta:\zeta < \theta\},E$ a two-place predicate, $P_\zeta$ a
  unary predicate.

We choose a $\tau$-model $M$ as follows:
\mn
\begin{enumerate}
\item[$(a)$]  its universe is $\theta \times \theta$
\sn
\item[$(b)$]  $E^M = \{((i,j_1),(i,j_2):i < \theta$ and 
$j_1,j_2 < \theta)\}$, an equivalence relation
\sn
\item[$(c)$]  $P^M_\zeta \subseteq |M|$ for $\zeta < \theta$
\sn
\item[$(d)$]  for $i < \theta$, letting $a_i = (i,0),A_i = a_i/E^M$, 
for every $\eta \in {}^i 2$ the following are equivalent:
\sn
\begin{enumerate}
\item[($\alpha$)]  there are $\theta$ elements 
$a \in A_i$ such that $(\forall \zeta < i)
(a \in P^M_\zeta \equiv \eta(\zeta)=1)$
\sn
\item[($\beta$)]  the set $\{a \in A_i$: if $\zeta < i$ then $a \in
  P^M_\zeta \equiv \eta(\zeta) =1\}$ has cardinality $\ne
  \partial_i$
\sn
\item[($\gamma$)]  the set $\{j < i:\eta(j)=1\}$ has cardinality $< 1 +|i|$.
\end{enumerate}
\end{enumerate}
\mn
We shall check that $T \defeq
\Th_{\bbL_{\theta,\theta}(\tau)}(M)$ is as required.

Let $A'_i \defeq \{a \in A_i$: if $\iota < i$ then $a \in P^M_\iota\}$; it 
is a subset of $A_i$ of cardinality exactly $\partial_i$ by clause
(d)($\alpha$) above
\mn
\begin{enumerate}
\item[$\boxplus_1$]  $T$ has global $\theta$-c.p.
\end{enumerate}
\mn
Why?  Let $\eps = 1,\bar y = \LL y_0,y_1\RR$ and
$\varphi_i = \varphi_i(x,\bar y) = x E y_0 \wedge P_i(x) \wedge x
\ne y_1$ for $i < \theta$ and let $\bar\varphi = \LL \varphi_i:i <
\theta\RR$.

For $i < \theta$ let $\Gamma_i = \{\varphi_j(x,\LL a_i,b\RR):b
\in A'_i$ and $j < i\}$
\mn
\begin{enumerate}
\item[$\bullet$]  $\Gamma_i$ is formally is as required for witnessing
  $\partial_i \in \spec(\bar\varphi \rest i,T)$ in particular
  $|\Gamma_i| = \partial_i$.
\end{enumerate}
\mn
[Why?  As $|A'_i| = \partial_i \ge i$.]
\mn
\begin{enumerate}
\item[$\bullet$]  $\Gamma_i$ is not realized.
\end{enumerate}
\mn
[Why?  As $\{x E a_i \wedge x \ne b \wedge P_\zeta(x):
b \in A'_i$ and $\zeta < i\}$ is not realized.]
\mn
\begin{enumerate}
\item[$\bullet$]  if $\Gamma \subseteq \Gamma_i$ has cardinality $<
  \partial_i$ then $\Gamma$ is realized.
\end{enumerate}
\mn
[Why?  As all but $< \partial_i$ members of $A'_i$ realize $\Gamma$.]

So $\boxplus_1$ holds indeed.
\mn
\begin{enumerate}
\item[$\boxplus_2$]  $T$ has the $\theta$-n.c.p.
\end{enumerate}
\mn
[Why?  Let $\varphi = \varphi(\bar x_{[\eps]},\bar y_{[\zeta]})$
and so for some $\kappa < \theta,\varphi$ belongs to
$\bbL_{\theta,\theta}(\{E,P_\zeta:\zeta < \kappa\})$, hence $M$
satisfies:
\mn
\begin{enumerate}
\item[$\bullet$]  if $a \in M,a \notin a_j/E^M$ for $j < \kappa^+$
  then for any $\eta \in {}^\kappa 2$ the set $\{b:b \in a/E^M$
  and $\zeta < \kappa \Rightarrow b \in P^M_\zeta \leftrightarrow
\eta(\zeta)=1\}$ has cardinality $\theta$.
\end{enumerate}
\mn
The rest should be clear.
\mn
\begin{enumerate}
\item[$\boxplus_3$]  $T$ is 1-stable.
\end{enumerate}
\mn
[Why?  Obvious.]

Together we are done.
\end{PROOF}

\begin{theorem}
\label{a52}
Assume $T$ is complete of cardinality $\theta$ and $T$ is definably stable
with global $\theta$-n.c.p. and $\lambda = \lambda^{< \theta}$.

\noindent
1) $T$ is $\triangleleft^{\ful}_{\lambda,\theta}$-minimal.

\noindent
2) Moreover, if $D \in \ruf_{\lambda,\theta}(I)$ and $\theta^I/D >
   \lambda$ and $M$ is a model of $T$  \underline{then}  $M^I/D$ is fully
   $(\lambda^+,\theta,\bbL_{\theta,\theta})$-saturated. 
\end{theorem}

\begin{PROOF}{\ref{a52}}
1) By part (2).

\noindent
2) As $T$ is definably stable we can use \ref{y15} and as $T$ has
$\theta$-$\ncp$ by \ref{a44}, we can use \ref{a17}, \ref{a19}.

Let $M \models T$ and $N = M^I/D$, let $\eps < \theta,
A \subseteq N,|A| \le \lambda$ and $p_0 \in 
\bfS^\eps(A,N)$ and we shall prove that $p_0(\bar
x_{[\eps]})$ is realized; by \ref{a8} and \ref{a30} without loss of generality, $M$ is
locally $(\lambda^+,\theta,\bbL_{\theta,\theta})$-saturated.
Let $\{\varphi(\bar x_{[\eps]},\bar y_{[\zeta]}):\varphi \in
\bbL_{\theta,\theta}(\tau_T)$ and $\zeta < \theta\}$ be listed
as $\LL \varphi_i(\bar x_{[\eps]},\bar y_{\zeta(i)}):i <
\theta\RR$.  Let $p_1(\bar x_{[\eps]}) \in \bfS^\eps(N)$ extends $p_0(\bar x_{[\eps]})$ and for each
$i < \theta$ let $\psi_i = \psi_i(\bar y_{\zeta(i)},\bar c_i)$ be a
formula from $\bbL_{\theta,\theta}(\tau_T)$ with parameters from $N$
defining $p_1(\bar x_{[\eps]}) \rest \varphi_i$ and let $\bar
c_\zeta = \LL \bar c_{\zeta,s}:s \in I\RR/D$.

As $D$ is a $(\lambda,\theta)$-regular ultrafilter, by \ref{x8}(2) there
is $\bar A = \LL A_s:s \in I\RR,A_s \in [M_s]^{< \theta}$
which is non-empty and $A = \{f_\alpha/D:\alpha < \lambda\}$ and $\alpha <
\lambda \Rightarrow f_\alpha \in \prod\limits_{s \in I} A_s$ 
and for $i \le \theta$ let $\Delta_i = \{\varphi_j(\bar x_{[\eps]},
\bar y_{\zeta(j)}):j < i\}$ and let 
$p_{s,i}(\bar x_{[\eps]}) = 
\{\varphi_j(\bar x_{[\eps]},\bar b):j<i,\bar b
\in A_s,M \models \psi_j(\bar b,\bar c_{j,s})\}$.

For each $i < \theta$ let $\partial_i = \sup(\spec(\Delta_i,T))$, see
\ref{a17}(3) so $\partial_i < \theta$ and let $I_i = \{s \in I$: there
is $p \in \bfS^\eps_{\Delta_i}(A_s)$ such that $\psi_j(\bar
y_{[\zeta(j)]},\bar c_{j,s})$ defines $p \rest \varphi_j$ for each
$j<i\}$.  

Now
\mn
\begin{enumerate}
\item[$(*)$]  $I_i \in D_i$.
\end{enumerate}
\mn
[Why?  Clear but we shall elaborate.  Clearly for every $\gamma <
  \theta$, letting $\bar y_{j,\gamma}$ be of length $\lh(\bar
  y_{\zeta(j)})$ the model $N$ satisfies $\vartheta_{i,\partial}(\ldots,\bar
c_j,\ldots)_{j<i}$ where

\begin{equation*}
\begin{array}{clcr}
\vartheta_{i,j} = \vartheta_{i,\partial}(\ldots,\bar
z^j,\ldots)_{j<i} \defeq &(\forall \ldots \bar y_{j,\gamma} \ldots)_{j<i,\gamma <
\partial} \Big[\bigwedge\limits_{j<i,\gamma < j} \psi_j(\bar
y_{j,\gamma},\bar z^j)^{\iif(\gamma\text{ is even})} \\
  &\Rightarrow (\exists x_{[\eps]}) \big(\bigwedge\limits_{j<i,\partial <j}
\varphi_i(\bar x_{[\eps]},
\bar y_{j,\gamma})^{\iif(\gamma\text{ is even})} \big) \Big].
\end{array}
\end{equation*}

\mn
Hence $I_i \supseteq \{s \in I:M \models
\vartheta_{i,\partial_i}(\ldots,\bar c_{j,s},\ldots)_{j<i}\}$ and so $I_i
\in D$.]  

Clearly $I_i \in D$ is decreasing with $i$.  Let 
$I'_\theta = \bigcap\{I_j : j < \theta\}$ and for $i < \theta$ let 
$I'_i = \bigcap\{I_j : j < i\} \setminus I_i$ for $i>0$ and let 
$I'_0 = I \setminus I_0$ and $\LL I'_i:i <
\theta\RR$ is a partition of $I \setminus I'_\theta$ 
to $\theta$ sets $= \varnothing \mod D$.

If $I'_\theta \in D$, recall that $M$ is
$(\lambda^+,\theta,\bbL_{\theta,\theta})$-saturated, hence we can find
$f \in {}^I\! M$ such that $s \in I'_\theta \Rightarrow f(s)$ realizes
$p_{s,\theta}$, clearly $f/D$ realizes $p$ in $N$ so we are done;
hence without loss of generality $I'_\theta = \varnothing$. 

Hence we can find $\bfh:I \rightarrow \theta$ such that $s \in I'_i
\Rightarrow \bfh(s) = i$.

Let $\bfh_* \in {}^I \theta$ be such that $\bfh_*/D$ is
the $\theta^\tthh$ member of $(\theta,<)^I/D$ and 
without loss of generality $\bfh_* \le \bfh$.
\bigskip

\noindent
\textbf{Case 1}:  $\bfh_* <_D \bfh$.

In this case we can prove that $p_0(\bar x_{[\eps]})$
is realized in $N$.
\bigskip

\noindent
\textbf{Case 2}:  Not Case 1.

In this case we can prove that $T$ has global $\theta$-c.p., contradicting an
assumption. 
\end{PROOF}

\begin{theorem}
\label{a53}
Assume $T$ is complete of cardinality $\theta$ and $T$ is 1-stable
with the global $\theta$-$\ncp$ and $\lambda = \lambda^{< \theta}$.
\underline{Then}  $T$ is $\blacktriangleleft^{*,\ful}_{\lambda,\theta}$-minimal.
\end{theorem}

\begin{question}
In the proof of \ref{a52} can we use ``$M$ is locally 
$(\lambda^+,\theta,\bbL_{\theta,\theta})$-saturated"?
We expect that we can prove this by combining the proofs 
of \ref{a52} and \ref{a32}.
\end{question}

\noindent
We now arrive to one of our main results.
\begin{conclusion}
\label{a54}
Assume $\lambda \ge 2^\theta,T$ is a complete
$\bbL_{\theta,\theta}(\tau_T)$-theory of cardinality $\theta$.  \underline{Then} 
$T$ is $\trianglelefteq^{\ful}_{\lambda,\theta}$-minimal \underline{iff} 
$T$ is definably stable and globally $\theta$-n.c.p.
\end{conclusion}

\begin{PROOF}{\ref{a54}}
Like the proof of \ref{a33g} by using \ref{a46}, \ref{a52} instead of 
\ref{a28} and \ref{a30} respectively.
\end{PROOF}

\begin{question}
\label{a58}
0) What are the implications between ``$T$ has $\theta$-$\ncp$" and
``$T$ has the global $\theta$-$\ncp$".  Debt.

\sn
1) For which $T$, for every $T_1 \supseteq T$, for every large enough
$\mu,\lambda = \lambda^\mu$ and $M_1 \ne T_2$ of cardinality $\lambda$,
there is a $(\mu^+,\theta,\bbL_{\theta,\theta})$-saturated $M_2$ of
cardinality $\lambda$ such that $M_1 \prec_{\bbL_{\theta,\theta}}
M_2$?

\sn
2) Can we fully characterize $(\lambda,\theta)$-minimal $T$ of
cardinality $\theta$?  We have to generalize superstable, say: 
every $p \in \bfS^\eps(M)$ is almost definable over some 
$A \in [M]^{< \theta}$, $\lambda = \lambda^{< \theta} \ge 2^\theta + |T|$, $T$ a complete $\bbL_{\theta,\theta}(\tau_T)$-theory.
\end{question}
\newpage


\bibliographystyle{amsalpha}
\bibliography{shlhetal}

\end{document}